\documentclass[12 pt]{amsart}
\usepackage{amsmath,amssymb,amsfonts,amsthm,enumerate,latexsym,floatflt,longtable}
\usepackage[dvips]{graphics}
\usepackage[centerlast, small]{caption2}

\title[The Table of the Structure Constants]{The Table of the Structure Constants for the Complex Simple Lie Algebra of Type $E_6$
and Chevalley Commutator Formulas in the Chevalley Group of Type $E_6$ over a Field}

\author[Anna I.~Polovinkina, Sergey G.~Kolesnikov]{Anna I.~Polovinkina, Sergey G.~Kolesnikov}

\date{\today}

\begin{document}

\sloppy



\maketitle
\begin{abstract}

This article is the third in the series. It is devoted the calculation of the structure constants for the complex simple Lie algebra of type $E_6$ and Chevalley commutator formulas.

\textit{Keywords:}  structure constant for the complex
simple Lie algebra, root system, Chevalley commutator formula.
\end{abstract}

\tableofcontents

\newpage

\section{General case}

\subsection{Information about the Root System of Type $E_6$}

In the Euclidean space $\mathbb{R}^8$ with the scalar product $(\,,)$ we choose
orthonormal basis  $e_2,$ $e_3,$ $e_4,$ $e_5,$ $e_6,$ $e_7,$ $e_8$. 
According to \cite{Bur72}, the vectors
$$
\pm e_i \pm e_j\ \ (1\leqslant i < j \leqslant 5),
$$
$$
\pm \frac12 \left(e_8 - e_7 - e_6 + \sum_{i=1}^{5} (-1)^{\nu(i)}e_i\right),
\text{ with even sum\ } \sum_{i=1}^{5} (-1)^{\nu(i),}
$$
from $\mathbb{R}^8$ form the root system of type $E_6$.
The roots
$$
a=(e_1 + e_8 - e_2 - e_3 - e_4 - e_5 - e_6 - e_7)/2, 
$$
$$
b = e_1 + e_2, \quad c = e_2 - e_1, \quad d = e_3 - e_2, \quad e = e_4 - e_3, \quad f = e_5 - e_4
$$
and respectively,
$$
\pm e_i + e_j\ (1\leqslant i < j \leqslant 5),
$$
$$
\frac12 \left( e_8 - e_7 - e_6 + \sum_{i=1}^{5} (-1)^{\nu(i)}e_i\right) 
\text{ with even sum\ } \sum_{i=1}^{5} (-1)^{\nu(i)} 
$$
form the system of fundamental $\Pi(E_6)$ and the system of positive roots $E_6^+$ in $E_6.$

Next, we will represent a positive root
$$
r=\alpha a + \beta b + \gamma c + \delta d + \epsilon e + \zeta  f
$$
as an ordered six of numbers $\alpha\beta\gamma\delta\epsilon\zeta.$ 
The sum $\alpha + \beta + \gamma + \delta + \epsilon + \zeta$
is called the height of the root $r$ and denoted by ${\rm ht}(r).$

Let $U$ be the subspace generated by the fundamental roots. We define the order $\prec$ on $U$ as follows:
for arbitrary vectors $x,y$ from $U$ we have
$$
x=\tau_1 f+\tau_2 e+\tau_3 b +\tau_4 d +\tau_5 c+\tau_6 a \prec y
 =\rho_1 f+\rho_2 e+\rho_3 b+\rho_4 d +\rho_5 c+\rho_6 a,
$$
if and only if the first coefficient in the expansion of the difference
$$
y-x=(\rho_1-\tau_1) f+(\rho_2-\tau_2) e+(\rho_3-\tau_3) b+(\rho_4-\tau_4) d +(\rho_5-\tau_5) c+(\rho_6-\tau_6) a 
$$
is positive. With respect to the order $\prec$, the positive roots of the root system $E_6$ are ordered in the following way:
$$
100000 \prec 001000 \prec 101000 \prec 000100 \prec 001100 \prec 101100 \prec 010000 \prec 010100 \prec 
$$
$$
011100 \prec 111100 \prec 000010 \prec 000110 \prec 001110 \prec 101110 \prec 010110 \prec 011110 \prec  
$$
$$
111110 \prec 011210 \prec 111210 \prec 112210 \prec 000001 \prec 000011 \prec 000111 \prec 001111 \prec  
$$
$$
101111 \prec 010111 \prec 011111 \prec 111111 \prec 011211 \prec 111211 \prec112211 \prec 011221 \prec 
$$
$$
111221 \prec 112221 \prec 112321 \prec 122321.
$$

The following table shows the correspondence between the number of a positive root in the ordered list and its coefficient of representation.
The table also shows root height. 

\begin{center}
\begin{tabular}{|r|l|l|c|r|l|l|c|}
    \hline
        N & $r$ & $\alpha\beta\gamma\delta\epsilon\zeta$ & ${\rm ht}(r)$ & N & $r$ & $\alpha\beta\gamma\delta\epsilon\zeta$ & ${\rm ht}(r)$ \\ \hline
        1 & $a$ & 100000 & 1  & 19 & $a+b+c+2d+e$ & 111210 & 6  \\ \hline
        2 & $c$ & 001000  & 1  & 20 & $a+b+2c+2d+e$ & 112210  & 7  \\ \hline
        3 & $a+c$ & 101000 & 2  & 21 & $f$ & 000001 & 1  \\ \hline
        4 & $d$ & 000100 & 1  & 22 & $e+f$ & 000011 & 2  \\ \hline
        5 & $c+d$ & 001100 & 2  & 23 & $d+e+f$ & 000111  & 3  \\ \hline
        6 & $a+c+d$ & 101100 & 3  & 24 & $c+d+e+f$ & 001111  & 4  \\ \hline
        7 & $b$ & 010000 & 1  & 25 & $a+c+d+e+f$ & 101111  & 5  \\ \hline
        8 & $b+d$ & 010100 & 2  & 26 & $b+d+e+f$ & 010111 & 4  \\ \hline
        9 & $b+c+d$ & 011100 & 3  & 27 & $b+c+d+e+f$ & 011111  & 5  \\ \hline
        10 & $a+b+c+d$ & 111100  & 4  & 28 & $a+b+c+d+e+f$ & 111111  & 6  \\ \hline
        11 & $e$ & 000010 & 1  & 29 & $b+c+2d+e+f$ & 011211  & 6  \\ \hline
        12 & $d+e$ & 000110 & 2  & 30 & $a+b+c+2d+e+f$ & 111211  & 7  \\ \hline
        13 & $c+d+e$ & 001110 & 3  & 31 & $a+b+2c+2d+e+f$ & 112211  & 8  \\ \hline
        14 & $a+c+d+e$ & 101110 & 4  & 32 & $b+c+2d+2e+f$ & 011221  & 7  \\ \hline
        15 & $b+d+e$ & 010110 & 3  & 33 & $a+b+c+2d+2e+f$ & 111221 & 8  \\ \hline
        16 & $b+c+d+e$ & 011110 & 4  & 34 & $a+b+2c+2d+2e+f$ & 112221 & 9  \\ \hline
        17 & $a+b+c+d+e$ & 111110 & 5  & 35 & $a+b+2c+3d+2e+f$ & 112321 & 10  \\ \hline
        18 & $b+c+2d+e$ & 011210 & 5  & 36 & $a+2b+2c+3d+2e+f$ & 122321 & 11  \\ \hline
    \end{tabular}
\end{center}
\vskip1mm
\centerline{Table 1.}
\vskip5mm

\subsection{Calculation of Structure Constants} In this section we prove the following

\textbf{Theorem 1.}
\textit{Let the structure constants $N_{r,s},$ corresponding to extraspecial pairs of roots be chosen as indicated in Table~2.}
\medskip

\begin{center}
\begin{tabular}{ |l |l | l| l | l| l | l |l| }
\hline
$N_{a,c}$         & $\alpha_1$ & $N_{a,b+c+2d+e+f}$  & $\alpha_9$ & $N_{c,b+d+e+f}$       & $\beta_7$  & $N_{d,a+b+2c+2d+2e+f}$ & $\gamma_6$   \\ 
$N_{a,c+d}$       & $\alpha_2$ & $N_{a,b+c+2d+2e+f}$ & $\alpha_0$ & $N_{c,a+b+c+2d+e+f}$  & $\beta_8$  & $N_{b,d+e}$            & $\delta_1$   \\ 
$N_{a,b+c+d}$     & $\alpha_3$ & $N_{c,d}$           & $\beta_1$  & $N_{c,a+b+c+2d+2e+f}$ &  $\beta_9$ & $N_{b,d+e+f}$          & $\delta_2$   \\ 
$N_{a,c+d+e}$     & $\alpha_4$ & $N_{c,b+d}$         & $\beta_2$  & $N_{d,b}$             & $\gamma_1$ & $N_{b,a+b+2c+3d+2e+f}$ & $\delta_3$   \\ 
$N_{a,b+c+d+e}$   & $\alpha_5$ & $N_{c,d+e}$         & $\beta_3$  & $N_{d,e}$             & $\gamma_2$ & $N_{e,f}$              & $\epsilon_1$ \\ 
$N_{a,b+c+2d+e}$  & $\alpha_6$ & $N_{c,b+d+e}$       & $\beta_4$  & $N_{d,b+c+d+e}$       & $\gamma_3$ & $N_{e,b+c+2d+e+f}$     & $\epsilon_2$ \\ 
$N_{a,c+d+e+f}$   & $\alpha_7$ & $N_{c,a+b+c+2d+e}$  & $\beta_5$  & $N_{d,e+f}$           & $\gamma_4$ & ~                      & ~            \\ 
$N_{a,b+c+d+e+f}$ & $\alpha_8$ & $N_{c,d+e+f}$       & $\beta_6$  & $N_{d,b+c+d+e+f}$     & $\gamma_5$ & ~                      & ~            \\ 
\hline
\end{tabular}
\end{center}
\vskip1mm
\centerline{Table 2.}
\vskip5mm

\textit{Then the values of non-zero structure constants $N_{r,s}$ and $N_{r,-s},$ where $0\prec r\prec s,$ will be as indicated in Tables~3 and 4 below.}
\vskip5mm

\noindent
\begin{tabular}{|c|c|c|c|c|c|}
\hline
 \boldmath$r$     & \boldmath$N_{a+c,r}$               & N   & \boldmath$r$    & \boldmath$N_{a+b+c+d,r}$                                  & N   \\[1mm]\hline
 $d$              & $\alpha_{12}\beta_1$               & 004 & $d+e$           & $-\alpha_{36}\beta_{24}\gamma_{13}\delta_1$               & 049 \\[1mm]\hline
 $b+d$            & $\alpha_{13}\beta_2$               & 011 & $c+d+e$         & $-\alpha_{36}\beta_{2345}\gamma_{13}\delta_1$             & 064 \\[1mm]\hline
 $d+e$            & $\alpha_{14}\beta_3$               & 022 & $e+f$           & $-\alpha_{38}\beta_{27}\gamma_{14}\delta_2$               & 097 \\[1mm]\hline
 $b+d+e$          & $\alpha_{15}\beta_4$               & 035 & $d+e+f$         & $-\alpha_{39}\beta_{27}\gamma_{15}\delta_2$               & 117 \\[1mm]\hline
 $b+c+2d+e$       & $-\alpha_{16}\beta_5$              & 058 & $c+d+e+f$       & $-\alpha_{39}\beta_{2678}\gamma_{15}\delta_2$             & 134 \\[1mm]\hline
 $d+e+f$          & $\alpha_{17}\beta_6$               & 080 & $b+c+2d+2e+f$   & $-\alpha_{30}\beta_{29}\gamma_{16}\delta_3$               & 203 \\[1mm]\hline
 $b+d+e+f$        & $\alpha_{18}\beta_7$               & 101 & $r$             & \boldmath$N_{e,r}$                                        & N   \\
\hline
\end{tabular}

\vskip5mm
\centerline{Table 3, part 1 of 2}

\noindent
\begin{tabular}{|c|c|c|c|c|c|}
\hline
$a+c+d+e$        & $\alpha_{45}\beta_{34}\delta_1$              & 037 & $\begin{array}{c} a+b+c+2d+\\+e+f\end{array}$ & $-\alpha_{90}\beta_{49}\gamma_{26}\delta_{13}\epsilon_2$    & 201 \\ [1mm]\hline
$c+d+e+f$        & $\beta_{67}\delta_2$                         & 090 & \boldmath$r$   & \boldmath$N_{a+b+c+d+e,r}$                                      & N   \\ [1mm]\hline 
$a+c+d+e+f$      & $\alpha_{78}\beta_{67}\delta_2$              & 099 & $f$            & $\alpha_{58}\beta_{47}\gamma_{24}\delta_{12}\epsilon_1$         & 095 \\ [1mm]\hline
\boldmath$r$     & \boldmath$N_{b+d,r}$                         & N   & $d+e+f$        & $-\alpha_{50}\beta_{47}\gamma_{25}\delta_{12}\epsilon_2$        & 154 \\ [1mm]\hline
$e$              & $-\gamma_{12}\delta_1$                       & 025 & $c+d+e+f$      & $-\alpha_{50}\beta_{4679}\gamma_{25}\delta_{12}\epsilon_2$      & 169 \\ [1mm]\hline
$c+d+e$          & $\beta_{34}\gamma_{13}\delta_1$              & 042 & $b+c+2d+e+f$   & $\alpha_{50}\beta_{49}\gamma_{26}\delta_{13}\epsilon_2$         & 207 \\ [1mm]\hline
$a+c+d+e$        & $\alpha_{46}\beta_{34}\gamma_{13}\delta_1$   & 051 & \boldmath$r$   & \boldmath$N_{b+c+2d+e,r}$                                       & N   \\ [1mm]\hline
$e+f$            & $-\gamma_{14}\delta_2$                       & 085 & $f$            & $\beta_{47}\gamma_{2345}\delta_{12}\epsilon_1$                  & 104 \\
\hline
\end{tabular}

\vskip3mm
\centerline{Table 3, part 2 of 2}

\newpage

\noindent
\begin{tabular}{|c|c|c|c|c|c|}
\hline
\boldmath$r$              &  \boldmath$N_{a,r}$    & N    & \boldmath$r$             & \boldmath$N_{a+c+d+e,r}$                                      & N    \\[1mm] \hline
$-a-c$                    & $-\alpha_1$            & 001  & $\begin{array}{c}-a-b-2c-2d-\\-2e-f\end{array}$ & $-\alpha_{40}\beta_{39}\gamma_{25}\epsilon_2$                 & 172  \\[1mm]\hline
$-a-c-d$                  & $-\alpha_2$            & 003  & $\begin{array}{c}-a-b-2c-3d-\\-2e-f\end{array}$ & $-\alpha_{40}\beta_{39}\gamma_{26}\epsilon_2$                 & 189  \\[1mm]\hline
$-a-b-c-d$                & $-\alpha_3$            & 010  & \boldmath$r$             & \boldmath$N_{b+d+e,r}$                                        & N    \\[1mm] \hline
$-a-c-d-e$                & $-\alpha_4$            & 019  & $-b-c-d-e$               & $\beta_4$                                                     & 027  \\[1mm] \hline
$-a-b-c-d-e$              & $-\alpha_5$            & 032  & $-a-b-c-d-e$             & $\alpha_{15}\beta_4$                                          & 036  \\[1mm]\hline
$-a-b-c-2d-e$             & $-\alpha_6$            & 046  & $-b-c-2d-e$              & $-\beta_{14}\gamma_3$                                         & 045  \\[1mm]\hline
$-a-c-d-e-f$              & $-\alpha_7$            & 075  & $-a-b-c-2d-e$            & $-\alpha_{26}\beta_{14}\gamma_3$                              & 054  \\
\hline
\end{tabular}
\vskip8pt

\centerline{Table 4, part 1 of 6}

\newpage

\noindent
\begin{tabular}{|c|c|c|c|c|c|}
\hline
$-a-b-c-d-e$      & $-\alpha_{15}\beta_4$            & 036  & $\begin{array}{c}-a-2b-2c-3d-\\-2e-f\end{array}$ & $-\alpha_{50}\beta_{49}\gamma_{26}\delta_{13}\epsilon_2$      & 208  \\[1mm]\hline
$\begin{array}{c}-a-b-2c-\\-2d-e\end{array}$ & $\alpha_{16}\beta_5$             & 059  & \boldmath$r$          & \boldmath$N_{b+c+2d+e,r}$                                           & N    \\[1mm] \hline
$-a-c-d-e-f$      & $-\alpha_{17}\beta_6$            & 081 & $-a-b-c-2d-e$      & $\alpha_6$                                                          & 046 \\[1mm] \hline
$\begin{array}{c}-a-b-c-d-\\-e-f\end{array}$ & $-\alpha_{18}\beta_7$   & 102 & $-a-b-2c-2d-e$      & $-\alpha_{16}\beta_5$                                               & 059 \\[1mm]\hline
$\begin{array}{c}-a-b-2c-2d-\\-e-f\end{array}$ & $\alpha_{19}\beta_8$  & 129 & $-b-c-2d-e-f$       & $-\beta_{47}\gamma_{2345}\delta_{12}\epsilon_1$                     & 105 \\[1mm]\hline
$\begin{array}{c}-a-b-2c-2d-\\-2e-f\end{array}$ & $\alpha_{10}\beta_9$ & 162 & $-b-c-2d-2e-f$      & $-\beta_{47}\gamma_{2345}\delta_{12}\epsilon_2$                     & 140 \\[1mm]\hline
\boldmath$r$      & \boldmath$N_{d,r}$               & N     & $\begin{array}{c}-a-b-2c-3d-\\-2e-f\end{array}$ & $-\dfrac{\alpha_{70}\beta_{4679}}{\gamma_{2356}\delta_{12}\epsilon_2}$ & 191 \\[3mm]\hline
$-c-d$            & $\beta_1$                        & 002 & $\begin{array}{c}-a-2b-2c-3d-\\-2e-f\end{array}$ & $-\dfrac{\alpha_{80}\beta_{49}}{\gamma_{2356}\delta_{13}\epsilon_2}$    & 206 \\
\hline
\end{tabular}
\vskip8pt

\centerline{Table 4, part 2 of 6}

\newpage

\noindent
\begin{tabular}{|c|c|c|c|c|c|}
\hline
$-b-c-2d-e-f$      & $\beta_{17}\gamma_5$             & 109 & $-b-d-e-f$          & $\gamma_{24}\delta_{12}\epsilon_1$                                                                 & 084 \\[1mm]\hline
$\begin{array}{c}-a-b-2c-2d-\\-e-f\end{array}$  & $-\alpha_{89}\beta_{18}\gamma_5$ & 127 & $-b-c-d-e-f$  & $\beta_{47}\gamma_{24}\delta_{12}\epsilon_1$                                & 089 \\[1mm]\hline
$\begin{array}{c}-a-b-2c-3d-\\-2e-f\end{array}$ & $\beta_{19}\gamma_6$             & 179 & $-a-b-c-d-e-f$      & $\dfrac{\alpha_{58}\beta_{47}}{\gamma_{24}\delta_{12}\epsilon_1}$      & 096 \\[3mm]\hline
\boldmath$r$       & \boldmath$N_{a+c+d,r}$           & N   & $-b-c-2d-e-f$       & $\beta_{47}\gamma_{2345}\delta_{12}\epsilon_1$                                                     & 105 \\[1mm]\hline
$-a-b-c-d$         & $-\alpha_{23}\beta_{12}\gamma_1$ & 014 & $\begin{array}{c}-a-b-c-2d-\\-e-f\end{array}$    & $\dfrac{\alpha_{69}\beta_{47}}{\gamma_{2345}\delta_{12}\epsilon_1}$    & 122 \\[3mm]\hline
$-a-c-d-e$         & $-\alpha_{24}\beta_{13}\gamma_2$                      & 021 & $\begin{array}{c}-a-b-2c-2d-\\-e-f\end{array}$  & $\dfrac{\alpha_{69}\beta_{4578}}{\gamma_{2345}\delta_{12}\epsilon_1}$    & 125 \\[3mm]\hline
$-a-b-c-2d-e$      & $\alpha_{26}\beta_{14}\gamma_3$                       & 054 & \boldmath$r$      & \boldmath$N_{e+f,r}$                                                     & N   \\[1mm]\hline
$-a-b-2c-2d-e$     & $\alpha_{26}\beta_{15}\gamma_3$                       & 061 & $-d-e-f$          & $\gamma_4$                                                               & 067 \\
\hline
\end{tabular}
\vskip8pt

\centerline{Table 4, part 3 of 6}

\newpage

\noindent
\begin{tabular}{|c|c|c|c|c|c|}
\hline
$-a-b-c-d$         & $\alpha_{13}\beta_2$                                  & 012 & \boldmath$r$      & \boldmath$N_{c+d+e+f,r}$                                        & N   \\[1mm]\hline
$-b-d-e$           & $\gamma_{12}\delta_1$                                 & 026 & $-a-c-d-e-f$      & $\alpha_7$                                                      & 075 \\[1mm]\hline
$-b-c-2d-e$        & $-\beta_{34}\gamma_{13}\delta_1$                      & 043 & $-b-c-d-e-f$      & $\beta_{67}\delta_2$                                            & 091 \\[1mm]\hline
$\begin{array}{c}-a-b-c-\\-2d-e\end{array}$ & $-\dfrac{\alpha_{46}\beta_{34}}{\gamma_{13}\delta_1}$   & 052 & $-b-c-2d-e-f$      & $\beta_{67}\gamma_{15}\delta_2$     & 107 \\[3mm]\hline
$-b-d-e-f$         & $\gamma_{14}\delta_2$                                  & 086 & $-a-b-2c-2d-e-f$  & $-\dfrac{\alpha_{39}\beta_{2678}}{\gamma_{15}\delta_2}$        & 135 \\[1mm]\hline
$\begin{array}{c}-b-c-2d-\\-e-f\end{array}$      & $-\beta_{67}\gamma_{15}\delta_2$  & 107 & $-b-c-2d-2e-f$    & $\beta_{67}\gamma_{25}\delta_{12}\epsilon_2$          & 144 \\[1mm]\hline
$\begin{array}{c}-a-b-c-\\-2d-e-f\end{array}$ & $-\dfrac{\alpha_{79}\beta_{67}}{\gamma_{15}\delta_2}$ & 116 & $-a-b-2c-2d-2e-f$ & $-\dfrac{\alpha_{50}\beta_{4679}\gamma_{25}}{\delta_{12}\epsilon_2}$               & 170 \\[3mm]\hline
$\begin{array}{c}-a-2b-2c-\\-3d-2e-f\end{array}$ & $\gamma_{16}\delta_3$     & 194 & $-a-b-2c-3d-2e-f$  & $-\dfrac{\alpha_{60}\beta_{4679}}{\gamma_{2356}\delta_{12}\epsilon_2}$ & 187 \\
\hline
\end{tabular}
\vskip8pt

\centerline{Table 4, part 4 of 6}

\newpage

\noindent
\begin{tabular}{|c|c|c|c|c|c|}
\hline
$\begin{array}{c}-a-b-2c-2d-\\-e-f\end{array}$   & $\dfrac{\alpha_{39}\beta_{2678}}{\gamma_{15}\delta_2}$  & 135 & \boldmath$r$ & \boldmath$N_{b+c+d+e+f,r}$                 & N   \\[3mm]\hline 
$\begin{array}{c}-a-2b-2c-3d-\\-2e-f\end{array}$ & $\dfrac{\alpha_{30}\beta_{29}}{\gamma_{16}\delta_3}$    & 204 & $\begin{array}{c}-a-b-c-d-\\-e-f\end{array}$ & $\alpha_8$ & 094 \\[3mm]\hline 
\boldmath$r$       & \boldmath$N_{e,r}$                  & N   & $-b-c-2d-e-f$                                    & $\gamma_5$                                               & 103 \\[3mm]\hline
$-d-e$             & $\gamma_2$                          & 015 & $\begin{array}{c}-a-b-2c-2d-\\-e-f\end{array}$   & $-\alpha_{29}\beta_{18}\gamma_5$                         & 131 \\[3mm]\hline
$-c-d-e$           & $\beta_{13}\gamma_2$                & 018 & $-b-c-2d-2e-f$                                   & $-\gamma_{25}\epsilon_2$                                 & 138 \\[3mm]\hline
$-a-c-d-e$         & $\dfrac{\alpha_{24}}{\beta_{13}\gamma_2}$  & 021 & $\begin{array}{c}-a-b-2c-2d-\\-2e-f\end{array}$  & $\alpha_{40}\beta_{39}\gamma_{25}\epsilon_2$      & 172 \\[3mm]\hline
$-b-d-e$           & $-\gamma_{12}\delta_1$              & 026 & $\begin{array}{c}-a-2b-2c-\\-3d-2e-f\end{array}$ & $-\dfrac{\alpha_{60}\beta_{49}}{\gamma_{2356}\delta_{13}\epsilon_2}$ & 210 \\[3mm]\hline
$-b-c-d-e$         & $-\dfrac{\beta_{24}}{\gamma_{12}\delta_1}$ & 029 & \boldmath$r$                              & \boldmath$N_{a+b+c+d+e+f,r}$                             & N   \\
\hline
\end{tabular}
\vskip3mm

\centerline{Table 4, part 5 of 6}

\newpage

\noindent
\begin{tabular}{|c|c|c|c|c|c|}
\hline
$\begin{array}{c}-a-b-c-\\-2d-2e-f\end{array}$   & $\alpha_{80}\gamma_{25}\epsilon_2$            & 151 & $\begin{array}{c}-a-b-c-\\-2d--2e-f\end{array}$     & $\alpha_{90}\epsilon_2$                                             & 149 \\[3mm]\hline
$\begin{array}{c}-a-b-2c-\\-3d-2e-f\end{array}$  & $-\dfrac{\alpha_{90}\beta_{89}}{\gamma_{26}\epsilon_2}$ & 177 & $\begin{array}{c}-a-b-2c-\\-3d--2e-f\end{array}$    & $-\alpha_{90}\beta_{39}\gamma_{26}\epsilon_2$                       & 183 \\[3mm]\hline
$r$                & \boldmath$N_{c+d+e,r}$                        & N   & $\begin{array}{c}-a-2b-2c-\\-3d-2e-f\end{array}$   & $-\dfrac{\alpha_{90}\beta_{49}}{\gamma_{26}\delta_{13}\epsilon_2}$  & 202 \\[3mm]\hline
$-a-c-d-e$         & $\alpha_4$                                    & 019 & $r$                  & \boldmath$N_{a+b+2c+2d+e+f,r}$                                      & N   \\[3mm]\hline
$-b-c-d-e$         & $\beta_{34}\delta_1$                          & 031 & $\begin{array}{c}-a-b-2c-\\-2d-2e-f\end{array}$    & $\alpha_{90}\beta_{89}\epsilon_2$                                   & 164 \\[3mm]\hline
$-b-c-2d-e$        & $\beta_{34}\gamma_{13}\delta_1$               & 043 & $\begin{array}{c}-a-b-2c-\\-3d-2e-f\end{array}$    & $\alpha_{90}\beta_{89}\gamma_{26}\epsilon_2$                        & 177 \\[3mm]\hline
$\begin{array}{c}-a-b-2c-\\-2d-e\end{array}$     & $-\dfrac{\alpha_{36}\beta_{2345}}{\gamma_{13}\delta_1}$ & 065 & $\begin{array}{c}-a-2b-2c-\\-3d-2e-f\end{array}$   & $\alpha_{90}\beta_{89}\gamma_{26}\delta_{13}\epsilon_2$             & 196 \\[3mm]\hline
$-c-d-e-f$         & $-\beta_{36}\gamma_{24}\epsilon_1$            & 072 & $r$                  & \boldmath$N_{b+c+2d+2e+f,r}$                                        & N   \\
\hline   
\end{tabular}
\vskip5mm

\centerline{Table 4, part 6 of 6}

\textbf{Proof.}
According to \cite[Theorem 4.1.2]{Car72}, the structure constants of a simple Lie algebra of type $\Phi$ over
$\mathbb{C}$ satisfy the following relations:

(i) $\displaystyle{N_{s,r}=-N_{r,s},\ r,s\in\Phi;}$
\bigskip

(ii) $\displaystyle{\frac{N_{r_1,r_2}}{(r_3,r_3)}=
\frac{N_{r_2,r_3}}{(r_1,r_1)}=\frac{N_{r_3,r_1}}{(r_2,r_2)},}$
\medskip

\noindent
if $\displaystyle{r_1,r_2,r_3\in\Phi}$ and
$\displaystyle{r_1+r_2+r_3=0;}$
\bigskip

(iii) $\displaystyle{N_{r,s}N_{-r,-s}=-(p+1)^2},$
\medskip

\noindent
if $r,s,r+s\in\Phi;$
\bigskip

(iv)
$\displaystyle{\frac{N_{r_1,r_2}N_{r_3,r_4}}{(r_1+r_2,r_1+r_2)}+
\frac{N_{r_2,r_3}N_{r_1,r_4}}{(r_2+r_3,r_2+r_3)}+
\frac{N_{r_3,r_1}N_{r_2,r_4}}{(r_3+r_1,r_3+r_1)}=0,}$
\medskip    

\noindent
if $\displaystyle{r_1,r_2,r_3,r_4\in \Phi}$ satisfy
$\displaystyle{r_1+r_2+r_3+r_4=0}$ and if no pairs are opposite.
\vskip3mm

Since the scalar square of any root of $E_6$ is equal to 2, the denominator in all relations can be omitted.
The expression $\theta_{ij\ldots}$ below always means the product $\theta_i \theta_j\ldots\,.$ For example, 
$\theta_{12}=\theta_1\theta_2.$
\medskip

001. We have $a+(c)+(-a-c) = 0,$ therefore (see Table 2)
$$
N_{a,c} =  N_{c,-a-c} = N_{-a-c,a} = \alpha_1.\eqno{(001)}
$$

002. We have $c+(d)+(-c-d) = 0,$ therefore (see Table 2)
$$
N_{c,d} =  N_{d,-c-d} = N_{-c-d,c} = \beta_1.\eqno{(002)}
$$

003. We have $a+(c+d)+(-a-c-d) = 0,$ therefore (see Table 2)
$$
N_{a,c+d} =  N_{c+d,-a-c-d} = N_{-a-c-d,a} = \alpha_2.\eqno{(003)}
$$

004. We have $(a+c)+d+(-a)+(-c-d)$ and $d+(-a)\notin E_6,$ hence
$$
S = N_{a+c,d} N_{-a,-c-d} + N_{-a,a+c} N_{d,-c-d} = 0.
$$
Since
$$
N_{a+c,d} N_{-a,-c-d}  = | \text{Table 2} | = -\alpha_2 N_{a+c,d} ,
$$
$$
N_{-a,a+c} N_{d,-c-d} = | \text{Formulas (001) and (002)} | = \alpha_1(\beta_1) = \alpha_1\beta_1,
$$
we have
$$
S = -\alpha_2 N_{a+c,d} + \alpha_1\beta_1 = 0
$$
and therefore
$$
N_{a+c,d} =\alpha_{12}\beta_1.\eqno{(004)}
$$

005. We have $(a+c)+d+(-a-c-d) = 0,$ therefore (see formula (004))
$$
N_{a+c,d} =  N_{d,-a-c-d} = N_{-a-c-d,a+c}=\alpha_{12}\beta_1.\eqno{(005)}
$$

006. We have $d+(b)+(-b-d) = 0,$ therefore (see Table 2)
$$
N_{d,b} =  N_{b,-b-d} = N_{-b-d,d}=\gamma_1.\eqno{(006)}
$$

007. We have $c+(b+d)+(-b-c-d) = 0,$ therefore (see Table 2)
$$
N_{c,b+d} =  N_{b+d,-b-c-d} = N_{-b-c-d,c}= \beta_2.\eqno{(007)}
$$

008. We have $(c+d)+b+(-c)+(-b-d)$ and $b+(-c) \notin E_6,$ hence 
$$
S = N_{c+d,b} N_{-c,-b-d} + N_{-c,c+d} N_{b,-b-d} = 0.
$$
Since
$$
N_{c+d,b} N_{-c,-b-d} = | \text{Table 2} | = -\beta_2 N_{c+d,b}, 
$$
$$
N_{-c,c+d} N_{b,-b-d} = | \text{Formulas (002) and (006)} | = \beta_1(\gamma_1) = \beta_1\gamma_1,
$$
we have
$$
S = -\beta_2 N_{c+d,b} + \beta_1\gamma_1 = 0
$$
and therefore
$$
N_{c+d,b} = \beta_{12}\gamma_1 .\eqno{(008)}
$$

009. We have $(c+d)+b+(-b-c-d) = 0,$ therefore (see formula (008))
$$
N_{c+d,b} =  N_{b,-b-c-d} = N_{-b-c-d,c+d} = \beta_{12}\gamma_1.\eqno{(009)}
$$

010. We have $a+(b+c+d)+(-a-b-c-d) = 0,$ therefore (see Table 2)
$$
N_{a,b+c+d} =  N_{b+c+d,-a-b-c-d} = N_{-a-b-c-d,a} = \alpha_3.\eqno{(010)}
$$

011. We have $(a+c)+(b+d)+(-a)+(-b-c-d) = 0$ and $(b+d)+(-a)\notin E_6,$ hence 
$$
S = N_{a+c,b+d}N_{-a,-b-c-d}+N_{-a,a+c}N_{b+d,-b-c-d} = 0.
$$
Since
$$
N_{a+c,b+d}N_{-a,-b-c-d} = |\text{Table 2}| = -\alpha_3 N_{a+c,b+d}, 
$$
$$
N_{-a,a+c}N_{b+d,-b-c-d} = |\text{Formulas (001) and (007)}| = \alpha_1\beta_2,
$$
we have
$$
S = -\alpha_3N_{a+c,b+d} +\alpha_1\beta_2 = 0
$$
and therefore
$$
N_{a+c,b+d} = \alpha_{13}\beta_2.\eqno{(011)}
$$

012. We have $(a+c)+(b+d)+(-a-b-c-d) = 0,$ therefore (see formula (011))
$$
N_{a+c,b+d} =  N_{b+d,-a-b-c-d} = N_{-a-b-c-d,a+c}=\alpha_{13}\beta_2.\eqno{(012)}
$$

013. We have $(a+c+d)+b+(-a)+(-b-c-d) = 0$ and $b+(-a)\notin E_6,$ hence 
$$
S = N_{a+c+d,b}N_{-a,-b-c-d}+N_{-a,a+c+d}N_{b,-b-c-d} = 0.
$$
Since
$$
N_{a+c+d,b}N_{-a,-b-c-d} = |\text{Table 2}| = -\alpha_3 N_{a+c+d,b},  
$$
$$
N_{-a,a+c+d}N_{b,-b-c-d} = |\text{Formulas (003) and (009)}| = \alpha_2(\beta_{12}\gamma_1) = \alpha_2\beta_{12}\gamma_1,
$$
we have
$$
S = -\alpha_3N_{a+c+d,b} + \alpha_2\beta_{12}\gamma_1 = 0
$$
and therefore 
$$
N_{a+c+d,b} = \alpha_{23}\beta_{12}\gamma_1.\eqno{(013)}
$$

014. We have $(a+c+d)+b+(-a-b-c-d) = 0,$ therefore (see formula (013))
$$
N_{a+c+d,b} =  N_{b,-a-b-c-d} = N_{-a-b-c-d,a+c+d}= \alpha_{23}\beta_{12}\gamma_1.\eqno{(014)}
$$

015. We have $d+(e)+(-d-e) = 0,$ therefore (see Table 2)
$$
N_{d,e} =  N_{e,-d-e} = N_{-d-e,d} = \gamma_2.\eqno{(015)}
$$

016. We have $c+(d+e)+(-c-d-e) = 0,$ therefore (see Table 2)
$$
N_{c,d+e} =  N_{d+e,-c-d-e} = N_{-c-d-e,c} = \beta_3.\eqno{(016)}
$$

017. We have $(c+d)+e+(-c)+(-d-e) = 0$ and $e+(-c)\notin E_6,$ hence 
$$
S = N_{c+d,e}N_{-c,-d-e} + N_{-c,c+d}N_{e,-d-e} = 0.
$$
Since
$$
N_{c+d,e}N_{-c,-d-e} = | \text{Table 2} | = -\beta_3 N_{c+d,e}, 
$$
$$
N_{-c,c+d}N_{e,-d-e} = | \text{Formulas (002) and (015)} | = \beta_1(\gamma_2) = \beta_1\gamma_2,
$$
we have
$$
S = -\beta_3 N_{c+d,e} + \beta_1\gamma_2 = 0
$$
and therefore
$$
N_{c+d,e} = \beta_{13}\gamma_2.\eqno{(017)}
$$

018. We have $(c+d)+e+(-c-d-e) = 0,$ therefore (see formula (017))
$$
N_{c+d,e} =  N_{e,-c-d-e} = N_{-c-d-e,c+d} = \beta_{13}\gamma_2.\eqno{(018)}
$$

019. We have $a+(c+d+e)+(-a-c-d-e) = 0,$ therefore (see Table 2)
$$
N_{a,c+d+e} =  N_{c+d+e,-a-c-d-e} = N_{-a-c-d-e,a} = \alpha_4.\eqno{(019)}
$$

020. We have $(a+c+d)+e+(-a)+(-c-d-e) = 0$ and $e+(-a)\notin E_6,$ hence 
$$
S = N_{a+c+d,e}N_{-a,-c-d-e} + N_{-a,a+c+d}N_{e,-c-d-e} = 0.
$$
Since
$$
N_{a+c+d,e}N_{-a,-c-d-e} = | \text{Table 2} | = -\alpha_4 N_{a+c+d,e},
$$
$$
N_{-a,a+c+d}N_{e,-c-d-e} = | \text{Formulas (003) and (018)} | = \alpha_2(\beta_{13}\gamma_2) = \alpha_2\beta_{13}\gamma_2,
$$
we have
$$
S = -\alpha_4N_{a+c+d,e} + \alpha_2\beta_{13}\gamma_2 = 0
$$
and therefore
$$
N_{a+c+d,e} = \alpha_{24}\beta_{13}\gamma_2.\eqno{(020)}
$$

021. We have $(a+c+d)+e+(-a-c-d-e) = 0,$ therefore (see formula (020))
$$
N_{a+c+d,e} =  N_{e,-a-c-d-e} = N_{-a-c-d-e,a+c+d} = \alpha_{24}\beta_{13}\gamma_2.\eqno{(021)}
$$

022. We have $(a+c)+(d+e)+(-a)+(-c-d-e) = 0$ and $(d+e)+(-a)\notin E_6,$ hence 
$$
S = N_{a+c,d+e} N_{-a,-c-d-e} + N_{-a,a+c} N_{d+e,-c-d-e} = 0.
$$
Since
$$
N_{a+c,d+e} N_{-a,-c-d-e} = | \text{Table 2} | = -\alpha_4 N_{a+c,d+e}, 
$$
$$
N_{-a,a+c} N_{d+e,-c-d-e} = | \text{Formulas (001) and (016)} | = \alpha_1(\beta_3) = \alpha_1\beta_3,
$$
we have
$$
S = -\alpha_4 N_{a+c,d+e} + \alpha_1\beta_3 = 0
$$
and therefore
$$
 N_{a+c,d+e} = \alpha_{14}\beta_3.\eqno{(022)}
$$

023. We have $(a+c)+(d+e)+(-a-c-d-e) = 0,$ therefore (see formula (022))
$$
N_{a+c,d+e} =  N_{d+e,-a-c-d-e} = N_{-a-c-d-e,a+c} = \alpha_{14}\beta_3.\eqno{(023)}
$$

024. We have $b+(d+e)+(-b-d-e) = 0,$ therefore (see Table 2)
$$
N_{b,d+e} =  N_{d+e,-b-d-e} = N_{-b-d-e,b} = \delta_1.\eqno{(024)}
$$

025. We have $(b+d)+e+(-b)+(-d-e)$ and $e+(-b) \notin E_6,$ hence 
$$
S = N_{b+d,e} N_{-b,-d-e} + N_{-b,b+d} N_{e,-d-e} = 0.
$$
Since
$$
N_{b+d,e} N_{-b,-d-e} = | \text{Table 2} | = -\delta_1 N_{b+d,e}, 
$$
$$
N_{-b,b+d} N_{e,-d-e} = | \text{Formulas (006) and (015)} | = -\gamma_1(\gamma_2) = -\gamma_{12},
$$
we have
$$
S = -\delta_1 N_{b+d,e} -\gamma_{12} = 0
$$
and therefore
$$
N_{b+d,e} = -\gamma_{12}\delta_1 .\eqno{(025)}
$$

026. We have $(b+d)+e+(-b-d-e) = 0,$ therefore (see formula (025))
$$
N_{b+d,e} =  N_{e,-b-d-e} = N_{-b-d-e,b+d} = -\gamma_{12}\delta_1.\eqno{(026)}
$$

027. We have $c+(b+d+e)+(-b-c-d-e) = 0,$ therefore (see Table 2)
$$
N_{c,b+d+e} =  N_{b+d+e,-b-c-d-e} = N_{-b-c-d-e,c} =\beta_4.\eqno{(027)}
$$

028. We have $(b+c+d)+e+(-c)+(-b-d-e) = 0$ and $e+(-c)\notin E_6,$ hence 
$$
S = N_{b+c+d,e}N_{-c,-b-d-e} + N_{-c,b+c+d}N_{e,-b-d-e} = 0.
$$
Since
$$
N_{b+c+d,e}N_{-c,-b-d-e} = | \text{Table 2} | = -\beta_4 N_{b+c+d,e},
$$
$$
N_{-c,b+c+d}N_{e,-b-d-e} = | \text{Formulas (007) and (026)} | = \beta_2(-\gamma_{12}\delta_1) = -\beta_2\gamma_{12}\delta_1,
$$
we have
$$
S = -\beta_4 N_{b+c+d,e} -\beta_2\gamma_{12}\delta_1 = 0
$$
and therefore
$$
N_{b+c+d,e} = -\beta_{24}\gamma_{12}\delta_1.\eqno{(028)}
$$

029. We have $(b+c+d)+e+(-b-c-d-e) = 0,$ therefore (see formula (028))
$$
N_{b+c+d,e} =  N_{e,-b-c-d-e} = N_{-b-c-d-e,b+c+d} = -\beta_{24}\gamma_{12}\delta_1.\eqno{(029)}
$$

030. We have $b+(c+d+e)+(-c)+(-b-d-e) = 0$ and $(-c)+b\notin E_6,$ hence 
$$
S = N_{b,c+d+e}N_{-c,-b-d-e}+N_{c+d+e,-c}N_{b,-b-d-e} = 0.
$$
Since
$$
N_{b,c+d+e}N_{-c,-b-d-e} = | \text{Table 2} | = -\beta_4 N_{b,c+d+e}, 
$$
$$
N_{c+d+e,-c}N_{b,-b-d-e} = | \text{Formulas (016) and (024)} | = -\beta_3(-\delta_1) = \beta_3\delta_1,
$$
we have
$$
S = -\beta_4 N_{b,c+d+e} + \beta_3\delta_1 = 0
$$
and therefore
$$
N_{b,c+d+e} = \beta_{34}\delta_1.\eqno{(030)}
$$

031. We have $b+(c+d+e)+(-b-c-d-e) = 0,$ therefore (see formula (030))
$$
N_{b,c+d+e} =  N_{c+d+e,-b-c-d-e} = N_{-b-c-d-e,b} = \beta_{34}\delta_1.\eqno{(031)}
$$

032. We have $a+(b+c+d+e)+(-a-b-c-d-e) = 0,$ therefore (see Table 2) 
$$
N_{a,b+c+d+e} =  N_{b+c+d+e,-a-b-c-d-e} = N_{-a-b-c-d-e,a} = \alpha_5.\eqno{(032)}
$$

033. We have $(a+b+c+d)+e+(-a)+(-b-c-d-e) = 0$ and $e+(-a)\notin E_6,$ hence 
$$
S = N_{a+b+c+d,e}N_{-a,-b-c-d-e} + N_{-a,a+b+c+d}N_{e,-b-c-d-e} = 0.
$$
Since
$$
N_{a+b+c+d,e}N_{-a,-b-c-d-e} = | \text{Table 2} | = -\alpha_5 N_{a+b+c+d,e},
$$
$$
N_{-a,a+b+c+d}N_{e,-b-c-d-e} = | \text{Formulas (010) and (029)} | = \alpha_3(-\beta_{24}\gamma_{12}\delta_1) = -\alpha_3\beta_{24}\gamma_{12}\delta_1,
$$
we have
$$
S = -\alpha_5 N_{a+b+c+d,e} - \alpha_3\beta_{24}\gamma_{12}\delta_1 = 0
$$
and therefore
$$
N_{a+b+c+d,e} = -\alpha_{35}\beta_{24}\gamma_{12}\delta_1.\eqno{(033)}
$$

034. We have $(a+b+c+d)+e+(-a-b-c-d-e) = 0,$ therefore (see formula (033)) 
$$
N_{a+b+c+d,e} =  N_{e,-a-b-c-d-e} = N_{-a-b-c-d-e,a+b+c+d} = -\alpha_{35}\beta_{24}\gamma_{12}\delta_1.\eqno{(034)}
$$

035. We have $(a+c)+(b+d+e)+(-a)+(-b-c-d-e) = 0$ and $(b+d+e)+(-a)\notin E_6,$ hence  
$$
S = N_{a+c,b+d+e} N_{-a,-b-c-d-e}+N_{-a,a+c} N_{b+d+e,-b-c-d-e} = 0.
$$
Since
$$
N_{a+c,b+d+e} N_{-a,-b-c-d-e} = | \text{Table 2}  | = -\alpha_5 N_{a+c,b+d+e}, 
$$
$$
N_{-a,a+c} N_{b+d+e,-b-c-d-e} = | \text{Formulas (001) and (027)} | = \alpha_1(\beta_4) = \alpha_1\beta_4,
$$
we have
$$
S = -\alpha_5 N_{a+c,b+d+e} + \alpha_1\beta_4 = 0
$$
and therefore
$$
N_{a+c,b+d+e} = \alpha_{15}\beta_4.\eqno{(035)}
$$

036. We have $(a+c)+(b+d+e)+(-a-b-c-d-e) = 0,$ therefore (see formula (035)) 
$$
N_{a+c,b+d+e} =  N_{b+d+e,-a-b-c-d-e} = N_{-a-b-c-d-e,a+c} = \alpha_{15}\beta_4.\eqno{(036)}
$$

037. We have $b+(a+c+d+e)+(-a)+(-b-c-d-e) = 0$ and $(-a)+b\notin E_6,$ hence  
$$
S = N_{b,a+c+d+e} N_{-a,-b-c-d-e} + N_{a+c+d+e,-a} N_{b,-b-c-d-e} = 0.
$$
Since
$$
N_{b,a+c+d+e} N_{-a,-b-c-d-e} = | \text{Table 2}  | = -\alpha_5 N_{b,a+c+d+e}, 
$$
$$
N_{a+c+d+e,-a} N_{b,-b-c-d-e} = | \text{Formulas (019) and (031)} | = -\alpha_4 (-\beta_{34}\delta_1) = \alpha_4\beta_{34}\delta_1,
$$
we have
$$
S = -\alpha_5 N_{b,a+c+d+e} + \alpha_4\beta_{34}\delta_1 = 0
$$
and therefore
$$
N_{b,a+c+d+e} = \alpha_{45}\beta_{34}\delta_1.\eqno{(037)}
$$

038. We have $b+(a+c+d+e)+(-a-b-c-d-e) = 0,$ therefore (see formula (037))  
$$
N_{b,a+c+d+e} =  N_{a+c+d+e,-a-b-c-d-e} = N_{-a-b-c-d-e,b} = \alpha_{45}\beta_{34}\delta_1.\eqno{(038)}
$$

039. We have $d+(b+c+d+e)+(-b-c-2d-e) = 0,$ therefore (see Table 2) 
$$
N_{d,b+c+d+e} =  N_{b+c+d+e,-b-c-2d-e} = N_{-b-c-2d-e,d} = \gamma_3.\eqno{(039)}
$$

040. We have $(b+c+d)+(d+e)+(-d)+(-b-c-d-e) = 0$ and $(-d)+(b+c+d)\notin E_6,$ hence 
$$
S = N_{b+c+d,d+e} N_{-d,-b-c-d-e} + N_{d+e,-d}N_{b+c+d,-b-c-d-e} = 0.
$$
Since
$$
N_{b+c+d,d+e} N_{-d,-b-c-d-e} = | \text{Table 2} | = -\gamma_3 N_{b+c+d,d+e}, 
$$
$$
N_{d+e,-d} N_{b+c+d,-b-c-d-e} = | \text{Formulas (015) and (029)} | = -\gamma_2(\beta_{24}\gamma_{12}\delta_1) = -\beta_{24}\gamma_{1}\delta_1,
$$
we have
$$
S = -\gamma_3 N_{b+c+d,d+e} -\beta_{24}\gamma_{1}\delta_1 = 0
$$
and therefore
$$
N_{b+c+d,d+e} = -\beta_{24}\gamma_{13}\delta_1.\eqno{(040)}
$$

041. We have $(b+c+d)+(d+e)+(-b-c-2d-e) = 0,$ therefore (see formula (040)) 
$$
N_{b+c+d,d+e} =  N_{d+e,-b-c-2d-e} = N_{-b-c-2d-e,b+c+d} = -\beta_{24}\gamma_{13}\delta_1.\eqno{(041)}
$$

042. We have $(b+d)+(c+d+e)+(-d)+(-b-c-d-e) = 0$ and $(c+d+e)+(-d)\notin E_6,$ hence 
$$
S = N_{b+d,c+d+e}N_{-d,-b-c-d-e} + N_{-d,b+d}N_{c+d+e,-b-c-d-e} = 0.
$$
Since
$$
N_{b+d,c+d+e}N_{-d,-b-c-d-e} = | \text{Table 2} | = -\gamma_3 N_{b+d,c+d+e}, 
$$
$$
N_{-d,b+d}N_{c+d+e,-b-c-d-e} = | \text{Formulas (006) and (031)} | = \gamma_1(\beta_{34}\delta_1) = \beta_{34}\gamma_1\delta_1,
$$
we have
$$
S = -\gamma_3 N_{b+d,c+d+e} + \beta_{34}\gamma_1\delta_1 = 0
$$
and therefore
$$
N_{b+d,c+d+e} = \beta_{34}\gamma_{13}\delta_1.\eqno{(042)}
$$

043. We have $(b+d)+(c+d+e)+(-b-c-2d-e) = 0,$ therefore (see formula (042)) 
$$
N_{b+d,c+d+e} =  N_{c+d+e,-b-c-2d-e} = N_{-b-c-2d-e,b+d}=\beta_{34}\gamma_{13}\delta_1.\eqno{(043)}
$$

044. We have $(c+d)+(b+d+e)+(-d)+(-b-c-d-e) = 0$ and $(b+d+e)+(-d)\notin E_6,$ hence 
$$
S = N_{c+d,b+d+e}N_{-d,-b-c-d-e} + N_{-d,c+d}N_{b+d+e,-b-c-d-e} = 0.
$$
Since
$$
N_{c+d,b+d+e}N_{-d,-b-c-d-e} = | \text{Table 2} | = -\gamma_3 N_{c+d,b+d+e}, 
$$
$$
N_{-d,c+d}N_{b+d+e,-b-c-d-e} = | \text{Formulas (002) and (027)} | = -\beta_1(\beta_4) = -\beta_{14},
$$
we have
$$
S = -\gamma_3 N_{c+d,b+d+e} - \beta_{14} = 0
$$
and therefore
$$
N_{c+d,b+d+e} = -\beta_{14}\gamma_3.\eqno{(044)}
$$

045. We have $(c+d)+(b+d+e)+(-b-c-2d-e) = 0,$ therefore (see formula (044)) 
$$
N_{c+d,b+d+e} =  N_{b+d+e,-b-c-2d-e} = N_{-b-c-2d-e,c+d}=-\beta_{14}\gamma_3.\eqno{(045)}
$$

046. We have $a+(b+c+2d+e)+(-a-b-c-2d-e) = 0,$ therefore (see Table 2)
$$
N_{a,b+c+2d+e} =  N_{b+c+2d+e,-a-b-c-2d-e} = N_{-a-b-c-2d-e,a} =\alpha_6.\eqno{(046)}
$$

047. We have $d+(a+b+c+d+e)+(-a)+(-b-c-2d-e) = 0$ and $d+(-a)\notin E_6,$ hence 
$$
S = N_{d,a+b+c+d+e}N_{-a,-b-c-2d-e} + N_{a+b+c+d+e,-a}N_{d,-b-c-2d-e} = 0.
$$
Since
$$
N_{d,a+b+c+d+e}N_{-a,-b-c-2d-e} = | \text{Table 2} | = -\alpha_6 N_{d,a+b+c+d+e}, 
$$
$$
N_{a+b+c+d+e,-a}N_{d,-b-c-2d-e} = | \text{Formulas (032) and (039)} | = -\alpha_5(-\gamma_3) = \alpha_5\gamma_3,
$$
we have
$$
S = -\alpha_6 N_{d,a+b+c+d+e} + \alpha_5\gamma_3 = 0
$$
and therefore
$$
N_{d,a+b+c+d+e} = \alpha_{56}\gamma_3.\eqno{(047)}
$$

048. We have $d+(a+b+c+d+e)+(-a-b-c-2d-e) = 0,$ therefore (see formula (047)) 
$$
N_{d,a+b+c+d+e} =  N_{a+b+c+d+e,-a-b-c-2d-e} = N_{-a-b-c-2d-e,d}=\alpha_{56}\gamma_3.\eqno{(048)}
$$

049. We have $(a+b+c+d)+(d+e)+(-a)+(-b-c-2d-e) = 0$ and $(d+e)+(-a)\notin E_6,$ hence 
$$
S = N_{a+b+c+d,d+e}N_{-a,-b-c-2d-e} + N_{-a,a+b+c+d}N_{d+e,-b-c-2d-e} = 0.
$$
Since
$$
N_{a+b+c+d,d+e}N_{-a,-b-c-2d-e} = | \text{Table 2} | = -\alpha_6 N_{a+b+c+d,d+e},
$$
$$
N_{-a,a+b+c+d}N_{d+e,-b-c-2d-e} = | \text{Formulas (010) and (041)} | = \alpha_3(-\beta_{24}\gamma_{13}\delta_1) = -\alpha_3\beta_{24}\gamma_{13}\delta_1,
$$
we have
$$
S = -\alpha_6 N_{a+b+c+d,d+e} - \alpha_3\beta_{24}\gamma_{13}\delta_1 = 0
$$
and therefore
$$
N_{a+b+c+d,d+e} = -\alpha_{36}\beta_{24}\gamma_{13}\delta_1.\eqno{(049)}
$$

050. We have $(a+b+c+d)+(d+e)+(-a-b-c-2d-e) = 0,$ therefore (see formula (049)) 
$$
N_{a+b+c+d,d+e} =  N_{d+e,-a-b-c-2d-e} = N_{-a-b-c-2d-e,a+b+c+d} = -\alpha_{36}\beta_{24}\gamma_{13}\delta_1.\eqno{(050)}
$$

051. We have $(b+d)+(a+c+d+e)+(-a)+(-b-c-2d-e) = 0$ and $(b+d)+(-a)\notin E_6,$ hence 
$$
S = N_{b+d,a+c+d+e}N_{-a,-b-c-2d-e} + N_{a+c+d+e,-a}N_{b+d,-b-c-2d-e} = 0.
$$
Since
$$
N_{b+d,a+c+d+e}N_{-a,-b-c-2d-e} = | \text{Table 2} | = -\alpha_6 N_{b+d,a+c+d+e}, 
$$
$$
N_{a+c+d+e,-a}N_{b+d,-b-c-2d-e} = | \text{Formulas (019) and (043)} | = -\alpha_4(-\beta_{34}\gamma_{13}\delta_1) = \alpha_4\beta_{34}\gamma_{13}\delta_1,
$$
we have
$$
S = -\alpha_6 N_{b+d,a+c+d+e} + \alpha_4\beta_{34}\gamma_{13}\delta_1 = 0
$$
and therefore
$$
N_{b+d,a+c+d+e} =\alpha_{46}\beta_{34}\gamma_{13}\delta_1.\eqno{(051)}
$$

052. We have $(b+d)+(a+c+d+e)+(-a-b-c-2d-e) = 0,$ therefore (see formula (051)) 
$$
N_{b+d,a+c+d+e} =  N_{a+c+d+e,-a-b-c-2d-e} = N_{-a-b-c-2d-e,b+d} = \alpha_{46}\beta_{34}\gamma_{13}\delta_1.\eqno{(052)}
$$

053. We have $(a+c+d)+(b+d+e)+(-a)+(-b-c-2d-e) = 0$ and $(b+d+e)+(-a)\notin E_6,$ hence 
$$
S = N_{a+c+d,b+d+e}N_{-a,-b-c-2d-e} + N_{-a,a+c+d}N_{b+d+e,-b-c-2d-e} = 0.
$$
Since
$$
N_{a+c+d,b+d+e}N_{-a,-b-c-2d-e} = | \text{Table 2} | = -\alpha_6 N_{a+c+d,b+d+e},
$$
$$
N_{-a,a+c+d}N_{b+d+e,-b-c-2d-e} = | \text{Formulas (003) and (045)} | = \alpha_2(-\beta_{14}\gamma_3) = -\alpha_2\beta_{14}\gamma_3,
$$
we have
$$
S = -\alpha_6 N_{a+c+d,b+d+e} - \alpha_2\beta_{14}\gamma_3 = 0
$$
and therefore
$$
N_{a+c+d,b+d+e} = -\alpha_{26}\beta_{14}\gamma_3.\eqno{(053)}
$$

054. We have $(a+c+d)+(b+d+e)+(-a-b-c-2d-e) = 0,$ therefore (see formula (053))
$$
N_{a+c+d,b+d+e} =  N_{b+d+e,-a-b-c-2d-e} = N_{-a-b-c-2d-e,a+c+d} = -\alpha_{26}\beta_{14}\gamma_3.\eqno{(054)}
$$

055. We have $c+(a+b+c+2d+e)+(-a-b-2c-2d-e) = 0,$ therefore (see Table 2) 
$$
N_{c,a+b+c+2d+e} =  N_{a+b+c+2d+e,-a-b-2c-2d-e} = N_{-a-b-2c-2d-e,c} = \beta_5.\eqno{(055)}
$$

056. We have $(c+d)+(a+b+c+d+e)+(-c)+(-a-b-c-2d-e) = 0$ and $(a+b+c+d+e)+(-c)\notin E_6,$ hence 
$$
S = N_{c+d,a+b+c+d+e}N_{-c,-a-b-c-2d-e} + N_{-c,c+d}N_{a+b+c+d+e,-a-b-c-2d-e} = 0.
$$
Since
$$
N_{c+d,a+b+c+d+e}N_{-c,-a-b-c-2d-e} = | \text{Table 2} | = -\beta_5 N_{c+d,a+b+c+d+e}, 
$$
$$
N_{-c,c+d}N_{a+b+c+d+e,-a-b-c-2d-e} = | \text{Formulas (002) and (048)} | = \beta_1(\alpha_{56}\gamma_3) = \alpha_{56}\beta_1\gamma_3,
$$
we have
$$
S = -\beta_5 N_{c+d,a+b+c+d+e} + \alpha_{56}\beta_1\gamma_3 = 0
$$
and therefore
$$
N_{c+d,a+b+c+d+e} = \alpha_{56}\beta_{15}\gamma_3.\eqno{(056)}
$$

057. We have $(c+d)+(a+b+c+d+e)+(-a-b-2c-2d-e) = 0,$ therefore (see formula (056)) 
$$
N_{c+d,a+b+c+d+e} =  N_{a+b+c+d+e,-a-b-2c-2d-e} = N_{-a-b-2c-2d-e,c+d} = \alpha_{56}\beta_{15}\gamma_3.\eqno{(057)}
$$

058. We have $(a+c)+(b+c+2d+e)+(-c)+(-a-b-c-2d-e) = 0$ and $(b+c+2d+e)+(-c)\notin E_6,$ hence 
$$
S = N_{a+c,b+c+2d+e}N_{-c,-a-b-c-2d-e} + N_{-c,a+c}N_{b+c+2d+e,-a-b-c-2d-e} = 0.
$$
Since
$$
N_{a+c,b+c+2d+e}N_{-c,-a-b-c-2d-e} = | \text{Table 2} | = -\beta_5 N_{a+c,b+c+2d+e}, 
$$
$$
N_{-c,a+c}N_{b+c+2d+e,-a-b-c-2d-e} = | \text{Formulas (001) and (046)} | = -\alpha_1(\alpha_6) = -\alpha_{16},
$$
we have
$$
S = -\beta_5 N_{a+c,b+c+2d+e} - \alpha_{16} = 0
$$
and therefore
$$
N_{a+c,b+c+2d+e} = -\alpha_{16}\beta_5.\eqno{(058)}
$$

059. We have $(a+c)+(b+c+2d+e)+(-a-b-2c-2d-e) = 0,$ therefore (see formula (058)) 
$$
N_{a+c,b+c+2d+e} =  N_{b+c+2d+e,-a-b-2c-2d-e} = N_{-a-b-2c-2d-e,a+c} = -\alpha_{16}\beta_5.\eqno{(059)}
$$

060. We have $(a+c+d)+(b+c+d+e)+(-c)+(-a-b-c-2d-e) = 0$ and $(a+c+d)+(-c)\notin E_6,$ hence 
$$
S = N_{a+c+d,b+c+d+e}N_{-c,-a-b-c-2d-e} + N_{b+c+d+e,-c}N_{a+c+d,-a-b-c-2d-e} = 0.
$$
Since
$$
N_{a+c+d,b+c+d+e}N_{-c,-a-b-c-2d-e} = | \text{Table 2} | = -\beta_5 N_{a+c+d,b+c+d+e}, 
$$
$$
N_{b+c+d+e,-c}N_{a+c+d,-a-b-c-2d-e} = | \text{Formulas (027) and (054)} | = -\beta_4(\alpha_{26}\beta_{14}\gamma_3) = -\alpha_{26}\beta_{1}\gamma_3,
$$
we have
$$
S = -\beta_5 N_{a+c+d,b+c+d+e} - \alpha_{26}\beta_{1}\gamma_3 = 0
$$
and therefore
$$
N_{a+c+d,b+c+d+e} = -\alpha_{26}\beta_{15}\gamma_3.\eqno{(060)}
$$

061. We have $(a+c+d)+(b+c+d+e)+(-a-b-2c-2d-e) = 0,$ therefore (see formula (060)) 
$$
N_{a+c+d,b+c+d+e} =  N_{b+c+d+e,-a-b-2c-2d-e} = N_{-a-b-2c-2d-e,a+c+d} = -\alpha_{26}\beta_{15}\gamma_3.\eqno{(061)}
$$

062. We have $(b+c+d)+(a+c+d+e)+(-c)+(-a-b-c-2d-e) = 0$ and $(a+c+d+e)+(-c)\notin E_6,$ hence 
$$
S = N_{b+c+d,a+c+d+e}N_{-c,-a-b-c-2d-e} + N_{-c,b+c+d}N_{a+c+d+e,-a-b-c-2d-e} = 0.
$$
Since
$$
N_{b+c+d,a+c+d+e}N_{-c,-a-b-c-2d-e} = | \text{Table 2} | = -\beta_5 N_{b+c+d,a+c+d+e}, 
$$
$$
N_{-c,b+c+d}N_{a+c+d+e,-a-b-c-2d-e} = | \text{Formulas (007) and (052)} | = \beta_2(\alpha_{46}\beta_{34}\gamma_{13}\delta_1) = 
$$
$$
=\alpha_{46}\beta_{234}\gamma_{13}\delta_1,
$$
we have
$$
S = -\beta_5 N_{b+c+d,a+c+d+e} + \alpha_{46}\beta_{234}\gamma_{13}\delta_1 = 0
$$
and therefore
$$
N_{b+c+d,a+c+d+e} = \alpha_{46}\beta_{2345}\gamma_{13}\delta_1.\eqno{(062)}
$$

063. We have $(b+c+d)+(a+c+d+e)+(-a-b-2c-2d-e) = 0,$ therefore (see formula (062)) 
$$
N_{b+c+d,a+c+d+e} =  N_{a+c+d+e,-a-b-2c-2d-e} = N_{-a-b-2c-2d-e,b+c+d} = \alpha_{46}\beta_{2345}\gamma_{13}\delta_1.\eqno{(063)}
$$

064. We have $(a+b+c+d)+(c+d+e)+(-c)+(-a-b-c-2d-e) = 0$ and $(a+b+c+d)+(-c)\notin E_6,$ hence 
$$
S = N_{a+b+c+d,c+d+e}N_{-c,-a-b-c-2d-e} + N_{c+d+e,-c}N_{a+b+c+d,-a-b-c-2d-e} = 0.
$$
Since
$$
N_{a+b+c+d,c+d+e}N_{-c,-a-b-c-2d-e} = | \text{Table 2} | = -\beta_5 N_{a+b+c+d,c+d+e},
$$
$$
N_{c+d+e,-c}N_{a+b+c+d,-a-b-c-2d-e} = | \text{Formulas (016) and (050)} | = -\beta_3(\alpha_{36}\beta_{24}\gamma_{13}\delta_1) = 
$$
$$
=-\alpha_{36}\beta_{234}\gamma_{13}\delta_1,
$$
we have
$$
S = -\beta_5 N_{a+b+c+d,c+d+e} - \alpha_{36}\beta_{234}\gamma_{13}\delta_1 = 0
$$
and therefore
$$
N_{a+b+c+d,c+d+e} = -\alpha_{36}\beta_{2345}\gamma_{13}\delta_1.\eqno{(064)}
$$

065. We have $(a+b+c+d)+(c+d+e)+(-a-b-2c-2d-e) = 0,$ therefore (see formula (064))
$$
N_{a+b+c+d,c+d+e} =  N_{c+d+e,-a-b-2c-2d-e} = N_{-a-b-2c-2d-e,a+b+c+d} = -\alpha_{36}\beta_{2345}\gamma_{13}\delta_1.\eqno{(065)}
$$

066. We have $e+(f)+(-e-f) = 0,$ therefore (see Table 2)
$$
N_{e,f} =  N_{f,-e-f} = N_{-e-f,e} =\epsilon_1.\eqno{(066)}
$$

067. We have $d+(e+f)+(-d-e-f) = 0,$ therefore (see Table 2)
$$
N_{d,e+f} =  N_{e+f,-d-e-f} = N_{-d-e-f,d} =\gamma_4.\eqno{(067)}
$$

068. We have $(d+e)+f+(-d)+(-e-f) = 0$ and $f+(-d)\notin E_6,$ hence 
$$
S = N_{d+e,f}N_{-d,-e-f} + N_{-d,d+e}N_{f,-e-f} = 0.
$$
Since
$$
N_{d+e,f}N_{-d,-e-f} = | \text{Table 2} | = -\gamma_4 N_{d+e,f}, 
$$
$$
N_{-d,d+e}N_{f,-e-f} = | \text{Formulas (015) and (066)} | = \gamma_2(\epsilon_1) = \gamma_2\epsilon_1,
$$
we have
$$
S = -\gamma_4 N_{d+e,f} + \gamma_2\epsilon_1 = 0
$$
and therefore
$$
N_{d+e,f} = \gamma_{24}\epsilon_1.\eqno{(068)}
$$

069. We have $(d+e)+f+(-d-e-f) = 0,$ therefore (see formula (068)) 
$$
N_{d+e,f} =  N_{f,-d-e-f} = N_{-d-e-f,d+e} = \gamma_{24}\epsilon_1.\eqno{(069)}
$$

070. We have $c+(d+e+f)+(-c-d-e-f) = 0,$ therefore (see Table 2)
$$
N_{c,d+e+f} =  N_{d+e+f,-c-d-e-f} = N_{-c-d-e-f,c} = \beta_6.\eqno{(070)}
$$

071. We have $(c+d+e)+f+(-c)+(-d-e-f) = 0$ and $f+(-c)\notin E_6,$ hence 
$$
S = N_{c+d+e,f}N_{-c,-d-e-f} + N_{-c,c+d+e}N_{f,-d-e-f} = 0.
$$
Since
$$
N_{c+d+e,f}N_{-c,-d-e-f} = | \text{Table 2} | = -\beta_6 N_{c+d+e,f}, 
$$
$$
N_{-c,c+d+e}N_{f,-d-e-f} = | \text{Formulas (016) and (069)} | = \beta_3(\gamma_{24}\epsilon_1) = \beta_3\gamma_{24}\epsilon_1,
$$
we have
$$
S = -\beta_6 N_{c+d+e,f} + \beta_3\gamma_{24}\epsilon_1 = 0
$$
and therefore
$$
N_{c+d+e,f} = \beta_{36}\gamma_{24}\epsilon_1.\eqno{(071)}
$$

072. We have $(c+d+e)+f+(-c-d-e-f) = 0,$ therefore (see formula (071))
$$
N_{c+d+e,f} =  N_{f,-c-d-e-f} = N_{-c-d-e-f,c+d+e} = \beta_{36}\gamma_{24}\epsilon_1.\eqno{(072)}
$$

073. We have $(c+d)+(e+f)+(-c)+(-d-e-f) = 0$ and $(e+f)+(-c)\notin E_6,$ hence 
$$
S = N_{c+d,e+f}N_{-c,-d-e-f} + N_{-c,c+d}N_{e+f,-d-e-f} = 0.
$$
Since
$$
N_{c+d,e+f}N_{-c,-d-e-f} = | \text{Table 2} | = -\beta_6 N_{c+d,e+f}, 
$$
$$
N_{-c,c+d}N_{e+f,-d-e-f} = | \text{Formulas (002) and (067)} | = \beta_1(\gamma_4) = \beta_1\gamma_4,
$$
we have
$$
S = -\beta_6 N_{c+d,e+f} + \beta_1\gamma_4 = 0
$$
and therefore
$$
N_{c+d,e+f} = \beta_{16}\gamma_4.\eqno{(073)}
$$

074. We have $(c+d)+(e+f)+(-c-d-e-f) = 0,$ therefore (see formula (073)) 
$$
N_{c+d,e+f} =  N_{e+f,-c-d-e-f} = N_{-c-d-e-f,c+d} = \beta_{16}\gamma_4.\eqno{(074)}
$$

075. We have $a+(c+d+e+f)+(-a-c-d-e-f) = 0,$ therefore (see Table 2)
$$
N_{a,c+d+e+f} =  N_{c+d+e+f,-a-c-d-e-f} = N_{-a-c-d-e-f,a} = \alpha_7.\eqno{(075)}
$$

076. We have $(a+c+d+e)+f+(-a)+(-c-d-e-f) = 0$ and $f+(-a)\notin E_6,$ hence 
$$
S = N_{a+c+d+e,f}N_{-a,-c-d-e-f} + N_{-a,a+c+d+e}N_{f,-c-d-e-f} = 0.
$$
Since
$$
N_{a+c+d+e,f}N_{-a,-c-d-e-f} = | \text{Table 2} | = -\alpha_7 N_{a+c+d+e,f}, 
$$
$$
N_{-a,a+c+d+e}N_{f,-c-d-e-f} = | \text{Formulas (019) and (072)} | = \alpha_4(\beta_{36}\gamma_{24}\epsilon_1) = \alpha_4\beta_{36}\gamma_{24}\epsilon_1,
$$
we have
$$
S = -\alpha_7 N_{a+c+d+e,f} + \alpha_4\beta_{36}\gamma_{24}\epsilon_1 = 0
$$
and therefore
$$
N_{a+c+d+e,f} = \alpha_{47}\beta_{36}\gamma_{24}\epsilon_1.\eqno{(076)}
$$

077. We have $(a+c+d+e)+f+(-a-c-d-e-f) = 0,$ therefore (see formula (076))
$$
N_{a+c+d+e,f} =  N_{f,-a-c-d-e-f} = N_{-a-c-d-e-f,a+c+d+e} = \alpha_{47}\beta_{36}\gamma_{24}\epsilon_1.\eqno{(077)}
$$

078. We have $(a+c+d)+(e+f)+(-a)+(-c-d-e-f) = 0$ and $(e+f)+(-a)\notin E_6,$ hence 
$$
S = N_{a+c+d,e+f}N_{-a,-c-d-e-f} + N_{-a,a+c+d}N_{e+f,-c-d-e-f} = 0.
$$
Since
$$
N_{a+c+d,e+f}N_{-a,-c-d-e-f} = | \text{Table 2} | = -\alpha_7 N_{a+c+d,e+f},  
$$
$$
N_{-a,a+c+d}N_{e+f,-c-d-e-f} = | \text{Formulas (003) and (074)} | = \alpha_2(\beta_{16}\gamma_4) = \alpha_2\beta_{16}\gamma_4,
$$
we have
$$
S = -\alpha_7 N_{a+c+d,e+f} + \alpha_2\beta_{16}\gamma_4 = 0
$$
and therefore
$$
N_{a+c+d,e+f} = \alpha_{27}\beta_{16}\gamma_4.\eqno{(078)}
$$

079. We have $(a+c+d)+(e+f)+(-a-c-d-e-f) = 0,$ therefore (see formula (078))
$$
N_{a+c+d,e+f} =  N_{e+f,-a-c-d-e-f} = N_{-a-c-d-e-f,a+c+d} = \alpha_{27}\beta_{16}\gamma_4.\eqno{(079)}
$$

080. We have $(a+c)+(d+e+f)+(-a)+(-c-d-e-f) = 0$ and $(d+e+f)+(-a)\notin E_6,$ hence 
$$
S = N_{a+c,d+e+f}N_{-a,-c-d-e-f} + N_{-a,a+c}N_{d+e+f,-c-d-e-f} = 0.
$$
Since
$$
N_{a+c,d+e+f}N_{-a,-c-d-e-f} = | \text{Table 2} | = -\alpha_7 N_{a+c,d+e+f}, 
$$
$$
N_{-a,a+c}N_{d+e+f,-c-d-e-f} = | \text{Formulas (001) and (070)} | = \alpha_1(\beta_6) = \alpha_1\beta_6,
$$
we have
$$
S = -\alpha_7 N_{a+c,d+e+f} + \alpha_1\beta_6 = 0
$$
and therefore
$$
N_{a+c,d+e+f} = \alpha_{17}\beta_6.\eqno{(080)}
$$

081. We have $(a+c)+(d+e+f)+(-a-c-d-e-f) = 0,$ therefore (see formula (080))
$$
N_{a+c,d+e+f} =  N_{d+e+f,-a-c-d-e-f} = N_{-a-c-d-e-f,a+c} = \alpha_{17}\beta_6.\eqno{(081)}
$$

082. We have $b+(d+e+f)+(-b-d-e-f) = 0,$ therefore (see Table 2) 
$$
N_{b,d+e+f} =  N_{d+e+f,-b-d-e-f} = N_{-b-d-e-f,b} = \delta_2.\eqno{(082)}
$$

083. We have $(b+d+e)+f+(-b)+(-d-e-f) = 0$ and $f+(-b)\notin E_6,$ hence 
$$
S = N_{b+d+e,f}N_{-b,-d-e-f} + N_{-b,b+d+e}N_{f,-d-e-f} = 0.
$$
Since
$$
N_{b+d+e,f}N_{-b,-d-e-f} = | \text{Table 2} | = -\delta_2 N_{b+d+e,f}, 
$$
$$
N_{-b,b+d+e}N_{f,-d-e-f} = | \text{Formulas (024) and (069)} | = \delta_1(\gamma_{24}\epsilon_1) = \gamma_{24}\delta_1\epsilon_1,
$$
we have
$$
S = -\delta_2 N_{b+d+e,f} + \gamma_{24}\delta_1\epsilon_1 = 0
$$
and therefore
$$
N_{b+d+e,f} =\gamma_{24}\delta_{12}\epsilon_1.\eqno{(083)}
$$

084. We have $(b+d+e)+f+(-b-d-e-f) = 0,$ therefore (see formula (083))
$$
N_{b+d+e,f} =  N_{f,-b-d-e-f} = N_{-b-d-e-f,b+d+e} = \gamma_{24}\delta_{12}\epsilon_1.\eqno{(084)}
$$

085. We have $(b+d)+(e+f)+(-b)+(-d-e-f) = 0$ and $(e+f)+(-b)\notin E_6,$ hence 
$$
S = N_{b+d,e+f}N_{-b,-d-e-f} + N_{-b,b+d}N_{e+f,-d-e-f} = 0.
$$
Since
$$
N_{b+d,e+f}N_{-b,-d-e-f} = | \text{Table 2} | = -\delta_2 N_{b+d,e+f}, 
$$
$$
N_{-b,b+d}N_{e+f,-d-e-f} = | \text{Formulas (006) and (067)} | = -\gamma_1(\gamma_4) = -\gamma_{14},
$$
we have
$$
S = -\delta_2 N_{b+d,e+f} - \gamma_{14} = 0
$$
and therefore
$$
N_{b+d,e+f} = -\gamma_{14}\delta_2.\eqno{(085)}
$$

086. We have $(b+d)+(e+f)+(-b-d-e-f) = 0,$ therefore (see formula (085))
$$
N_{b+d,e+f} =  N_{e+f,-b-d-e-f} = N_{-b-d-e-f,b+d} = -\gamma_{14}\delta_2.\eqno{(086)}
$$

087. We have $c+(b+d+e+f)+(-b-c-d-e-f) = 0,$ therefore (see Table 2)
$$
N_{c,b+d+e+f} =  N_{b+d+e+f,-b-c-d-e-f} = N_{-b-c-d-e-f,c} = \beta_7.\eqno{(087)}
$$

088. We have $(b+c+d+e)+f+(-c)+(-b-d-e-f) = 0$ and $f+(-c)\notin E_6,$ hence 
$$
S = N_{b+c+d+e,f}N_{-c,-b-d-e-f} + N_{-c,b+c+d+e}N_{f,-b-d-e-f} = 0.
$$
Since
$$
N_{b+c+d+e,f}N_{-c,-b-d-e-f} = | \text{Table 2} | = -\beta_7 N_{b+c+d+e,f}, 
$$
$$
N_{-c,b+c+d+e}N_{f,-b-d-e-f} = | \text{Formulas (027) and (084)}  | = \beta_4(\gamma_{24}\delta_{12}\epsilon_1) = \beta_4\gamma_{24}\delta_{12}\epsilon_1,
$$
we have
$$
S = -\beta_7 N_{b+c+d+e,f} + \beta_4\gamma_{24}\delta_{12}\epsilon_1 = 0
$$
and therefore
$$
N_{b+c+d+e,f} = \beta_{47}\gamma_{24}\delta_{12}\epsilon_1.\eqno{(088)}
$$

089. We have $(b+c+d+e)+f+(-b-c-d-e-f) = 0,$ therefore (see formula (088))
$$
N_{b+c+d+e,f} =  N_{f,-b-c-d-e-f} = N_{-b-c-d-e-f,b+c+d+e} = \beta_{47}\gamma_{24}\delta_{12}\epsilon_1.\eqno{(089)}
$$

090. We have $b+(c+d+e+f)+(-c)+(-b-d-e-f) = 0$ and $b+(-c)\notin E_6,$ hence 
$$
S = N_{b,c+d+e+f}N_{-c,-b-d-e-f} + N_{c+d+e+f,-c}N_{b,-b-d-e-f} = 0.
$$
Since
$$
N_{b,c+d+e+f}N_{-c,-b-d-e-f} = | \text{Table 2} | = -\beta_7 N_{b,c+d+e+f}, 
$$
$$
N_{c+d+e+f,-c}N_{b,-b-d-e-f} = | \text{Formulas (070) and (082)} | = -\beta_6(-\delta_2) = \beta_6\delta_2,
$$
we have
$$
S = -\beta_7 N_{b,c+d+e+f} + \beta_6\delta_2 = 0
$$
and therefore
$$
N_{b,c+d+e+f} = \beta_{67}\delta_2.\eqno{(090)}
$$

091. We have $b+(c+d+e+f)+(-b-c-d-e-f) = 0,$ therefore (see formula (090)) 
$$
N_{b,c+d+e+f} =  N_{c+d+e+f,-b-c-d-e-f} = N_{-b-c-d-e-f,b}=\beta_{67}\delta_2.\eqno{(091)}
$$

092. We have $(b+c+d)+(e+f)+(-c)+(-b-d-e-f) = 0$ and $(e+f)+(-c)\notin E_6,$ hence 
$$
S = N_{b+c+d,e+f}N_{-c,-b-d-e-f} + N_{-c,b+c+d}N_{e+f,-b-d-e-f} = 0.
$$
Since
$$
N_{b+c+d,e+f}N_{-c,-b-d-e-f} = | \text{Table 2} | = -\beta_7 N_{b+c+d,e+f},
$$
$$
N_{-c,b+c+d}N_{e+f,-b-d-e-f} = | \text{Formulas (007) and (086)} | = \beta_2(-\gamma_{14}\delta_2) = -\beta_2\gamma_{14}\delta_2,
$$
we have
$$
S = -\beta_7 N_{b+c+d,e+f} - \beta_2\gamma_{14}\delta_2 = 0
$$
and therefore
$$
N_{b+c+d,e+f} = -\beta_{27}\gamma_{14}\delta_2.\eqno{(092)}
$$

093. We have $(b+c+d)+(e+f)+(-b-c-d-e-f) = 0,$ therefore (see formula (092))
$$
N_{b+c+d,e+f} = N_{e+f,-b-c-d-e-f} = N_{-b-c-d-e-f,b+c+d} = -\beta_{27}\gamma_{14}\delta_2.\eqno{(093)}
$$

094. We have $a+(b+c+d+e+f)+(-a-b-c-d-e-f) = 0,$ therefore (see Table 2)
$$
N_{a,b+c+d+e+f} = N_{b+c+d+e+f,-a-b-c-d-e-f} = N_{-a-b-c-d-e-f,a} = \alpha_8.\eqno{(094)}
$$

095. We have $(a+b+c+d+e)+f+(-a)+(-b-c-d-e-f) = 0$ and $f+(-a)\notin E_6,$ hence 
$$
S = N_{a+b+c+d+e,f}N_{-a,-b-c-d-e-f} + N_{-a,a+b+c+d+e}N_{f,-b-c-d-e-f} = 0.
$$
Since
$$
N_{a+b+c+d+e,f}N_{-a,-b-c-d-e-f} = | \text{Table 2} | = -\alpha_8 N_{a+b+c+d+e,f},
$$
$$
N_{-a,a+b+c+d+e}N_{f,-b-c-d-e-f} = | \text{Formulas (032) and (089)} | = \alpha_5(\beta_{47}\gamma_{24}\delta_{12}\epsilon_1) = 
$$
$$
= \alpha_5\beta_{47}\gamma_{24}\delta_{12}\epsilon_1,
$$
we have
$$
S = -\alpha_8 N_{a+b+c+d+e,f} + \alpha_5\beta_{47}\gamma_{24}\delta_{12}\epsilon_1 = 0
$$
and therefore
$$
N_{a+b+c+d+e,f} = \alpha_{58}\beta_{47}\gamma_{24}\delta_{12}\epsilon_1.\eqno{(095)}
$$

096. We have $(a+b+c+d+e)+f+(-a-b-c-d-e-f) = 0,$ therefore (see formula (095))
$$
N_{a+b+c+d+e,f} = N_{f,-a-b-c-d-e-f} = N_{-a-b-c-d-e-f,a+b+c+d+e} = \alpha_{58}\beta_{47}\gamma_{24}\delta_{12}\epsilon_1.\eqno{(096)}
$$

097. We have $(a+b+c+d)+(e+f)+(-a)+(-b-c-d-e-f) = 0$ and $(e+f)+(-a)\notin E_6,$ hence 
$$
S = N_{a+b+c+d,e+f}N_{-a,-b-c-d-e-f} + N_{-a,a+b+c+d}N_{e+f,-b-c-d-e-f} = 0.
$$
Since
$$
N_{a+b+c+d,e+f}N_{-a,-b-c-d-e-f} = | \text{Table 2} | = -\alpha_8 N_{a+b+c+d,e+f},
$$
$$
N_{-a,a+b+c+d}N_{e+f,-b-c-d-e-f} = | \text{Formulas (010) and (093)} | = \alpha_3(-\beta_{27}\gamma_{14}\delta_2) = -\alpha_3\beta_{27}\gamma_{14}\delta_2,
$$
we have
$$
S = -\alpha_8 N_{a+b+c+d,e+f} - \alpha_3\beta_{27}\gamma_{14}\delta_2 = 0
$$
and therefore
$$
N_{a+b+c+d,e+f} = -\alpha_{38}\beta_{27}\gamma_{14}\delta_2.\eqno{(097)}
$$

098. We have $(a+b+c+d)+(e+f)+(-a-b-c-d-e-f) = 0,$ therefore (see formula (097)) 
$$
N_{a+b+c+d,e+f} = N_{e+f,-a-b-c-d-e-f} = N_{-a-b-c-d-e-f,a+b+c+d} = -\alpha_{38}\beta_{27}\gamma_{14}\delta_2.\eqno{(098)}
$$

099. We have $b+(a+c+d+e+f)+(-a)+(-b-c-d-e-f) = 0$ and $b+(-a)\notin E_6,$ hence 
$$
S = N_{b,a+c+d+e+f}N_{-a,-b-c-d-e-f} + N_{a+c+d+e+f,-a}N_{b,-b-c-d-e-f} = 0.
$$
Since
$$
N_{b,a+c+d+e+f}N_{-a,-b-c-d-e-f} = | \text{Table 2} | = -\alpha_8 N_{b,a+c+d+e+f},
$$
$$
N_{a+c+d+e+f,-a}N_{b,-b-c-d-e-f} = | \text{Formulas (075) and (091)} | = -\alpha_7(-\beta_{67}\delta_2) = \alpha_7\beta_{67}\delta_2,
$$
we have
$$
S = -\alpha_8 N_{b,a+c+d+e+f} + \alpha_7\beta_{67}\delta_2 = 0
$$
and therefore
$$
N_{b,a+c+d+e+f} = \alpha_{78}\beta_{67}\delta_2.\eqno{(099)}
$$

100. We have $b+(a+c+d+e+f)+(-a-b-c-d-e-f) = 0,$ therefore (see formula (099))
$$
N_{b,a+c+d+e+f} = N_{a+c+d+e+f,-a-b-c-d-e-f} = N_{-a-b-c-d-e-f,b} = \alpha_{78}\beta_{67}\delta_2.\eqno{(100)}
$$

101. We have $(a+c)+(b+d+e+f)+(-a)+(-b-c-d-e-f) = 0$ and $(b+d+e+f)+(-a)\notin E_6,$ hence 
$$
S = N_{a+c,b+d+e+f}N_{-a,-b-c-d-e-f} + N_{-a,a+c}N_{b+d+e+f,-b-c-d-e-f} = 0.
$$
Since
$$
N_{a+c,b+d+e+f}N_{-a,-b-c-d-e-f} = | \text{Table 2} | = -\alpha_8 N_{a+c,b+d+e+f},
$$
$$
N_{-a,a+c}N_{b+d+e+f,-b-c-d-e-f} = | \text{Formulas (001) and (087)} | = \alpha_1(\beta_7) = \alpha_1\beta_7,
$$
we have
$$
S = -\alpha_8 N_{a+c,b+d+e+f} + \alpha_1\beta_7 = 0
$$
and therefore
$$
N_{a+c,b+d+e+f} = \alpha_{18}\beta_7.\eqno{(101)}
$$

102. We have $(a+c)+(b+d+e+f)+(-a-b-c-d-e-f) = 0,$ therefore (see formula (101))
$$
N_{a+c,b+d+e+f} = N_{b+d+e+f,-a-b-c-d-e-f} = N_{-a-b-c-d-e-f,a+c} = \alpha_{18}\beta_7.\eqno{(102)}
$$

103. We have $d+(b+c+d+e+f)+(-b-c-2d-e-f) = 0,$ therefore (see Table 2)
$$
N_{d,b+c+d+e+f} = N_{b+c+d+e+f,-b-c-2d-e-f} = N_{-b-c-2d-e-f,d} = \gamma_5.\eqno{(103)}
$$

104. We have $(b+c+2d+e)+f+(-d)+(-b-c-d-e-f) = 0$ and $f+(-d)\notin E_6,$ hence 
$$
S = N_{b+c+2d+e,f}N_{-d,-b-c-d-e-f} + N_{-d,b+c+2d+e}N_{f,-b-c-d-e-f} = 0.
$$
Since
$$
N_{b+c+2d+e,f}N_{-d,-b-c-d-e-f} = | \text{Table 2} | = -\gamma_5 N_{b+c+2d+e,f},
$$
$$
N_{-d,b+c+2d+e}N_{f,-b-c-d-e-f} = | \text{Formulas (039) and (089)} | = \gamma_3(\beta_{47}\gamma_{24}\delta_{12}\epsilon_1) = \beta_{47}\gamma_{234}\delta_{12}\epsilon_1,
$$
we have
$$
S = -\gamma_5 N_{b+c+2d+e,f} + \beta_{47}\gamma_{234}\delta_{12}\epsilon_1 = 0
$$
and therefore
$$
N_{b+c+2d+e,f} = \beta_{47}\gamma_{2345}\delta_{12}\epsilon_1.\eqno{(104)}
$$

105. We have $(b+c+2d+e)+f+(-b-c-2d-e-f) = 0,$ therefore (see formula (104))
$$
N_{b+c+2d+e,f} = N_{f,-b-c-2d-e-f} = N_{-b-c-2d-e-f,b+c+2d+e}=\beta_{47}\gamma_{2345}\delta_{12}\epsilon_1.\eqno{(105)}
$$

106. We have $(b+d)+(c+d+e+f)+(-d)+(-b-c-d-e-f) = 0$ and $(c+d+e+f)+(-d)\notin E_6,$ hence 
$$
S = N_{b+d,c+d+e+f}N_{-d,-b-c-d-e-f} + N_{-d,b+d}N_{c+d+e+f,-b-c-d-e-f} = 0.
$$
Since
$$
N_{b+d,c+d+e+f}N_{-d,-b-c-d-e-f} = | \text{Table 2} | = -\gamma_5 N_{b+d,c+d+e+f},
$$
$$
N_{-d,b+d}N_{c+d+e+f,-b-c-d-e-f} = | \text{Formulas (006) and (091)} | = \gamma_1(\beta_{67}\delta_2) = \beta_{67}\gamma_1\delta_2,
$$
we have
$$
S = -\gamma_5 N_{b+d,c+d+e+f} + \beta_{67}\gamma_1\delta_2 = 0
$$
and therefore
$$
N_{b+d,c+d+e+f} = \beta_{67}\gamma_{15}\delta_2.\eqno{(106)}
$$

107. We have $(b+d)+(c+d+e+f)+(-b-c-2d-e-f) = 0,$ therefore (see formula (106))
$$
N_{b+d,c+d+e+f} = N_{c+d+e+f,-b-c-2d-e-f} = N_{-b-c-2d-e-f,b+d} = \beta_{67}\gamma_{15}\delta_2.\eqno{(107)}
$$

108. We have $(c+d)+(b+d+e+f)+(-d)+(-b-c-d-e-f) = 0$ and $(b+d+e+f)+(-d)\notin E_6,$ hence 
$$
S = N_{c+d,b+d+e+f}N_{-d,-b-c-d-e-f} + N_{-d,c+d}N_{b+d+e+f,-b-c-d-e-f} = 0.
$$
Since
$$
N_{c+d,b+d+e+f}N_{-d,-b-c-d-e-f} = | \text{Table 2} | = -\gamma_5 N_{c+d,b+d+e+f},
$$
$$
N_{-d,c+d}N_{b+d+e+f,-b-c-d-e-f} = | \text{Formulas (002) and (087)} | = -\beta_1(\beta_7) = -\beta_{17},
$$
we have
$$
S = -\gamma_5 N_{c+d,b+d+e+f} - \beta_{17} = 0
$$
and therefore
$$
N_{c+d,b+d+e+f} = -\beta_{17}\gamma_5.\eqno{(108)}
$$

109. We have $(c+d)+(b+d+e+f)+(-b-c-2d-e-f) = 0,$ therefore (see formula (108)) 
$$
N_{c+d,b+d+e+f} = N_{b+d+e+f,-b-c-2d-e-f} = N_{-b-c-2d-e-f,c+d} = -\beta_{17}\gamma_5.\eqno{(109)}
$$

110. We have $(b+c+d)+(d+e+f)+(-d)+(-b-c-d-e-f) = 0$ and $(b+c+d)+(-d)\notin E_6,$ hence 
$$
S = N_{b+c+d,d+e+f}N_{-d,-b-c-d-e-f} + N_{d+e+f,-d}N_{b+c+d,-b-c-d-e-f} = 0.
$$
Since
$$
N_{b+c+d,d+e+f}N_{-d,-b-c-d-e-f} = | \text{Table 2} | = -\gamma_5 N_{b+c+d,d+e+f},
$$
$$
N_{d+e+f,-d}N_{b+c+d,-b-c-d-e-f} = | \text{Formulas (067) and (093)} | = -\gamma_4(\beta_{27}\gamma_{14}\delta_2) = -\beta_{27}\gamma_1\delta_2,
$$
we have
$$
S = -\gamma_5 N_{b+c+d,d+e+f} - \beta_{27}\gamma_1\delta_2 = 0
$$
and therefore
$$
N_{b+c+d,d+e+f} = -\beta_{27}\gamma_{15}\delta_2.\eqno{(110)}
$$

111. We have $(b+c+d)+(d+e+f)+(-b-c-2d-e-f) = 0,$ therefore (see formula (110))
$$
N_{b+c+d,d+e+f} = N_{d+e+f,-b-c-2d-e-f} = N_{-b-c-2d-e-f,b+c+d} = -\beta_{27}\gamma_{15}\delta_2.\eqno{(111)}
$$

112. We have $a+(b+c+2d+e+f)+(-a-b-c-2d-e-f) = 0,$ therefore (see Table 2)
$$
N_{a,b+c+2d+e+f} = N_{b+c+2d+e+f,-a-b-c-2d-e-f} = N_{-a-b-c-2d-e-f,a} = \alpha_9.\eqno{(112)}
$$

113. We have $d+(a+b+c+d+e+f)+(-a)+(-b-c-2d-e-f) = 0$ and $d+(-a)\notin E_6,$ hence 
$$
S = N_{d,a+b+c+d+e+f}N_{-a,-b-c-2d-e-f} + N_{a+b+c+d+e+f,-a}N_{d,-b-c-2d-e-f} = 0.
$$
Since
$$
N_{d,a+b+c+d+e+f}N_{-a,-b-c-2d-e-f} = | \text{Table 2} | = -\alpha_9 N_{d,a+b+c+d+e+f},
$$
$$
N_{a+b+c+d+e+f,-a}N_{d,-b-c-2d-e-f} = | \text{Formulas (094) and (103)} | = -\alpha_8(-\gamma_5) = \alpha_8\gamma_5,
$$
we have
$$
S = -\alpha_9 N_{d,a+b+c+d+e+f} + \alpha_8\gamma_5 = 0
$$
and therefore
$$
N_{d,a+b+c+d+e+f} = \alpha_{89}\gamma_5.\eqno{(113)}
$$

114. We have $d+(a+b+c+d+e+f)+(-a-b-c-2d-e-f) = 0,$ therefore (see formula (113))
$$
N_{d,a+b+c+d+e+f} = N_{a+b+c+d+e+f,-a-b-c-2d-e-f} = N_{-a-b-c-2d-e-f,d} = \alpha_{89}\gamma_5.\eqno{(114)}
$$

115. We have $(b+d)+(a+c+d+e+f)+(-a)+(-b-c-2d-e-f) = 0$ and $(b+d)+(-a)\notin E_6,$ hence 
$$
S = N_{b+d,a+c+d+e+f}N_{-a,-b-c-2d-e-f} + N_{a+c+d+e+f,-a}N_{b+d,-b-c-2d-e-f} = 0.
$$
Since
$$
N_{b+d,a+c+d+e+f}N_{-a,-b-c-2d-e-f} = | \text{Table 2} | = -\alpha_9 N_{b+d,a+c+d+e+f},
$$
$$
N_{a+c+d+e+f,-a}N_{b+d,-b-c-2d-e-f} = | \text{Formulas (075) and (107)} | = -\alpha_7(-\beta_{67}\gamma_{15}\delta_2) =
$$
$$
= \alpha_7\beta_{67}\gamma_{15}\delta_2,
$$
we have
$$
S = -\alpha_9 N_{b+d,a+c+d+e+f} + \alpha_7\beta_{67}\gamma_{15}\delta_2 = 0
$$
and therefore
$$
N_{b+d,a+c+d+e+f} = \alpha_{79}\beta_{67}\gamma_{15}\delta_2.\eqno{(115)}
$$

116. We have $(b+d)+(a+c+d+e+f)+(-a-b-c-2d-e-f) = 0,$ therefore (see formula (115))
$$
N_{b+d,a+c+d+e+f} = N_{a+c+d+e+f,-a-b-c-2d-e-f} = N_{-a-b-c-2d-e-f,b+d} = \alpha_{79}\beta_{67}\gamma_{15}\delta_2.\eqno{(116)}
$$

117. We have $(a+b+c+d)+(d+e+f)+(-a)+(-b-c-2d-e-f) = 0$ and $(d+e+f)+(-a)\notin E_6,$ hence 
$$
S = N_{a+b+c+d,d+e+f}N_{-a,-b-c-2d-e-f} + N_{-a,a+b+c+d}N_{d+e+f,-b-c-2d-e-f} = 0.
$$
Since
$$
N_{a+b+c+d,d+e+f}N_{-a,-b-c-2d-e-f} = | \text{Table 2} | = -\alpha_9 N_{a+b+c+d,d+e+f},
$$
$$
N_{-a,a+b+c+d}N_{d+e+f,-b-c-2d-e-f} = | \text{Formulas (010) and (111)} | = \alpha_3(-\beta_{27}\gamma_{15}\delta_2) = 
$$
$$
=-\alpha_3\beta_{27}\gamma_{15}\delta_2,
$$
we have
$$
S = -\alpha_9 N_{a+b+c+d,d+e+f} - \alpha_3\beta_{27}\gamma_{15}\delta_2 = 0
$$
and therefore
$$
N_{a+b+c+d,d+e+f} = -\alpha_{39}\beta_{27}\gamma_{15}\delta_2.\eqno{(117)}
$$

118. We have $(a+b+c+d)+(d+e+f)+(-a-b-c-2d-e-f) = 0,$ therefore (see formula (117))
$$
N_{a+b+c+d,d+e+f} = N_{d+e+f,-a-b-c-2d-e-f} = N_{-a-b-c-2d-e-f,a+b+c+d} = -\alpha_{39}\beta_{27}\gamma_{15}\delta_2.\eqno{(118)}
$$

119. We have $(a+c+d)+(b+d+e+f)+(-a)+(-b-c-2d-e-f) = 0$ and $(b+d+e+f)+(-a)\notin E_6,$ hence 
$$
S = N_{a+c+d,b+d+e+f}N_{-a,-b-c-2d-e-f} + N_{-a,a+c+d}N_{b+d+e+f,-b-c-2d-e-f} = 0.
$$
Since
$$
N_{a+c+d,b+d+e+f}N_{-a,-b-c-2d-e-f} = | \text{Table 2} | = -\alpha_9 N_{a+c+d,b+d+e+f},
$$
$$
N_{-a,a+c+d}N_{b+d+e+f,-b-c-2d-e-f} = | \text{Formulas (003) and (109)} | = \alpha_2(-\beta_{17}\gamma_5) = -\alpha_2\beta_{17}\gamma_5,
$$
we have
$$
S = -\alpha_9 N_{a+c+d,b+d+e+f} - \alpha_2\beta_{17}\gamma_5 = 0
$$
and therefore
$$
N_{a+c+d,b+d+e+f} = -\alpha_{29}\beta_{17}\gamma_5.\eqno{(119)}
$$

120. We have $(a+c+d)+(b+d+e+f)+(-a-b-c-2d-e-f) = 0,$ therefore (see formula (119)) 
$$
N_{a+c+d,b+d+e+f} = N_{b+d+e+f,-a-b-c-2d-e-f} = N_{-a-b-c-2d-e-f,a+c+d} = -\alpha_{29}\beta_{17}\gamma_5.\eqno{(120)}
$$

121. We have $(a+b+c+2d+e)+f+(-a)+(-b-c-2d-e-f) = 0$ and $f+(-a)\notin E_6,$ hence 
$$
S = N_{a+b+c+2d+e,f}N_{-a,-b-c-2d-e-f} + N_{-a,a+b+c+2d+e}N_{f,-b-c-2d-e-f} = 0.
$$
Since
$$
N_{a+b+c+2d+e,f}N_{-a,-b-c-2d-e-f} = | \text{Table 2} | = -\alpha_9 N_{a+b+c+2d+e,f},
$$
$$
N_{-a,a+b+c+2d+e}N_{f,-b-c-2d-e-f} = | \text{Formulas (046) and (105)} | = \alpha_6(\beta_{47}\gamma_{2345}\delta_{12}\epsilon_1) = 
$$
$$
=\alpha_6\beta_{47}\gamma_{2345}\delta_{12}\epsilon_1,
$$
we have
$$
S = -\alpha_9 N_{a+b+c+2d+e,f} + \alpha_6\beta_{47}\gamma_{2345}\delta_{12}\epsilon_1 = 0
$$
and therefore
$$
N_{a+b+c+2d+e,f} = \alpha_{69}\beta_{47}\gamma_{2345}\delta_{12}\epsilon_1.\eqno{(121)}
$$

122. We have $(a+b+c+2d+e)+f+(-a-b-c-2d-e-f) = 0,$ therefore (see formula (121))
$$
N_{a+b+c+2d+e,f} = N_{f,-a-b-c-2d-e-f} = N_{-a-b-c-2d-e-f,a+b+c+2d+e} = \alpha_{69}\beta_{47}\gamma_{2345}\delta_{12}\epsilon_1.\eqno{(122)}
$$

123. We have $c+(a+b+c+2d+e+f)+(-a-b-2c-2d-e-f) = 0,$ therefore (see Table 2)
$$
N_{c,a+b+c+2d+e+f} = N_{a+b+c+2d+e+f,-a-b-2c-2d-e-f} = N_{-a-b-2c-2d-e-f,c} = \beta_8.\eqno{(123)}
$$

124. We have $(a+b+2c+2d+e)+f+(-c)+(-a-b-c-2d-e-f) = 0$ and $f+(-c)\notin E_6,$ hence 
$$
S = N_{a+b+2c+2d+e,f}N_{-c,-a-b-c-2d-e-f} + N_{-c,a+b+2c+2d+e}N_{f,-a-b-c-2d-e-f} = 0.
$$
Since
$$
N_{a+b+2c+2d+e,f}N_{-c,-a-b-c-2d-e-f} = | \text{Table 2} | = -\beta_8 N_{a+b+2c+2d+e,f},
$$
$$
N_{-c,a+b+2c+2d+e}N_{f,-a-b-c-2d-e-f} = | \text{Formulas (055) and (122)} |  =
$$
$$
= \beta_5(\alpha_{69}\beta_{47}\gamma_{2345}\delta_{12}\epsilon_1) = \alpha_{69}\beta_{457}\gamma_{2345}\delta_{12}\epsilon_1,
$$
we have
$$
S = -\beta_8 N_{a+b+2c+2d+e,f} + \alpha_{69}\beta_{457}\gamma_{2345}\delta_{12}\epsilon_1 = 0
$$
and therefore
$$
N_{a+b+2c+2d+e,f} = \alpha_{69}\beta_{4578}\gamma_{2345}\delta_{12}\epsilon_1.\eqno{(124)}
$$

125. We have $(a+b+2c+2d+e)+f+(-a-b-2c-2d-e-f) = 0,$ therefore (see formula (124))
$$
N_{a+b+2c+2d+e,f} = N_{f,-a-b-2c-2d-e-f} = N_{-a-b-2c-2d-e-f,a+b+2c+2d+e} = \alpha_{69}\beta_{4578}\gamma_{2345}\delta_{12}\epsilon_1.\eqno{(125)}
$$

126. We have $(c+d)+(a+b+c+d+e+f)+(-c)+(-a-b-c-2d-e-f) = 0$ and $(a+b+c+d+e+f)+(-c)\notin E_6,$ hence 
$$
S = N_{c+d,a+b+c+d+e+f}N_{-c,-a-b-c-2d-e-f} + N_{-c,c+d}N_{a+b+c+d+e+f,-a-b-c-2d-e-f} = 0.
$$
Since
$$
N_{c+d,a+b+c+d+e+f}N_{-c,-a-b-c-2d-e-f} = | \text{Table 2} | = -\beta_8 N_{c+d,a+b+c+d+e+f},
$$
$$
N_{-c,c+d}N_{a+b+c+d+e+f,-a-b-c-2d-e-f} = | \text{Formulas (002) and (114)} | = \beta_1(\alpha_{89}\gamma_5) = \alpha_{89}\beta_1\gamma_5,
$$
we have
$$
S = -\beta_8 N_{c+d,a+b+c+d+e+f} + \alpha_{89}\beta_1\gamma_5 = 0
$$
and therefore
$$
N_{c+d,a+b+c+d+e+f} = \alpha_{89}\beta_{18}\gamma_5.\eqno{(126)}
$$

127. We have $(c+d)+(a+b+c+d+e+f)+(-a-b-2c-2d-e-f) = 0,$ therefore (see formula (126)) 
$$
N_{c+d,a+b+c+d+e+f} = N_{a+b+c+d+e+f,-a-b-2c-2d-e-f} = N_{-a-b-2c-2d-e-f,c+d} = \alpha_{89}\beta_{18}\gamma_5.\eqno{(127)}
$$

128. We have $(a+c)+(b+c+2d+e+f)+(-c)+(-a-b-c-2d-e-f) = 0$ and $(b+c+2d+e+f)+(-c)\notin E_6,$ hence 
$$
S = N_{a+c,b+c+2d+e+f}N_{-c,-a-b-c-2d-e-f} + N_{-c,a+c}N_{b+c+2d+e+f,-a-b-c-2d-e-f} = 0.
$$
Since
$$
N_{a+c,b+c+2d+e+f}N_{-c,-a-b-c-2d-e-f} = | \text{Table 2} | = -\beta_8 N_{a+c,b+c+2d+e+f},
$$
$$
N_{-c,a+c}N_{b+c+2d+e+f,-a-b-c-2d-e-f} = | \text{Formulas (001) and (112)} | = -\alpha_1(\alpha_9) = -\alpha_{19},
$$
we have
$$
S = -\beta_8 N_{a+c,b+c+2d+e+f} - \alpha_{19} = 0
$$
and therefore
$$
N_{a+c,b+c+2d+e+f} = -\alpha_{19}\beta_8.\eqno{(128)}
$$

129. We have $(a+c)+(b+c+2d+e+f)+(-a-b-2c-2d-e-f) = 0,$ therefore (see formula (128))
$$
N_{a+c,b+c+2d+e+f} = N_{b+c+2d+e+f,-a-b-2c-2d-e-f} = N_{-a-b-2c-2d-e-f,a+c} = -\alpha_{19}\beta_8.\eqno{(129)}
$$

130. We have $(a+c+d)+(b+c+d+e+f)+(-c)+(-a-b-c-2d-e-f) = 0$ and $(a+c+d)+(-c)\notin E_6,$ hence 
$$
S = N_{a+c+d,b+c+d+e+f}N_{-c,-a-b-c-2d-e-f} + N_{b+c+d+e+f,-c}N_{a+c+d,-a-b-c-2d-e-f} = 0.
$$
Since
$$
N_{a+c+d,b+c+d+e+f}N_{-c,-a-b-c-2d-e-f} = | \text{Table 2} | = -\beta_8 N_{a+c+d,b+c+d+e+f},
$$
$$
N_{b+c+d+e+f,-c}N_{a+c+d,-a-b-c-2d-e-f} = | \text{Formulas (087) and (120)} | = -\beta_7(\alpha_{29}\beta_{17}\gamma_5) = 
$$
$$
=-\alpha_{29}\beta_1\gamma_5,
$$
we have
$$
S = -\beta_8 N_{a+c+d,b+c+d+e+f} - \alpha_{29}\beta_1\gamma_5 = 0
$$
and therefore
$$
N_{a+c+d,b+c+d+e+f} = -\alpha_{29}\beta_{18}\gamma_5.\eqno{(130)}
$$

131. We have $(a+c+d)+(b+c+d+e+f)+(-a-b-2c-2d-e-f) = 0,$ therefore (see formula (130))
$$
N_{a+c+d,b+c+d+e+f} = N_{b+c+d+e+f,-a-b-2c-2d-e-f} = N_{-a-b-2c-2d-e-f,a+c+d} = -\alpha_{29}\beta_{18}\gamma_5.\eqno{(131)}
$$

132. We have $(b+c+d)+(a+c+d+e+f)+(-c)+(-a-b-c-2d-e-f) = 0$ and $(a+c+d+e+f)+(-c)\notin E_6,$ hence 
$$
S = N_{b+c+d,a+c+d+e+f}N_{-c,-a-b-c-2d-e-f} + N_{-c,b+c+d}N_{a+c+d+e+f,-a-b-c-2d-e-f} = 0.
$$
Since
$$
N_{b+c+d,a+c+d+e+f}N_{-c,-a-b-c-2d-e-f} = | \text{Table 2} | = -\beta_8 N_{b+c+d,a+c+d+e+f},
$$
$$
N_{-c,b+c+d}N_{a+c+d+e+f,-a-b-c-2d-e-f} = | \text{Formulas (007) and (116)} | = \beta_2(\alpha_{79}\beta_{67}\gamma_{15}\delta_2) = 
$$
$$
=\alpha_{79}\beta_{267}\gamma_{15}\delta_2,
$$
we have
$$
S = -\beta_8 N_{b+c+d,a+c+d+e+f} + \alpha_{79}\beta_{267}\gamma_{15}\delta_2 = 0
$$
and therefore
$$
N_{b+c+d,a+c+d+e+f} = \alpha_{79}\beta_{2678}\gamma_{15}\delta_2.\eqno{(132)}
$$

133. We have $(b+c+d)+(a+c+d+e+f)+(-a-b-2c-2d-e-f) = 0,$ therefore (see formula (132))
$$
N_{b+c+d,a+c+d+e+f} = N_{a+c+d+e+f,-a-b-2c-2d-e-f} = N_{-a-b-2c-2d-e-f,b+c+d} = \alpha_{79}\beta_{2678}\gamma_{15}\delta_2.\eqno{(133)}
$$

134. We have $(a+b+c+d)+(c+d+e+f)+(-c)+(-a-b-c-2d-e-f) = 0$ and $(a+b+c+d)+(-c)\notin E_6,$ hence 
$$
S = N_{a+b+c+d,c+d+e+f}N_{-c,-a-b-c-2d-e-f} + N_{c+d+e+f,-c}N_{a+b+c+d,-a-b-c-2d-e-f} = 0.
$$
Since
$$
N_{a+b+c+d,c+d+e+f}N_{-c,-a-b-c-2d-e-f} = | \text{Table 2} | = -\beta_8 N_{a+b+c+d,c+d+e+f},
$$
$$
N_{c+d+e+f,-c}N_{a+b+c+d,-a-b-c-2d-e-f} = | \text{Formulas (070) and (118)} | =  -\beta_6(\alpha_{39}\beta_{27}\gamma_{15}\delta_2) = 
$$
$$
=-\alpha_{39}\beta_{267}\gamma_{15}\delta_2,
$$
we have
$$
S = -\beta_8 N_{a+b+c+d,c+d+e+f} - \alpha_{39}\beta_{267}\gamma_{15}\delta_2 = 0
$$
and therefore
$$
N_{a+b+c+d,c+d+e+f} = -\alpha_{39}\beta_{2678}\gamma_{15}\delta_2.\eqno{(134)}
$$

135. We have $(a+b+c+d)+(c+d+e+f)+(-a-b-2c-2d-e-f) = 0,$ therefore (see formula (134))
$$
N_{a+b+c+d,c+d+e+f} = N_{c+d+e+f,-a-b-2c-2d-e-f} = N_{-a-b-2c-2d-e-f,a+b+c+d} = -\alpha_{39}\beta_{2678}\gamma_{15}\delta_2.\eqno{(135)}
$$

136. We have $e+(b+c+2d+e+f)+(-b-c-2d-2e-f) = 0,$ therefore (see Table 2)
$$
N_{e,b+c+2d+e+f} = N_{b+c+2d+e+f,-b-c-2d-2e-f} = N_{-b-c-2d-2e-f,e} = \epsilon_2.\eqno{(136)}
$$

137. We have $(d+e)+(b+c+d+e+f)+(-e)+(-b-c-2d-e-f) = 0$ and $(b+c+d+e+f)+(-e)\notin E_6,$ hence 
$$
S = N_{d+e,b+c+d+e+f}N_{-e,-b-c-2d-e-f} + N_{-e,d+e}N_{b+c+d+e+f,-b-c-2d-e-f} = 0.
$$
Since
$$
N_{d+e,b+c+d+e+f}N_{-e,-b-c-2d-e-f} = | \text{Table 2} | = -\epsilon_2 N_{d+e,b+c+d+e+f},
$$
$$
N_{-e,d+e}N_{b+c+d+e+f,-b-c-2d-e-f} = | \text{Formulas (015) and (103)} | = -\gamma_2(\gamma_5) = -\gamma_{25},
$$
we have
$$
S = -\epsilon_2 N_{d+e,b+c+d+e+f} - \gamma_{25} = 0
$$
and therefore
$$
N_{d+e,b+c+d+e+f} = -\gamma_{25}\epsilon_2.\eqno{(137)}
$$

138. We have $(d+e)+(b+c+d+e+f)+(-b-c-2d-2e-f) = 0,$ therefore (see formula (137))
$$
N_{d+e,b+c+d+e+f} = N_{b+c+d+e+f,-b-c-2d-2e-f} = N_{-b-c-2d-2e-f,d+e} = -\gamma_{25}\epsilon_2.\eqno{(138)}
$$

139. We have $(b+c+2d+e)+(e+f)+(-e)+(-b-c-2d-e-f) = 0$ and $(b+c+2d+e)+(-e)\notin E_6,$ hence 
$$
S = N_{b+c+2d+e,e+f}N_{-e,-b-c-2d-e-f} + N_{e+f,-e}N_{b+c+2d+e,-b-c-2d-e-f} = 0.
$$
Since
$$
N_{b+c+2d+e,e+f}N_{-e,-b-c-2d-e-f} = | \text{Table 2} | = -\epsilon_2 N_{b+c+2d+e,e+f},
$$
$$
N_{e+f,-e}N_{b+c+2d+e,-b-c-2d-e-f} = | \text{Formulas (066) and (105)} | = -\epsilon_1(-\beta_{47}\gamma_{2345}\delta_{12}\epsilon_1) = 
$$
$$
=\beta_{47}\gamma_{2345}\delta_{12},
$$
we have
$$
S = -\epsilon_2 N_{b+c+2d+e,e+f} + \beta_{47}\gamma_{2345}\delta_{12} = 0
$$
and therefore
$$
N_{b+c+2d+e,e+f} = \beta_{47}\gamma_{2345}\delta_{12}\epsilon_2.\eqno{(139)}
$$

140. We have $(b+c+2d+e)+(e+f)+(-b-c-2d-2e-f) = 0,$ therefore (see formula (139))
$$
N_{b+c+2d+e,e+f} = N_{e+f,-b-c-2d-2e-f} = N_{-b-c-2d-2e-f,b+c+2d+e} = \beta_{47}\gamma_{2345}\delta_{12}\epsilon_2.\eqno{(140)}
$$

141. We have $(b+c+d+e)+(d+e+f)+(-e)+(-b-c-2d-e-f) = 0$ and $(d+e+f)+(-e)\notin E_6,$ hence 
$$
S = N_{b+c+d+e,d+e+f}N_{-e,-b-c-2d-e-f} + N_{-e,b+c+d+e}N_{d+e+f,-b-c-2d-e-f} = 0.
$$
Since
$$
N_{b+c+d+e,d+e+f}N_{-e,-b-c-2d-e-f} = | \text{Table 2} | = -\epsilon_2 N_{b+c+d+e,d+e+f},
$$
$$
N_{-e,b+c+d+e}N_{d+e+f,-b-c-2d-e-f} = | \text{Formulas (029) and (111)} | =  \beta_{24}\gamma_{12}\delta_1(-\beta_{27}\gamma_{15}\delta_2) =
$$
$$
= -\beta_{47}\gamma_{25}\delta_{12},
$$
we have
$$
S = -\epsilon_2 N_{b+c+d+e,d+e+f} - \beta_{47}\gamma_{25}\delta_{12} = 0
$$
and therefore
$$
N_{b+c+d+e,d+e+f} = -\beta_{47}\gamma_{25}\delta_{12}\epsilon_2.\eqno{(141)}
$$

142. We have $(b+c+d+e)+(d+e+f)+(-b-c-2d-2e-f) = 0,$ therefore (see formula (141))
$$
N_{b+c+d+e,d+e+f} = N_{d+e+f,-b-c-2d-2e-f} = N_{-b-c-2d-2e-f,b+c+d+e} = -\beta_{47}\gamma_{25}\delta_{12}\epsilon_2.\eqno{(142)}
$$

143. We have $(b+d+e)+(c+d+e+f)+(-e)+(-b-c-2d-e-f) = 0$ and $(c+d+e+f)+(-e)\notin E_6,$ hence 
$$
S = N_{b+d+e,c+d+e+f}N_{-e,-b-c-2d-e-f} + N_{-e,b+d+e}N_{c+d+e+f,-b-c-2d-e-f} = 0.
$$
Since
$$
N_{b+d+e,c+d+e+f}N_{-e,-b-c-2d-e-f} = | \text{Table 2} | = -\epsilon_2 N_{b+d+e,c+d+e+f},
$$
$$
N_{-e,b+d+e}N_{c+d+e+f,-b-c-2d-e-f} = | \text{Formulas (026) and (107)} | = \gamma_{12}\delta_1(\beta_{67}\gamma_{15}\delta_2) = \beta_{67}\gamma_{25}\delta_{12},
$$
we have
$$
S = -\epsilon_2 N_{b+d+e,c+d+e+f} + \beta_{67}\gamma_{25}\delta_{12} = 0
$$
and therefore
$$
N_{b+d+e,c+d+e+f} = \beta_{67}\gamma_{25}\delta_{12}\epsilon_2.\eqno{(143)}
$$

144. We have $(b+d+e)+(c+d+e+f)+(-b-c-2d-2e-f) = 0,$ therefore (see formula (143))
$$
N_{b+d+e,c+d+e+f} = N_{c+d+e+f,-b-c-2d-2e-f} = N_{-b-c-2d-2e-f,b+d+e} = \beta_{67}\gamma_{25}\delta_{12}\epsilon_2.\eqno{(144)}
$$

145. We have $(c+d+e)+(b+d+e+f)+(-e)+(-b-c-2d-e-f) = 0$ and $(b+d+e+f)+(-e)\notin E_6,$ hence 
$$
S = N_{c+d+e,b+d+e+f}N_{-e,-b-c-2d-e-f} + N_{-e,c+d+e}N_{b+d+e+f,-b-c-2d-e-f} = 0.
$$
Since
$$
N_{c+d+e,b+d+e+f}N_{-e,-b-c-2d-e-f} = | \text{Table 2} | = -\epsilon_2 N_{c+d+e,b+d+e+f},
$$
$$
N_{-e,c+d+e}N_{b+d+e+f,-b-c-2d-e-f} = | \text{Formulas (018) and (109)} | = -\beta_{13}\gamma_2(-\beta_{17}\gamma_5) = \beta_{37}\gamma_{25},
$$
we have
$$
S = -\epsilon_2 N_{c+d+e,b+d+e+f} + \beta_{37}\gamma_{25} = 0
$$
and therefore
$$
N_{c+d+e,b+d+e+f} = \beta_{37}\gamma_{25}\epsilon_2.\eqno{(145)}
$$

146. We have $(c+d+e)+(b+d+e+f)+(-b-c-2d-2e-f) = 0,$ therefore (see formula (145))
$$
N_{c+d+e,b+d+e+f} = N_{b+d+e+f,-b-c-2d-2e-f} = N_{-b-c-2d-2e-f,c+d+e} = \beta_{37}\gamma_{25}\epsilon_2.\eqno{(146)}
$$

147. We have $a+(b+c+2d+2e+f)+(-a-b-c-2d-2e-f) = 0,$ therefore (see Table 2)
$$
N_{a,b+c+2d+2e+f} = N_{b+c+2d+2e+f,-a-b-c-2d-2e-f} = N_{-a-b-c-2d-2e-f,a} = \alpha_0.\eqno{(147)}
$$

148. We have $e+(a+b+c+2d+e+f)+(-a)+(-b-c-2d-2e-f) = 0$ and $e+(-a)\notin E_6,$ hence 
$$
S = N_{e,a+b+c+2d+e+f}N_{-a,-b-c-2d-2e-f} + N_{a+b+c+2d+e+f,-a}N_{e,-b-c-2d-2e-f} = 0.
$$
Since
$$
N_{e,a+b+c+2d+e+f}N_{-a,-b-c-2d-2e-f} = | \text{Table 2} | = -\alpha_0 N_{e,a+b+c+2d+e+f},
$$
$$
N_{a+b+c+2d+e+f,-a}N_{e,-b-c-2d-2e-f} = | \text{Formulas (112) and (136)} | = -\alpha_9(-\epsilon_2) = \alpha_9\epsilon_2,
$$
we have
$$
S = -\alpha_0 N_{e,a+b+c+2d+e+f} + \alpha_9\epsilon_2 = 0
$$
and therefore
$$
N_{e,a+b+c+2d+e+f} = \alpha_{90}\epsilon_2.\eqno{(148)}
$$

149. We have $e+(a+b+c+2d+e+f)+(-a-b-c-2d-2e-f) = 0,$ therefore (see formula (148))
$$
N_{e,a+b+c+2d+e+f} = N_{a+b+c+2d+e+f,-a-b-c-2d-2e-f} = N_{-a-b-c-2d-2e-f,e} = \alpha_{90}\epsilon_2.\eqno{(149)}
$$

150. We have $(d+e)+(a+b+c+d+e+f)+(-a)+(-b-c-2d-2e-f) = 0$ and $(d+e)+(-a)\notin E_6,$ hence 
$$
S = N_{d+e,a+b+c+d+e+f}N_{-a,-b-c-2d-2e-f} + N_{a+b+c+d+e+f,-a}N_{d+e,-b-c-2d-2e-f} = 0.
$$
Since
$$
N_{d+e,a+b+c+d+e+f}N_{-a,-b-c-2d-2e-f} = | \text{Table 2} | = -\alpha_0 N_{d+e,a+b+c+d+e+f},
$$
$$
N_{a+b+c+d+e+f,-a}N_{d+e,-b-c-2d-2e-f} = | \text{Formulas (094) and (138)} | = -\alpha_8(\gamma_{25}\epsilon_2) = -\alpha_8\gamma_{25}\epsilon_2,
$$
we have
$$
S = -\alpha_0 N_{d+e,a+b+c+d+e+f} - \alpha_8\gamma_{25}\epsilon_2 = 0
$$
and therefore
$$
N_{d+e,a+b+c+d+e+f} = -\alpha_{80}\gamma_{25}\epsilon_2.\eqno{(150)}
$$

151. We have $(d+e)+(a+b+c+d+e+f)+(-a-b-c-2d-2e-f) = 0,$ therefore (see formula (150))
$$
N_{d+e,a+b+c+d+e+f} = N_{a+b+c+d+e+f,-a-b-c-2d-2e-f} = N_{-a-b-c-2d-2e-f,d+e} = -\alpha_{80}\gamma_{25}\epsilon_2.\eqno{(151)}
$$

152. We have $(a+b+c+2d+e)+(e+f)+(-a)+(-b-c-2d-2e-f) = 0$ and $(e+f)+(-a)\notin E_6,$ hence 
$$
S = N_{a+b+c+2d+e,e+f}N_{-a,-b-c-2d-2e-f} + N_{-a,a+b+c+2d+e}N_{e+f,-b-c-2d-2e-f} = 0.
$$
Since
$$
N_{a+b+c+2d+e,e+f}N_{-a,-b-c-2d-2e-f} = | \text{Table 2} | = -\alpha_0 N_{a+b+c+2d+e,e+f},
$$
$$
N_{-a,a+b+c+2d+e}N_{e+f,-b-c-2d-2e-f} = | \text{Formulas (046) and (140)} | = \alpha_6(\beta_{47}\gamma_{2345}\delta_{12}\epsilon_2) = 
$$
$$
=\alpha_6\beta_{47}\gamma_{2345}\delta_{12}\epsilon_2,
$$
we have
$$
S = -\alpha_0 N_{a+b+c+2d+e,e+f} + \alpha_6\beta_{47}\gamma_{2345}\delta_{12}\epsilon_2 = 0
$$
and therefore
$$
N_{a+b+c+2d+e,e+f} = \alpha_{60}\beta_{47}\gamma_{2345}\delta_{12}\epsilon_2.\eqno{(152)}
$$

153. We have $(a+b+c+2d+e)+(e+f)+(-a-b-c-2d-2e-f) = 0,$ therefore (see formula (152))
$$
N_{a+b+c+2d+e,e+f} = N_{e+f,-a-b-c-2d-2e-f} = N_{-a-b-c-2d-2e-f,a+b+c+2d+e} = \alpha_{60}\beta_{47}\gamma_{2345}\delta_{12}\epsilon_2.\eqno{(153)}
$$

154. We have $(a+b+c+d+e)+(d+e+f)+(-a)+(-b-c-2d-2e-f) = 0$ and $(d+e+f)+(-a)\notin E_6,$ hence 
$$
S = N_{a+b+c+d+e,d+e+f}N_{-a,-b-c-2d-2e-f} + N_{-a,a+b+c+d+e}N_{d+e+f,-b-c-2d-2e-f} = 0.
$$
Since
$$
N_{a+b+c+d+e,d+e+f}N_{-a,-b-c-2d-2e-f} = | \text{Table 2} | = -\alpha_0 N_{a+b+c+d+e,d+e+f},
$$
$$
N_{-a,a+b+c+d+e}N_{d+e+f,-b-c-2d-2e-f} = | \text{Formulas (032) and (142)} | = \alpha_5(-\beta_{47}\gamma_{25}\delta_{12}\epsilon_2) = 
$$
$$
=-\alpha_5\beta_{47}\gamma_{25}\delta_{12}\epsilon_2,
$$
we have
$$
S = -\alpha_0 N_{a+b+c+d+e,d+e+f} - \alpha_5\beta_{47}\gamma_{25}\delta_{12}\epsilon_2 = 0
$$
and therefore
$$
N_{a+b+c+d+e,d+e+f} = -\alpha_{50}\beta_{47}\gamma_{25}\delta_{12}\epsilon_2.\eqno{(154)}
$$

155. We have $(a+b+c+d+e)+(d+e+f)+(-a-b-c-2d-2e-f) = 0,$ therefore (see formula (154))
$$
N_{a+b+c+d+e,d+e+f} = N_{d+e+f,-a-b-c-2d-2e-f} = N_{-a-b-c-2d-2e-f,a+b+c+d+e} = 
$$
$$
=-\alpha_{50}\beta_{47}\gamma_{25}\delta_{12}\epsilon_2.\eqno{(155)}
$$

156. We have $(b+d+e)+(a+c+d+e+f)+(-a)+(-b-c-2d-2e-f) = 0$ and $(b+d+e)+(-a)\notin E_6,$ hence 
$$
S = N_{b+d+e,a+c+d+e+f}N_{-a,-b-c-2d-2e-f} + N_{a+c+d+e+f,-a}N_{b+d+e,-b-c-2d-2e-f} = 0.
$$
Since
$$
N_{b+d+e,a+c+d+e+f}N_{-a,-b-c-2d-2e-f} = | \text{Table 2} | = -\alpha_0 N_{b+d+e,a+c+d+e+f},
$$
$$
N_{a+c+d+e+f,-a}N_{b+d+e,-b-c-2d-2e-f} = | \text{Formulas (075) and (144)} | = -\alpha_7(-\beta_{67}\gamma_{25}\delta_{12}\epsilon_2) = 
$$
$$
=\alpha_7\beta_{67}\gamma_{25}\delta_{12}\epsilon_2,
$$
we have
$$
S = -\alpha_0 N_{b+d+e,a+c+d+e+f} + \alpha_7\beta_{67}\gamma_{25}\delta_{12}\epsilon_2 = 0
$$
and therefore
$$
N_{b+d+e,a+c+d+e+f} = \alpha_{70}\beta_{67}\gamma_{25}\delta_{12}\epsilon_2.\eqno{(156)}
$$

157. We have $(b+d+e)+(a+c+d+e+f)+(-a-b-c-2d-2e-f) = 0,$ therefore (see formula (156))
$$
N_{b+d+e,a+c+d+e+f} = N_{a+c+d+e+f,-a-b-c-2d-2e-f} = N_{-a-b-c-2d-2e-f,b+d+e} = \alpha_{70}\beta_{67}\gamma_{25}\delta_{12}\epsilon_2.\eqno{(157)}
$$

158. We have $(a+c+d+e)+(b+d+e+f)+(-a)+(-b-c-2d-2e-f) = 0$ and $(b+d+e+f)+(-a)\notin E_6,$ hence 
$$
S = N_{a+c+d+e,b+d+e+f}N_{-a,-b-c-2d-2e-f} + N_{-a,a+c+d+e}N_{b+d+e+f,-b-c-2d-2e-f} = 0.
$$
Since
$$
N_{a+c+d+e,b+d+e+f}N_{-a,-b-c-2d-2e-f} = | \text{Table 2} | = -\alpha_0 N_{a+c+d+e,b+d+e+f},
$$
$$
N_{-a,a+c+d+e}N_{b+d+e+f,-b-c-2d-2e-f} = | \text{Formulas (019) and (146)} | = \alpha_4(\beta_{37}\gamma_{25}\epsilon_2) = 
$$
$$
=\alpha_4\beta_{37}\gamma_{25}\epsilon_2,
$$
we have
$$
S = -\alpha_0 N_{a+c+d+e,b+d+e+f} + \alpha_4\beta_{37}\gamma_{25}\epsilon_2 = 0
$$
and therefore
$$
N_{a+c+d+e,b+d+e+f} = \alpha_{40}\beta_{37}\gamma_{25}\epsilon_2.\eqno{(158)}
$$

159. We have $(a+c+d+e)+(b+d+e+f)+(-a-b-c-2d-2e-f) = 0,$ therefore (see formula (158))
$$
N_{a+c+d+e,b+d+e+f} = N_{b+d+e+f,-a-b-c-2d-2e-f} = N_{-a-b-c-2d-2e-f,a+c+d+e} = \alpha_{40}\beta_{37}\gamma_{25}\epsilon_2.\eqno{(159)}
$$

160. We have $c+(a+b+c+2d+2e+f)+(-a-b-2c-2d-2e-f) = 0,$ therefore (see Table 2)
$$
N_{c,a+b+c+2d+2e+f} = N_{a+b+c+2d+2e+f,-a-b-2c-2d-2e-f} = N_{-a-b-2c-2d-2e-f,c} = \beta_9.\eqno{(160)}
$$

161. We have $(a+c)+(b+c+2d+2e+f)+(-c)+(-a-b-c-2d-2e-f) = 0$ and $(b+c+2d+2e+f)+(-c)\notin E_6,$ hence 
$$
S = N_{a+c,b+c+2d+2e+f}N_{-c,-a-b-c-2d-2e-f} + N_{-c,a+c}N_{b+c+2d+2e+f,-a-b-c-2d-2e-f} = 0.
$$
Since
$$
N_{a+c,b+c+2d+2e+f}N_{-c,-a-b-c-2d-2e-f} = | \text{Table 2} | = -\beta_9 N_{a+c,b+c+2d+2e+f},
$$
$$
N_{-c,a+c}N_{b+c+2d+2e+f,-a-b-c-2d-2e-f} = | \text{Formulas (001) and (147)} | = -\alpha_1(\alpha_0) = -\alpha_{10},
$$
we have
$$
S = -\beta_9 N_{a+c,b+c+2d+2e+f} - \alpha_{10} = 0
$$
and therefore
$$
N_{a+c,b+c+2d+2e+f} = -\alpha_{10}\beta_9.\eqno{(161)}
$$

162. We have $(a+c)+(b+c+2d+2e+f)+(-a-b-2c-2d-2e-f) = 0,$ therefore (see formula (161)) 
$$
N_{a+c,b+c+2d+2e+f} = N_{b+c+2d+2e+f,-a-b-2c-2d-2e-f} = N_{-a-b-2c-2d-2e-f,a+c} = -\alpha_{10}\beta_9.\eqno{(162)}
$$

163. We have $e+(a+b+2c+2d+e+f)+(-c)+(-a-b-c-2d-2e-f) = 0$ and $e+(-c)\notin E_6,$ hence 
$$
S = N_{e,a+b+2c+2d+e+f}N_{-c,-a-b-c-2d-2e-f} + N_{a+b+2c+2d+e+f,-c}N_{e,-a-b-c-2d-2e-f} = 0.
$$
Since
$$
N_{e,a+b+2c+2d+e+f}N_{-c,-a-b-c-2d-2e-f} = | \text{Table 2} | = -\beta_9 N_{e,a+b+2c+2d+e+f},
$$
$$
N_{a+b+2c+2d+e+f,-c}N_{e,-a-b-c-2d-2e-f} = | \text{Formulas (123) and (149)} | = -\beta_8(-\alpha_{90}\epsilon_2) = 
$$
$$
=\alpha_{90}\beta_8\epsilon_2,
$$
we have
$$
S = -\beta_9 N_{e,a+b+2c+2d+e+f} + \alpha_{90}\beta_8\epsilon_2 = 0
$$
and therefore
$$
N_{e,a+b+2c+2d+e+f} = \alpha_{90}\beta_{89}\epsilon_2.\eqno{(163)}
$$

164. We have $e+(a+b+2c+2d+e+f)+(-a-b-2c-2d-2e-f) = 0,$ therefore (see formula (163))
$$
N_{e,a+b+2c+2d+e+f} = N_{a+b+2c+2d+e+f,-a-b-2c-2d-2e-f} = N_{-a-b-2c-2d-2e-f,e} = \alpha_{90}\beta_{89}\epsilon_2.\eqno{(164)}
$$

165. We have $(a+b+2c+2d+e)+(e+f)+(-c)+(-a-b-c-2d-2e-f) = 0$ and $(e+f)+(-c)\notin E_6,$ hence 
$$
S = N_{a+b+2c+2d+e,e+f}N_{-c,-a-b-c-2d-2e-f} + N_{-c,a+b+2c+2d+e}N_{e+f,-a-b-c-2d-2e-f} = 0.
$$
Since
$$
N_{a+b+2c+2d+e,e+f}N_{-c,-a-b-c-2d-2e-f} = | \text{Table 2} | = -\beta_9 N_{a+b+2c+2d+e,e+f},
$$
$$
N_{-c,a+b+2c+2d+e}N_{e+f,-a-b-c-2d-2e-f} = | \text{Formulas (055) and (153)} | = 
$$
$$
= \beta_5(\alpha_{60}\beta_{47}\gamma_{2345}\delta_{12}\epsilon_2) = \alpha_{60}\beta_{457}\gamma_{2345}\delta_{12}\epsilon_2,
$$
we have
$$
S = -\beta_9 N_{a+b+2c+2d+e,e+f} + \alpha_{60}\beta_{457}\gamma_{2345}\delta_{12}\epsilon_2 = 0
$$
and therefore
$$
N_{a+b+2c+2d+e,e+f} = \alpha_{60}\beta_{4579}\gamma_{2345}\delta_{12}\epsilon_2.\eqno{(165)}
$$

166. We have $(a+b+2c+2d+e)+(e+f)+(-a-b-2c-2d-2e-f) = 0,$ therefore (see formula (165)) 
$$
N_{a+b+2c+2d+e,e+f} = N_{e+f,-a-b-2c-2d-2e-f} = N_{-a-b-2c-2d-2e-f,a+b+2c+2d+e} = 
$$
$$
=\alpha_{60}\beta_{4579}\gamma_{2345}\delta_{12}\epsilon_2.\eqno{(166)}
$$

167. We have $(c+d+e)+(a+b+c+d+e+f)+(-c)+(-a-b-c-2d-2e-f) = 0$ and $(a+b+c+d+e+f)+(-c)\notin E_6,$ hence 
$$
S = N_{c+d+e,a+b+c+d+e+f}N_{-c,-a-b-c-2d-2e-f} + N_{-c,c+d+e}N_{a+b+c+d+e+f,-a-b-c-2d-2e-f} = 0.
$$
Since
$$
N_{c+d+e,a+b+c+d+e+f}N_{-c,-a-b-c-2d-2e-f} = | \text{Table 2} | = -\beta_9 N_{c+d+e,a+b+c+d+e+f},
$$
$$
N_{-c,c+d+e}N_{a+b+c+d+e+f,-a-b-c-2d-2e-f} = | \text{Formulas (016) and (151)} | = \beta_3(-\alpha_{80}\gamma_{25}\epsilon_2) = 
$$
$$
=-\alpha_{80}\beta_3\gamma_{25}\epsilon_2,
$$
we have
$$
S = -\beta_9 N_{c+d+e,a+b+c+d+e+f} - \alpha_{80}\beta_3\gamma_{25}\epsilon_2 = 0
$$
and therefore
$$
N_{c+d+e,a+b+c+d+e+f} = -\alpha_{80}\beta_{39}\gamma_{25}\epsilon_2.\eqno{(167)}
$$

168. We have $(c+d+e)+(a+b+c+d+e+f)+(-a-b-2c-2d-2e-f) = 0,$ therefore (see formula (167))
$$
N_{c+d+e,a+b+c+d+e+f} = N_{a+b+c+d+e+f,-a-b-2c-2d-2e-f} = N_{-a-b-2c-2d-2e-f,c+d+e} = 
$$
$$ = -\alpha_{80}\beta_{39}\gamma_{25}\epsilon_2.\eqno{(168)}
$$

169. We have $(a+b+c+d+e)+(c+d+e+f)+(-c)+(-a-b-c-2d-2e-f) = 0$ and $(a+b+c+d+e)+(-c)\notin E_6,$ hence 
$$
S = N_{a+b+c+d+e,c+d+e+f}N_{-c,-a-b-c-2d-2e-f} + N_{c+d+e+f,-c}N_{a+b+c+d+e,-a-b-c-2d-2e-f} = 0.
$$
Since
$$
N_{a+b+c+d+e,c+d+e+f}N_{-c,-a-b-c-2d-2e-f} = | \text{Table 2} | = -\beta_9 N_{a+b+c+d+e,c+d+e+f},
$$
$$
N_{c+d+e+f,-c}N_{a+b+c+d+e,-a-b-c-2d-2e-f} = | \text{Formulas (070) and (155)} | = 
$$
$$
= -\beta_6(\alpha_{50}\beta_{47}\gamma_{25}\delta_{12}\epsilon_2) = -\alpha_{50}\beta_{467}\gamma_{25}\delta_{12}\epsilon_2,
$$
we have
$$
S = -\beta_9 N_{a+b+c+d+e,c+d+e+f} - \alpha_{50}\beta_{467}\gamma_{25}\delta_{12}\epsilon_2 = 0
$$
and therefore
$$
N_{a+b+c+d+e,c+d+e+f} = -\alpha_{50}\beta_{4679}\gamma_{25}\delta_{12}\epsilon_2.\eqno{(169)}
$$

170. We have $(a+b+c+d+e)+(c+d+e+f)+(-a-b-2c-2d-2e-f) = 0,$ therefore (see formula (169))
$$
N_{a+b+c+d+e,c+d+e+f} = N_{c+d+e+f,-a-b-2c-2d-2e-f} = N_{-a-b-2c-2d-2e-f,a+b+c+d+e} = 
$$
$$
= -\alpha_{50}\beta_{4679}\gamma_{25}\delta_{12}\epsilon_2.\eqno{(170)}
$$

171. We have $(a+c+d+e)+(b+c+d+e+f)+(-c)+(-a-b-c-2d-2e-f) = 0$ and $(a+c+d+e)+(-c)\notin E_6,$ hence 
$$
S = N_{a+c+d+e,b+c+d+e+f}N_{-c,-a-b-c-2d-2e-f} + N_{b+c+d+e+f,-c}N_{a+c+d+e,-a-b-c-2d-2e-f} = 0.
$$
Since
$$
N_{a+c+d+e,b+c+d+e+f}N_{-c,-a-b-c-2d-2e-f} = | \text{Table 2} | = -\beta_9 N_{a+c+d+e,b+c+d+e+f},
$$
$$
N_{b+c+d+e+f,-c}N_{a+c+d+e,-a-b-c-2d-2e-f} = | \text{Formulas (087) and (159)} | = 
$$
$$
= -\beta_7(-\alpha_{40}\beta_{37}\gamma_{25}\epsilon_2) = \alpha_{40}\beta_{3}\gamma_{25}\epsilon_2,
$$
we have
$$
S = -\beta_9 N_{a+c+d+e,b+c+d+e+f} + \alpha_{40}\beta_{3}\gamma_{25}\epsilon_2 = 0
$$
and therefore
$$
N_{a+c+d+e,b+c+d+e+f} = \alpha_{40}\beta_{39}\gamma_{25}\epsilon_2.\eqno{(171)}
$$

172. We have $(a+c+d+e)+(b+c+d+e+f)+(-a-b-2c-2d-2e-f) = 0,$ therefore (see formula (171))
$$
N_{a+c+d+e,b+c+d+e+f} = N_{b+c+d+e+f,-a-b-2c-2d-2e-f} = N_{-a-b-2c-2d-2e-f,a+c+d+e} = 
$$
$$
=\alpha_{40}\beta_{39}\gamma_{25}\epsilon_2.\eqno{(172)}
$$

173. We have $(b+c+d+e)+(a+c+d+e+f)+(-c)+(-a-b-c-2d-2e-f) = 0$ and $(a+c+d+e+f)+(-c)\notin E_6,$ hence 
$$
S = N_{b+c+d+e,a+c+d+e+f}N_{-c,-a-b-c-2d-2e-f} + N_{-c,b+c+d+e}N_{a+c+d+e+f,-a-b-c-2d-2e-f} = 0.
$$
Since
$$
N_{b+c+d+e,a+c+d+e+f}N_{-c,-a-b-c-2d-2e-f} = | \text{Table 2} | = -\beta_9 N_{b+c+d+e,a+c+d+e+f},
$$
$$
N_{-c,b+c+d+e}N_{a+c+d+e+f,-a-b-c-2d-2e-f} = | \text{Formulas (027) and (157)} | = 
$$
$$ 
= \beta_4(\alpha_{70}\beta_{67}\gamma_{25}\delta_{12}\epsilon_2) = \alpha_{70}\beta_{467}\gamma_{25}\delta_{12}\epsilon_2,
$$
we have
$$
S = -\beta_9 N_{b+c+d+e,a+c+d+e+f} + \alpha_{70}\beta_{467}\gamma_{25}\delta_{12}\epsilon_2 = 0
$$
and therefore
$$
N_{b+c+d+e,a+c+d+e+f} = \alpha_{70}\beta_{4679}\gamma_{25}\delta_{12}\epsilon_2.\eqno{(173)}
$$

174. We have $(b+c+d+e)+(a+c+d+e+f)+(-a-b-2c-2d-2e-f) = 0,$ therefore (see formula (173))
$$
N_{b+c+d+e,a+c+d+e+f} = N_{a+c+d+e+f,-a-b-2c-2d-2e-f} = N_{-a-b-2c-2d-2e-f,b+c+d+e} =
$$
$$
= \alpha_{70}\beta_{4679}\gamma_{25}\delta_{12}\epsilon_2.\eqno{(174)}
$$

175. We have $d+(a+b+2c+2d+2e+f)+(-a-b-2c-3d-2e-f) = 0,$ therefore (see Table 2)
$$
N_{d,a+b+2c+2d+2e+f} = N_{a+b+2c+2d+2e+f,-a-b-2c-3d-2e-f} = N_{-a-b-2c-3d-2e-f,d} =\gamma_6.\eqno{(175)}
$$

176. We have $(d+e)+(a+b+2c+2d+e+f)+(-d)+(-a-b-2c-2d-2e-f) = 0$ and $(a+b+2c+2d+e+f)+(-d)\notin E_6,$ hence 
$$
S = N_{d+e,a+b+2c+2d+e+f}N_{-d,-a-b-2c-2d-2e-f} + N_{-d,d+e}N_{a+b+2c+2d+e+f,-a-b-2c-2d-2e-f} = 0.
$$
Since
$$
N_{d+e,a+b+2c+2d+e+f}N_{-d,-a-b-2c-2d-2e-f} = | \text{Table 2} | = -\gamma_6 N_{d+e,a+b+2c+2d+e+f},
$$
$$
N_{-d,d+e}N_{a+b+2c+2d+e+f,-a-b-2c-2d-2e-f} = | \text{Formulas (015) and (164)} | = \gamma_2(\alpha_{90}\beta_{89}\epsilon_2) = 
$$
$$
=\alpha_{90}\beta_{89}\gamma_2\epsilon_2,
$$
we have
$$
S = -\gamma_6 N_{d+e,a+b+2c+2d+e+f} + \alpha_{90}\beta_{89}\gamma_2\epsilon_2 = 0
$$
and therefore
$$
N_{d+e,a+b+2c+2d+e+f} = \alpha_{90}\beta_{89}\gamma_{26}\epsilon_2.\eqno{(176)}
$$

177. We have $(d+e)+(a+b+2c+2d+e+f)+(-a-b-2c-3d-2e-f) = 0,$ therefore (see formula (176))
$$
N_{d+e,a+b+2c+2d+e+f} = N_{a+b+2c+2d+e+f,-a-b-2c-3d-2e-f} = N_{-a-b-2c-3d-2e-f,d+e}=
$$
$$
=\alpha_{90}\beta_{89}\gamma_{26}\epsilon_2.\eqno{(177)}
$$

178. We have $(c+d)+(a+b+c+2d+2e+f)+(-d)+(-a-b-2c-2d-2e-f) = 0$ and $(a+b+c+2d+2e+f)+(-d)\notin E_6,$ hence 
$$
S = N_{c+d,a+b+c+2d+2e+f}N_{-d,-a-b-2c-2d-2e-f} + N_{-d,c+d}N_{a+b+c+2d+2e+f,-a-b-2c-2d-2e-f} = 0.
$$
Since
$$
N_{c+d,a+b+c+2d+2e+f}N_{-d,-a-b-2c-2d-2e-f} = | \text{Table 2} | = -\gamma_6 N_{c+d,a+b+c+2d+2e+f},
$$
$$
N_{-d,c+d}N_{a+b+c+2d+2e+f,-a-b-2c-2d-2e-f} = | \text{Formulas (002) and (160)} | = -\beta_1(\beta_9) = -\beta_{19},
$$
we have
$$
S = -\gamma_6 N_{c+d,a+b+c+2d+2e+f} - \beta_{19} = 0
$$
and therefore
$$
N_{c+d,a+b+c+2d+2e+f} = -\beta_{19}\gamma_6.\eqno{(178)}
$$

179. We have $(c+d)+(a+b+c+2d+2e+f)+(-a-b-2c-3d-2e-f) = 0,$ therefore (see formula (178)) 
$$
N_{c+d,a+b+c+2d+2e+f} = N_{a+b+c+2d+2e+f,-a-b-2c-3d-2e-f} = N_{-a-b-2c-3d-2e-f,c+d} = -\beta_{19}\gamma_6.\eqno{(179)}
$$

180. We have $(a+b+2c+2d+e)+(d+e+f)+(-d)+(-a-b-2c-2d-2e-f) = 0$ and $(a+b+2c+2d+e)+(-d)\notin E_6,$ hence 
$$
S = N_{a+b+2c+2d+e,d+e+f}N_{-d,-a-b-2c-2d-2e-f} + N_{d+e+f,-d}N_{a+b+2c+2d+e,-a-b-2c-2d-2e-f} = 0.
$$
Since
$$
N_{a+b+2c+2d+e,d+e+f}N_{-d,-a-b-2c-2d-2e-f} = | \text{Table 2} | = -\gamma_6 N_{a+b+2c+2d+e,d+e+f},
$$
$$
N_{d+e+f,-d}N_{a+b+2c+2d+e,-a-b-2c-2d-2e-f} = | \text{Formulas (067) and (166)} | = 
$$
$$
= -\gamma_4(-\alpha_{60}\beta_{4579}\gamma_{2345}\delta_{12}\epsilon_2) = \alpha_{60}\beta_{4579}\gamma_{235}\delta_{12}\epsilon_2,
$$
we have
$$
S = -\gamma_6 N_{a+b+2c+2d+e,d+e+f} + \alpha_{60}\beta_{4579}\gamma_{235}\delta_{12}\epsilon_2 = 0
$$
and therefore
$$
N_{a+b+2c+2d+e,d+e+f} = \alpha_{60}\beta_{4579}\gamma_{2356}\delta_{12}\epsilon_2.\eqno{(180)}
$$

181. We have $(a+b+2c+2d+e)+(d+e+f)+(-a-b-2c-3d-2e-f) = 0,$ therefore (see formula (180))
$$
N_{a+b+2c+2d+e,d+e+f} = N_{d+e+f,-a-b-2c-3d-2e-f} = N_{-a-b-2c-3d-2e-f,a+b+2c+2d+e} = 
$$
$$
=\alpha_{60}\beta_{4579}\gamma_{2356}\delta_{12}\epsilon_2.\eqno{(181)}
$$

182. We have $(c+d+e)+(a+b+c+2d+e+f)+(-d)+(-a-b-2c-2d-2e-f) = 0$ and $(c+d+e)+(-d)\notin E_6,$ hence 
$$
S = N_{c+d+e,a+b+c+2d+e+f}N_{-d,-a-b-2c-2d-2e-f} + \ \ \ \ \ \ \ \ \ \ \ \ \ \ \ \ 
$$
$$
\ \ \ \ \ \ \ \ \ \ \ \ \ + N_{a+b+c+2d+e+f,-d}N_{c+d+e,-a-b-2c-2d-2e-f} = 0. 
$$
Since
$$
N_{c+d+e,a+b+c+2d+e+f}N_{-d,-a-b-2c-2d-2e-f} = | \text{Table 2} | = -\gamma_6 N_{c+d+e,a+b+c+2d+e+f},
$$
$$
N_{a+b+c+2d+e+f,-d}N_{c+d+e,-a-b-2c-2d-2e-f} = | \text{Formulas (114) and (168)} | = 
$$
$$
= -\alpha_{89}\gamma_5(\alpha_{80}\beta_{39}\gamma_{25}\epsilon_2) = -\alpha_{90}\beta_{39}\gamma_2\epsilon_2,
$$
we have
$$
S = -\gamma_6 N_{c+d+e,a+b+c+2d+e+f} - \alpha_{90}\beta_{39}\gamma_2\epsilon_2 = 0
$$
and therefore
$$
N_{c+d+e,a+b+c+2d+e+f} = -\alpha_{90}\beta_{39}\gamma_{26}\epsilon_2.\eqno{(182)}
$$

183. We have $(c+d+e)+(a+b+c+2d+e+f)+(-a-b-2c-3d-2e-f) = 0,$ therefore (see formula (182)) 
$$
N_{c+d+e,a+b+c+2d+e+f} = N_{a+b+c+2d+e+f,-a-b-2c-3d-2e-f} = N_{-a-b-2c-3d-2e-f,c+d+e} = 
$$
$$
= -\alpha_{90}\beta_{39}\gamma_{26}\epsilon_2.\eqno{(183)}
$$

184. We have $(a+c+d)+(b+c+2d+2e+f)+(-d)+(-a-b-2c-2d-2e-f) = 0$ and $(b+c+2d+2e+f)+(-d)\notin E_6,$ hence 
$$
S = N_{a+c+d,b+c+2d+2e+f}N_{-d,-a-b-2c-2d-2e-f} + N_{-d,a+c+d}N_{b+c+2d+2e+f,-a-b-2c-2d-2e-f} = 0.
$$
Since
$$
N_{a+c+d,b+c+2d+2e+f}N_{-d,-a-b-2c-2d-2e-f} = | \text{Table 2} | = -\gamma_6 N_{a+c+d,b+c+2d+2e+f},
$$
$$
N_{-d,a+c+d}N_{b+c+2d+2e+f,-a-b-2c-2d-2e-f} = | \text{Formulas (005) and (162)} | = -\alpha_{12}\beta_1(-\alpha_{10}\beta_9) = 
$$
$$
=\alpha_{20}\beta_{19},
$$
we have
$$
S = -\gamma_6 N_{a+c+d,b+c+2d+2e+f} + \alpha_{20}\beta_{19} = 0
$$
and therefore
$$
N_{a+c+d,b+c+2d+2e+f} = \alpha_{20}\beta_{19}\gamma_6.\eqno{(184)}
$$

185. We have $(a+c+d)+(b+c+2d+2e+f)+(-a-b-2c-3d-2e-f) = 0,$ therefore (see formula (184))
$$
N_{a+c+d,b+c+2d+2e+f} = N_{b+c+2d+2e+f,-a-b-2c-3d-2e-f} = N_{-a-b-2c-3d-2e-f,a+c+d} = \alpha_{20}\beta_{19}\gamma_6.\eqno{(185)}
$$

186. We have $(a+b+c+2d+e)+(c+d+e+f)+(-d)+(-a-b-2c-2d-2e-f) = 0$ and $(c+d+e+f)+(-d)\notin E_6,$ hence 
$$
S = N_{a+b+c+2d+e,c+d+e+f}N_{-d,-a-b-2c-2d-2e-f} + \ \ \ \ \ \ \ \ \ \ \ \ \ \ 
$$
$$
 \ \ \ \ \ \ \ \ \ \ \ \ \ \ \ \ \ \ \ \ + N_{-d,a+b+c+2d+e}N_{c+d+e+f,-a-b-2c-2d-2e-f} = 0.
$$
Since
$$
N_{a+b+c+2d+e,c+d+e+f}N_{-d,-a-b-2c-2d-2e-f} = | \text{Table 2} | = -\gamma_6 N_{a+b+c+2d+e,c+d+e+f},
$$
$$
N_{-d,a+b+c+2d+e}N_{c+d+e+f,-a-b-2c-2d-2e-f} = | \text{Formulas (048) and (170)} | = 
$$
$$
= \alpha_{56}\gamma_3(-\alpha_{50}\beta_{4679}\gamma_{25}\delta_{12}\epsilon_2) = -\alpha_{60}\beta_{4679}\gamma_{235}\delta_{12}\epsilon_2,
$$
we have
$$
S = -\gamma_6 N_{a+b+c+2d+e,c+d+e+f} - \alpha_{60}\beta_{4679}\gamma_{235}\delta_{12}\epsilon_2 = 0
$$
and therefore
$$
N_{a+b+c+2d+e,c+d+e+f} = -\alpha_{60}\beta_{4679}\gamma_{2356}\delta_{12}\epsilon_2.\eqno{(186)}
$$

187. We have $(a+b+c+2d+e)+(c+d+e+f)+(-a-b-2c-3d-2e-f) = 0,$ therefore (see formula (186))
$$
N_{a+b+c+2d+e,c+d+e+f} = N_{c+d+e+f,-a-b-2c-3d-2e-f} = N_{-a-b-2c-3d-2e-f,a+b+c+2d+e} = 
$$
$$
=-\alpha_{60}\beta_{4679}\gamma_{2356}\delta_{12}\epsilon_2.\eqno{(187)}
$$

188. We have $(a+c+d+e)+(b+c+2d+e+f)+(-d)+(-a-b-2c-2d-2e-f) = 0$ and $(a+c+d+e)+(-d)\notin E_6,$ hence 
$$
S = N_{a+c+d+e,b+c+2d+e+f}N_{-d,-a-b-2c-2d-2e-f} + \ \ \ \ \ \ \ \ \ \ \ \ \ \ \ \ \ \ \ \ \ \ 
$$
$$
\ \ \ \ \ \ \ \ \ \ \ \ \ \ \ \ \ \ \ \ \ \  +N_{b+c+2d+e+f,-d}N_{a+c+d+e,-a-b-2c-2d-2e-f} = 0.
$$
Since
$$
N_{a+c+d+e,b+c+2d+e+f}N_{-d,-a-b-2c-2d-2e-f} = | \text{Table 2} | = -\gamma_6 N_{a+c+d+e,b+c+2d+e+f},
$$
$$
N_{b+c+2d+e+f,-d}N_{a+c+d+e,-a-b-2c-2d-2e-f} = | \text{Formulas (103) and (172)} | = 
$$
$$
= -\gamma_5(-\alpha_{40}\beta_{39}\gamma_{25}\epsilon_2) = \alpha_{40}\beta_{39}\gamma_{2}\epsilon_2,
$$
we have
$$
S = -\gamma_6 N_{a+c+d+e,b+c+2d+e+f} + \alpha_{40}\beta_{39}\gamma_{2}\epsilon_2 = 0
$$
and therefore
$$
N_{a+c+d+e,b+c+2d+e+f} = \alpha_{40}\beta_{39}\gamma_{26}\epsilon_2.\eqno{(188)}
$$

189. We have $(a+c+d+e)+(b+c+2d+e+f)+(-a-b-2c-3d-2e-f) = 0,$ therefore (see formula (188))
$$
N_{a+c+d+e,b+c+2d+e+f} = N_{b+c+2d+e+f,-a-b-2c-3d-2e-f} = N_{-a-b-2c-3d-2e-f,a+c+d+e} = 
$$
$$
= \alpha_{40}\beta_{39}\gamma_{26}\epsilon_2.\eqno{(189)}
$$

190. We have $(b+c+2d+e)+(a+c+d+e+f)+(-d)+(-a-b-2c-2d-2e-f) = 0$ and $(a+c+d+e+f)+(-d)\notin E_6,$ hence 
$$
S = N_{b+c+2d+e,a+c+d+e+f}N_{-d,-a-b-2c-2d-2e-f} + \ \ \ \ \ \ \ \ \ \ \ \ \ \ \ \ \ \ \ \ \ \ 
$$
$$
\ \ \ \ \ \ \ \ \ \ \ \ \ \ \ \ \ \ \ \ \ \ +N_{-d,b+c+2d+e}N_{a+c+d+e+f,-a-b-2c-2d-2e-f} = 0.
$$
Since
$$
N_{b+c+2d+e,a+c+d+e+f}N_{-d,-a-b-2c-2d-2e-f} = | \text{Table 2} | = -\gamma_6 N_{b+c+2d+e,a+c+d+e+f},
$$
$$
N_{-d,b+c+2d+e}N_{a+c+d+e+f,-a-b-2c-2d-2e-f} = | \text{Formulas (039) and (174)} | = 
$$
$$
= \gamma_3(\alpha_{70}\beta_{4679}\gamma_{25}\delta_{12}\epsilon_2) = \alpha_{70}\beta_{4679}\gamma_{235}\delta_{12}\epsilon_2,
$$
we have
$$
S = -\gamma_6 N_{b+c+2d+e,a+c+d+e+f} + \alpha_{70}\beta_{4679}\gamma_{235}\delta_{12}\epsilon_2 = 0
$$
and therefore
$$
N_{b+c+2d+e,a+c+d+e+f} = \alpha_{70}\beta_{4679}\gamma_{2356}\delta_{12}\epsilon_2.\eqno{(190)}
$$

191. We have $(b+c+2d+e)+(a+c+d+e+f)+(-a-b-2c-3d-2e-f) = 0,$ therefore (see formula (190))
$$
N_{b+c+2d+e,a+c+d+e+f} = N_{a+c+d+e+f,-a-b-2c-3d-2e-f} = N_{-a-b-2c-3d-2e-f,b+c+2d+e} = 
$$
$$
=\alpha_{70}\beta_{4679}\gamma_{2356}\delta_{12}\epsilon_2.\eqno{(191)}
$$

192. We have $b+(a+b+2c+3d+2e+f)+(-a-2b-2c-3d-2e-f) = 0,$ therefore (see Table 2)
$$
N_{b,a+b+2c+3d+2e+f} = N_{a+b+2c+3d+2e+f,-a-2b-2c-3d-2e-f} = N_{-a-2b-2c-3d-2e-f,b} = \delta_3.\eqno{(192)}
$$

193. We have $(b+d)+(a+b+2c+2d+2e+f)+(-b)+(-a-b-2c-3d-2e-f) = 0$ and $(a+b+2c+2d+2e+f)+(-b)\notin E_6,$ hence 
$$
S = N_{b+d,a+b+2c+2d+2e+f}N_{-b,-a-b-2c-3d-2e-f} + N_{-b,b+d}N_{a+b+2c+2d+2e+f,-a-b-2c-3d-2e-f} = 0.
$$
Since
$$
N_{b+d,a+b+2c+2d+2e+f}N_{-b,-a-b-2c-3d-2e-f} = | \text{Table 2} | = -\delta_3 N_{b+d,a+b+2c+2d+2e+f},
$$
$$
N_{-b,b+d}N_{a+b+2c+2d+2e+f,-a-b-2c-3d-2e-f} = | \text{Formulas (006) and (175)} | = -\gamma_1(\gamma_6) = -\gamma_{16},
$$
we have
$$
S = -\delta_3 N_{b+d,a+b+2c+2d+2e+f} - \gamma_{16} = 0
$$
and therefore
$$
N_{b+d,a+b+2c+2d+2e+f} = -\gamma_{16}\delta_3.\eqno{(193)}
$$

194. We have $(b+d)+(a+b+2c+2d+2e+f)+(-a-2b-2c-3d-2e-f) = 0,$ therefore (see formula (193))
$$
N_{b+d,a+b+2c+2d+2e+f} = N_{a+b+2c+2d+2e+f,-a-2b-2c-3d-2e-f} = N_{-a-2b-2c-3d-2e-f,b+d} = -\gamma_{16}\delta_3.\eqno{(194)}
$$

195. We have $(b+d+e)+(a+b+2c+2d+e+f)+(-b)+(-a-b-2c-3d-2e-f) = 0$ and $(a+b+2c+2d+e+f)+(-b)\notin E_6,$ hence 
$$
S = N_{b+d+e,a+b+2c+2d+e+f}N_{-b,-a-b-2c-3d-2e-f} + N_{-b,b+d+e}N_{a+b+2c+2d+e+f,-a-b-2c-3d-2e-f} = 0.
$$
Since
$$
N_{b+d+e,a+b+2c+2d+e+f}N_{-b,-a-b-2c-3d-2e-f} = | \text{Table 2} | = -\delta_3 N_{b+d+e,a+b+2c+2d+e+f},
$$
$$
N_{-b,b+d+e}N_{a+b+2c+2d+e+f,-a-b-2c-3d-2e-f} = | \text{Formulas (024) and (177)} | = 
$$
$$
= \delta_1(\alpha_{90}\beta_{89}\gamma_{26}\epsilon_2) = \alpha_{90}\beta_{89}\gamma_{26}\delta_1\epsilon_2,
$$
we have
$$
S = -\delta_3 N_{b+d+e,a+b+2c+2d+e+f} + \alpha_{90}\beta_{89}\gamma_{26}\delta_1\epsilon_2 = 0
$$
and therefore
$$
N_{b+d+e,a+b+2c+2d+e+f} =  \alpha_{90}\beta_{89}\gamma_{26}\delta_{13}\epsilon_2.\eqno{(195)}
$$

196. We have $(b+d+e)+(a+b+2c+2d+e+f)+(-a-2b-2c-3d-2e-f) = 0,$ therefore (see formula (195))
$$
N_{b+d+e,a+b+2c+2d+e+f} = N_{a+b+2c+2d+e+f,-a-2b-2c-3d-2e-f} = N_{-a-2b-2c-3d-2e-f,b+d+e} = 
$$
$$
=\alpha_{90}\beta_{89}\gamma_{26}\delta_{13}\epsilon_2.\eqno{(196)}
$$

197. We have $(b+c+d)+(a+b+c+2d+2e+f)+(-b)+(-a-b-2c-3d-2e-f) = 0$ and $(a+b+c+2d+2e+f)+(-b)\notin E_6,$ hence 
$$
S = N_{b+c+d,a+b+c+2d+2e+f}N_{-b,-a-b-2c-3d-2e-f} + N_{-b,b+c+d}N_{a+b+c+2d+2e+f,-a-b-2c-3d-2e-f} = 0.
$$
Since
$$
N_{b+c+d,a+b+c+2d+2e+f}N_{-b,-a-b-2c-3d-2e-f} = | \text{Table 2} | = -\delta_3 N_{b+c+d,a+b+c+2d+2e+f},
$$
$$
N_{-b,b+c+d}N_{a+b+c+2d+2e+f,-a-b-2c-3d-2e-f} = | \text{Formulas (009) and (179)} | = -\beta_{12}\gamma_1(-\beta_{19}\gamma_6) = 
$$
$$
=\beta_{29}\gamma_{16},
$$
we have
$$
S = -\delta_3 N_{b+c+d,a+b+c+2d+2e+f} + \beta_{29}\gamma_{16} = 0
$$
and therefore
$$
N_{b+c+d,a+b+c+2d+2e+f} = \beta_{29}\gamma_{16}\delta_3.\eqno{(197)}
$$

198. We have $(b+c+d)+(a+b+c+2d+2e+f)+(-a-2b-2c-3d-2e-f) = 0,$ therefore (see formula (197))
$$
N_{b+c+d,a+b+c+2d+2e+f} = N_{a+b+c+2d+2e+f,-a-2b-2c-3d-2e-f} = N_{-a-2b-2c-3d-2e-f,b+c+d} =
$$
$$
= \beta_{29}\gamma_{16}\delta_3.\eqno{(198)}
$$

199. We have $(a+b+2c+2d+e)+(b+d+e+f)+(-b)+(-a-b-2c-3d-2e-f) = 0$ and $(a+b+2c+2d+e)+(-b)\notin E_6,$ hence 
$$
S = N_{a+b+2c+2d+e,b+d+e+f}N_{-b,-a-b-2c-3d-2e-f} + N_{b+d+e+f,-b}N_{a+b+2c+2d+e,-a-b-2c-3d-2e-f} = 0.
$$
Since
$$
N_{a+b+2c+2d+e,b+d+e+f}N_{-b,-a-b-2c-3d-2e-f} = | \text{Table 2} | = -\delta_3 N_{a+b+2c+2d+e,b+d+e+f},
$$
$$
N_{b+d+e+f,-b}N_{a+b+2c+2d+e,-a-b-2c-3d-2e-f} = | \text{Formulas (082) and (181)} | = 
$$
$$
= -\delta_2(-\alpha_{60}\beta_{4579}\gamma_{2356}\delta_{12}\epsilon_2) = \alpha_{60}\beta_{4579}\gamma_{2356}\delta_{1}\epsilon_2,
$$
we have
$$
S = -\delta_3 N_{a+b+2c+2d+e,b+d+e+f} + \alpha_{60}\beta_{4579}\gamma_{2356}\delta_{1}\epsilon_2 = 0
$$
and therefore
$$
N_{a+b+2c+2d+e,b+d+e+f} = \alpha_{60}\beta_{4579}\gamma_{2356}\delta_{13}\epsilon_2.\eqno{(199)}
$$

200. We have $(a+b+2c+2d+e)+(b+d+e+f)+(-a-2b-2c-3d-2e-f) = 0,$ therefore (see formula (199))
$$
N_{a+b+2c+2d+e,b+d+e+f} = N_{b+d+e+f,-a-2b-2c-3d-2e-f} = N_{-a-2b-2c-3d-2e-f,a+b+2c+2d+e} =
$$
$$
= \alpha_{60}\beta_{4579}\gamma_{2356}\delta_{13}\epsilon_2.\eqno{(200)}
$$

201. We have $(b+c+d+e)+(a+b+c+2d+e+f)+(-b)+(-a-b-2c-3d-2e-f) = 0$ and $(a+b+c+2d+e+f)+(-b)\notin E_6,$ hence 
$$
S = N_{b+c+d+e,a+b+c+2d+e+f}N_{-b,-a-b-2c-3d-2e-f} + 
$$
$$
+N_{-b,b+c+d+e}N_{a+b+c+2d+e+f,-a-b-2c-3d-2e-f} = 0.
$$
Since
$$
N_{b+c+d+e,a+b+c+2d+e+f}N_{-b,-a-b-2c-3d-2e-f} = | \text{Table 2} | = -\delta_3 N_{b+c+d+e,a+b+c+2d+e+f},
$$
$$
N_{-b,b+c+d+e}N_{a+b+c+2d+e+f,-a-b-2c-3d-2e-f} = | \text{Formulas (031) and (183)} | = 
$$
$$
= \beta_{34}\delta_1(-\alpha_{90}\beta_{39}\gamma_{26}\epsilon_2) = -\alpha_{90}\beta_{49}\gamma_{26}\delta_1\epsilon_2,
$$
we have
$$
S = -\delta_3 N_{b+c+d+e,a+b+c+2d+e+f} - \alpha_{90}\beta_{49}\gamma_{26}\delta_1\epsilon_2 = 0
$$
and therefore
$$
N_{b+c+d+e,a+b+c+2d+e+f} = - \alpha_{90}\beta_{49}\gamma_{26}\delta_{13}\epsilon_2.\eqno{(201)}
$$

202. We have $(b+c+d+e)+(a+b+c+2d+e+f)+(-a-2b-2c-3d-2e-f) = 0,$ therefore (see formula (201))
$$
N_{b+c+d+e,a+b+c+2d+e+f} = N_{a+b+c+2d+e+f,-a-2b-2c-3d-2e-f} = N_{-a-2b-2c-3d-2e-f,b+c+d+e} =
$$
$$
= -\alpha_{90}\beta_{49}\gamma_{26}\delta_{13}\epsilon_2.\eqno{(202)}
$$

203. We have $(a+b+c+d)+(b+c+2d+2e+f)+(-b)+(-a-b-2c-3d-2e-f) = 0$ and $(b+c+2d+2e+f)+(-b)\notin E_6,$ hence 
$$
S = N_{a+b+c+d,b+c+2d+2e+f}N_{-b,-a-b-2c-3d-2e-f} + 
$$
$$
+N_{-b,a+b+c+d}N_{b+c+2d+2e+f,-a-b-2c-3d-2e-f} = 0.
$$
Since
$$
N_{a+b+c+d,b+c+2d+2e+f}N_{-b,-a-b-2c-3d-2e-f} = | \text{Table 2} | = -\delta_3 N_{a+b+c+d,b+c+2d+2e+f},
$$
$$
N_{-b,a+b+c+d}N_{b+c+2d+2e+f,-a-b-2c-3d-2e-f} = | \text{Formulas (014) and (185)} | = 
$$
$$
= -\alpha_{23}\beta_{12}\gamma_1(\alpha_{20}\beta_{19}\gamma_6) = -\alpha_{30}\beta_{29}\gamma_{16},
$$
we have
$$
S = -\delta_3 N_{a+b+c+d,b+c+2d+2e+f} - \alpha_{30}\beta_{29}\gamma_{16} = 0
$$
and therefore
$$
N_{a+b+c+d,b+c+2d+2e+f} = -\alpha_{30}\beta_{29}\gamma_{16}\delta_3.\eqno{(203)}
$$

204. We have $(a+b+c+d)+(b+c+2d+2e+f)+(-a-2b-2c-3d-2e-f) = 0,$ therefore (see formula (203))
$$
N_{a+b+c+d,b+c+2d+2e+f} = N_{b+c+2d+2e+f,-a-2b-2c-3d-2e-f} = N_{-a-2b-2c-3d-2e-f,a+b+c+d} = 
$$
$$
=-\alpha_{30}\beta_{29}\gamma_{16}\delta_3.\eqno{(204)}
$$

205. We have $(b+c+2d+e)+(a+b+c+d+e+f)+(-b)+(-a-b-2c-3d-2e-f) = 0$ and $(b+c+2d+e)+(-b)\notin E_6,$ hence 
$$
S = N_{b+c+2d+e,a+b+c+d+e+f}N_{-b,-a-b-2c-3d-2e-f} + 
$$
$$
+N_{a+b+c+d+e+f,-b}N_{b+c+2d+e,-a-b-2c-3d-2e-f} = 0.
$$
Since
$$
N_{b+c+2d+e,a+b+c+d+e+f}N_{-b,-a-b-2c-3d-2e-f} = | \text{Table 2} | = -\delta_3 N_{b+c+2d+e,a+b+c+d+e+f},
$$
$$
N_{a+b+c+d+e+f,-b}N_{b+c+2d+e,-a-b-2c-3d-2e-f} = | \text{Formulas (100) and (191)} | =
$$
$$
= -\alpha_{78}\beta_{67}\delta_2(-\alpha_{70}\beta_{4679}\gamma_{2356}\delta_{12}\epsilon_2) = \alpha_{80}\beta_{49}\gamma_{2356}\delta_{1}\epsilon_2,
$$
we have
$$
S = -\delta_3 N_{b+c+2d+e,a+b+c+d+e+f} + \alpha_{80}\beta_{49}\gamma_{2356}\delta_{1}\epsilon_2 = 0
$$
and therefore
$$
N_{b+c+2d+e,a+b+c+d+e+f} = \alpha_{80}\beta_{49}\gamma_{2356}\delta_{13}\epsilon_2.\eqno{(205)}
$$

206. We have $(b+c+2d+e)+(a+b+c+d+e+f)+(-a-2b-2c-3d-2e-f) = 0,$ therefore (see formula (205))
$$
N_{b+c+2d+e,a+b+c+d+e+f} = N_{a+b+c+d+e+f,-a-2b-2c-3d-2e-f} = N_{-a-2b-2c-3d-2e-f,b+c+2d+e} = 
$$
$$
= \alpha_{80}\beta_{49}\gamma_{2356}\delta_{13}\epsilon_2.\eqno{(206)}
$$

207. We have $(a+b+c+d+e)+(b+c+2d+e+f)+(-b)+(-a-b-2c-3d-2e-f) = 0$ and $(b+c+2d+e+f)+(-b)\notin E_6,$ hence 
$$
S = N_{a+b+c+d+e,b+c+2d+e+f}N_{-b,-a-b-2c-3d-2e-f} + 
$$
$$
+N_{-b,a+b+c+d+e}N_{b+c+2d+e+f,-a-b-2c-3d-2e-f} = 0.
$$
Since
$$
N_{a+b+c+d+e,b+c+2d+e+f}N_{-b,-a-b-2c-3d-2e-f} = | \text{Table 2} | = -\delta_3 N_{a+b+c+d+e,b+c+2d+e+f},
$$
$$
N_{-b,a+b+c+d+e}N_{b+c+2d+e+f,-a-b-2c-3d-2e-f} = | \text{Formulas (038) and (189)} | = 
$$
$$
= \alpha_{45}\beta_{34}\delta_1(\alpha_{40}\beta_{39}\gamma_{26}\epsilon_2) = \alpha_{50}\beta_{49}\gamma_{26}\delta_1\epsilon_2,
$$
we have
$$
S = -\delta_3 N_{a+b+c+d+e,b+c+2d+e+f} + \alpha_{50}\beta_{49}\gamma_{26}\delta_1\epsilon_2 = 0
$$
and therefore
$$
N_{a+b+c+d+e,b+c+2d+e+f} = \alpha_{50}\beta_{49}\gamma_{26}\delta_{13}\epsilon_2.\eqno{(207)}
$$

208. We have $(a+b+c+d+e)+(b+c+2d+e+f)+(-a-2b-2c-3d-2e-f) = 0,$ therefore (see formula (207))
$$
N_{a+b+c+d+e,b+c+2d+e+f} = N_{b+c+2d+e+f,-a-2b-2c-3d-2e-f} = N_{-a-2b-2c-3d-2e-f,a+b+c+d+e} = 
$$
$$
= \alpha_{50}\beta_{49}\gamma_{26}\delta_{13}\epsilon_2.\eqno{(208)}
$$

209. We have $(a+b+c+2d+e)+(b+c+d+e+f)+(-b)+(-a-b-2c-3d-2e-f) = 0$ and $(a+b+c+2d+e)+(-b)\notin E_6,$ hence 
$$
S = N_{a+b+c+2d+e,b+c+d+e+f}N_{-b,-a-b-2c-3d-2e-f} + 
$$
$$
+N_{b+c+d+e+f,-b}N_{a+b+c+2d+e,-a-b-2c-3d-2e-f} = 0.
$$
Since
$$
N_{a+b+c+2d+e,b+c+d+e+f}N_{-b,-a-b-2c-3d-2e-f} = | \text{Table 2} | = -\delta_3 N_{a+b+c+2d+e,b+c+d+e+f},
$$
$$
N_{b+c+d+e+f,-b}N_{a+b+c+2d+e,-a-b-2c-3d-2e-f} = | \text{Formulas (091) and (187)} | = 
$$
$$
= -\beta_{67}\delta_2(\alpha_{60}\beta_{4679}\gamma_{2356}\delta_{12}\epsilon_2) = -\alpha_{60}\beta_{49}\gamma_{2356}\delta_{1}\epsilon_2,
$$
we have
$$
S = -\delta_3 N_{a+b+c+2d+e,b+c+d+e+f} - \alpha_{60}\beta_{49}\gamma_{2356}\delta_{1}\epsilon_2 = 0
$$
and therefore
$$
N_{a+b+c+2d+e,b+c+d+e+f} = -\alpha_{60}\beta_{49}\gamma_{2356}\delta_{13}\epsilon_2.\eqno{(209)}
$$

210. We have $(a+b+c+2d+e)+(b+c+d+e+f)+(-a-2b-2c-3d-2e-f) = 0,$ therefore (see formula (209))
$$
N_{a+b+c+2d+e,b+c+d+e+f} = N_{b+c+d+e+f,-a-2b-2c-3d-2e-f} = N_{-a-2b-2c-3d-2e-f,a+b+c+2d+e} = 
$$ 
$$
=- \alpha_{60}\beta_{49}\gamma_{2356}\delta_{13}\epsilon_2.\eqno{(210)}
$$

To conclude the proof, we present the diagrams that show the dependencies of the formulas.

\newpage

\thicklines
{

}

The theorem is completely proven. \hfill$\square$

\subsection{Chevalley's Commutator Formulas}

Let $\Phi$ be a reduced indecomposable root system, $K$ is a
field. According to \cite[Theorem 5.2.2]{Car72} the commutator
$$
[x_s(u),x_r(t)]=x_s(u)^{-1}x_r(t)^{-1}x_s(u)x_r(t)
$$
where $r,s\in\Phi$ and $u,t\in K,$ is equal to one if $r+s\notin\Phi$ and $r\ne-s,$ and expands
into the product of root elements according to the formula
$$
\left[x_s(u),x_r(t)\right]=
 \prod_{{\scriptsize\begin{array}{c} ir+js\in\Phi,\\ i,j>0\end{array}}}x_{ir+js}(C_{ij,rs}( -t)^iu^j),
$$
if $r+s\in\Phi.$ The product is taken over all pairs of positive integers $i,j$ for which $ir+js$ is a root, in order of increasing $i+j.$ The constant $C_{ij,rs}$
are integers and are determined by the formulas \cite[theorem 5.2.2]{Car72}:
\begin{align*}
    C_{i1,rs} &= M_{r,s,i}, \\[3mm]
    C_{1j,rs} &= (-1)^j M_{s,r,j}, \\[3mm]
    C_{32,rs} &= \frac{1}{3}M_{r+s,r,2}, \\[3mm]
    C_{23,rs} &= -\frac{2}{3}M_{s+r,s,2}.
\end{align*}
In turn, the numbers $M_{r,s,i}$ are expressed in terms of structure
constants $N_{r,s}$ of the corresponding Lie algebra according to the formula \cite[p. 61]{Car72}
$$
M_{r,s,i}=\frac{1}{i!}N_{r,s}N_{r,r+s}\ldots N_{r,(i-1)r+s},
$$

Let $r,s,r+s\in\Phi.$ We have 
$$
C_{11,rs} = N_{rs} = -N_{sr}
$$
and
$$
C_{11,-r,-s} = N_{-r,-s} = N_{sr}.
$$

It follows that when $2r+s,r+2s\notin\Phi,$ we have
$$
[x_s(u),x_r(t)] = x_{r+s}(C_{11,rs}(-t)u) = x_{r+s}(N_{sr}tu)
$$ 
and
$$
[x_{-s}(u),x_{-r}(t)] = x_{-r-s}(C_{11,-r,-s}(-t)u) = x_{r+s}(-N_{sr}tu).
$$

\subsection{List of Formulas}

\begin{center}
\textit{Positive roots}
\end{center}

{\boldmath$a:$}

$$
\begin{array}{rcl}
[x_a(u),x_{c}(t)] &=& x_{a+c}(N_{a,c}tu)\\[3mm]
                  &=& x_{a+c}( \alpha_1 tu).
\end{array}\eqno{(001)}
$$

$$
\begin{array}{rcl}
[x_a(u),x_{c+d}(t)] &=& x_{a+c+d}(N_{a,c+d}tu)\\[3mm]
                    &=& x_{a+c+d}( \alpha_2 tu).
\end{array}\eqno{(002)}
$$

$$
\begin{array}{rcl}
[x_a(u),x_{b+c+d}(t)] &=& x_{a+b+c+d}(N_{a,b+c+d}tu)\\[3mm] 
                      &=& x_{a+b+c+d}( \alpha_3 tu).
\end{array}\eqno{(003)}
$$

$$
\begin{array}{rcl}
[x_a(u),x_{c+d+e}(t)] &=& x_{a+c+d+e}(N_{a,c+d+e}tu)\\[3mm] 
                      &=& x_{a+c+d+e}( \alpha_4 tu).
\end{array}\eqno{(004)}
$$

$$
\begin{array}{rcl}
[x_a(u),x_{b+c+d+e}(t)] &=& x_{a+b+c+d+e}(N_{a,b+c+d+e}tu)\\[3mm]
                        &=& x_{a+b+c+d+e}( \alpha_5 tu).
\end{array}\eqno{(005)}
$$

$$
\begin{array}{rcl}
[x_a(u),x_{b+c+2d+e}(t)] &=& x_{a+b+c+2d+e}(N_{a,b+c+2d+e}tu)\\[3mm]
                         &=& x_{a+b+c+2d+e}( \alpha_6 tu).
\end{array}\eqno{(006)}
$$

$$
\begin{array}{rcl}
[x_a(u),x_{c+d+e+f}(t)] &=& x_{a+c+d+e+f}(N_{a,c+d+e+f}tu)\\[3mm]
                        &=& x_{a+c+d+e+f}( \alpha_7 tu).
\end{array}\eqno{(007)}
$$

$$
\begin{array}{rcl}
[x_a(u),x_{b+c+d+e+f}(t)] &=& x_{a+b+c+d+e+f}(N_{a,b+c+d+e+f}tu)\\[3mm]
                          &=& x_{a+b+c+d+e+f}( \alpha_8 tu).
\end{array}\eqno{(008)}
$$

$$
\begin{array}{rcl}
[x_a(u),x_{b+c+2d+e+f}(t)] &=& x_{a+b+c+2d+e+f}(N_{a,b+c+2d+e+f}tu)\\[3mm]
                           &=& x_{a+b+c+2d+e+f}( \alpha_9 tu).
\end{array}\eqno{(009)}
$$

$$
\begin{array}{rcl}
[x_a(u),x_{b+c+2d+2e+f}(t)] &=& x_{a+b+c+2d+2e+f}(N_{a,b+c+2d+2e+f}tu)\\[3mm]
                            &=& x_{a+b+c+2d+2e+f}( \alpha_0 tu).
\end{array}\eqno{(010)}
$$

{\boldmath$c:$}

$$
\begin{array}{rcl}
[x_c(u),x_{d}(t)] &=& x_{c+d}(N_{c,d} tu)\\[3mm]
                  &=& x_{c+d}( \beta_1 tu).
\end{array}\eqno{(011)}
$$

$$
\begin{array}{rcl}
[x_c(u),x_{b+d}(t)] &=& x_{b+c+d}(N_{c,b+d} tu)\\[3mm]
                    &=& x_{b+c+d}( \beta_2 tu).
\end{array}\eqno{(012)}
$$

$$
\begin{array}{rcl}
[x_c(u),x_{d+e}(t)] &=& x_{c+d+e}(N_{c,d+e} tu)\\[3mm] 
                    &=& x_{c+d+e}( \beta_3 tu).
\end{array}\eqno{(013)}
$$

$$
\begin{array}{rcl}
[x_c(u),x_{b+d+e}(t)] &=& x_{b+c+d+e}(N_{c,b+d+e} tu)\\[3mm]
                      &=& x_{b+c+d+e}( \beta_4 tu).
\end{array}\eqno{(014)}
$$

$$
\begin{array}{rcl}
[x_c(u),x_{a+b+c+2d+e}(t)] &=& x_{a+b+2c+2d+e}(N_{c,a+b+c+2d+e} tu)\\[3mm]
                           &=& x_{a+b+2c+2d+e}( \beta_5 tu).
\end{array}\eqno{(015)}
$$

$$
\begin{array}{rcl}
[x_c(u),x_{d+e+f}(t)] &=& x_{c+d+e+f}(N_{c,d+e+f} tu)\\[3mm]
                      &=& x_{c+d+e+f}( \beta_6 tu).
\end{array}\eqno{(016)}
$$

$$
\begin{array}{rcl}
[x_c(u),x_{b+d+e+f}(t)] &=& x_{b+c+d+e+f}(N_{c,b+d+e+f} tu)\\[3mm]
                        &=& x_{b+c+d+e+f}( \beta_7 tu).
\end{array}\eqno{(017)}
$$

$$
\begin{array}{rcl}
[x_c(u),x_{a+b+c+2d+e+f}(t)] &=& x_{a+b+2c+2d+e+f}(N_{c,a+b+c+2d+e+f} tu)\\[3mm]
                             &=& x_{a+b+2c+2d+e+f}( \beta_8 tu).
\end{array}\eqno{(018)}
$$

$$
\begin{array}{rcl}
[x_c(u),x_{a+b+c+2d+2e+f}(t)] &=& x_{a+b+2c+2d+2e+f}(N_{c,a+b+c+2d+2e+f} tu)\\[3mm]
                              &=& x_{a+b+2c+2d+2e+f}( \beta_9 tu).
\end{array}\eqno{(019)}
$$

{\boldmath$a+c:$}

$$
\begin{array}{rcl}
[x_{a+c}(u),x_{d}(t)] &=& x_{a+c+d}(N_{a+c,d} tu)\\[3mm]
                      &=& x_{a+c+d}( \alpha_{12}\beta_1 tu).
\end{array}\eqno{(020)}
$$

$$
\ \ \ \begin{array}{rcl}
[x_{a+c}(u),x_{b+d}(t)] &=& x_{a+b+c+d}(N_{a+c,b+d} tu)\\[3mm]
                        &=& x_{a+b+c+d}( \alpha_{13}\beta_2 tu).
\end{array}\eqno{(021)}
$$

$$
\begin{array}{rcl}
[x_{a+c}(u),x_{d+e}(t)] &=& x_{a+c+d+e}(N_{a+c,d+e} tu)\\[3mm]
                        &=& x_{a+c+d+e}( \alpha_{14}\beta_3 tu).
\end{array}\eqno{(022)}
$$

$$
\begin{array}{rcl}
[x_{a+c}(u),x_{b+d+e}(t)] &=& x_{a+b+c+d+e}(N_{a+c,b+d+e} tu)\\[3mm]
                          &=& x_{a+b+c+d+e}( \alpha_{15}\beta_4 tu).
\end{array}\eqno{(023)}
$$

$$
\begin{array}{rcl}
[x_{a+c}(u),x_{b+c+2d+e}(t)] &=& x_{a+b+2c+2d+e}(N_{a+c,b+c+2d+e} tu)\\[3mm]
                             &=& x_{a+b+2c+2d+e}( -\alpha_{16}\beta_5 tu).
\end{array}\eqno{(024)}
$$

$$
\begin{array}{rcl}
[x_{a+c}(u),x_{d+e+f}(t)] &=& x_{a+c+d+e+f}(N_{a+c,d+e+f} tu)\\[3mm]
                          &=& x_{a+c+d+e+f}( \alpha_{17}\beta_6 tu).
\end{array}\eqno{(025)}
$$

$$
\begin{array}{rcl}
[x_{a+c}(u),x_{b+d+e+f}(t)] &=& x_{a+b+c+d+e+f}(N_{a+c,b+d+e+f} tu)\\[3mm]
                            &=& x_{a+b+c+d+e+f}( \alpha_{18}\beta_7 tu).
\end{array}\eqno{(026)}
$$

$$
\begin{array}{rcl}
[x_{a+c}(u),x_{b+c+2d+e+f}(t)] &=& x_{a+b+2c+2d+e+f}(N_{a+c,b+c+2d+e+f} tu)\\[3mm] 
                               &=& x_{a+b+2c+2d+e+f}( -\alpha_{19}\beta_8 tu).
\end{array}\eqno{(027)}
$$

$$
\begin{array}{rcl}
[x_{a+c}(u),x_{b+c+2d+2e+f}(t)] &=& x_{a+b+2c+2d+2e+f}(N_{a+c,b+c+2d+2e+f} tu)\\[3mm]
                                &=& x_{a+b+2c+2d+2e+f}( -\alpha_{10}\beta_9 tu).
\end{array}\eqno{(028)}
$$

{\boldmath$d:$}

$$
\begin{array}{rcl}
[x_d(u),x_{b}(t)] &=& x_{b+d}(N_{d,b} tu)\\[3mm]
                  &=& x_{b+d}( \gamma_1 tu).
\end{array}\eqno{(029)}
$$

$$
\begin{array}{rcl}
[x_d(u),x_{e}(t)] &=& x_{d+e}(N_{d,e} tu)\\[3mm] 
                  &=& x_{d+e}( \gamma_2 tu).
\end{array}\eqno{(030)}
$$

$$
\begin{array}{rcl}
[x_d(u),x_{b+c+d+e}(t)] &=& x_{b+c+2d+e}(N_{d,b+c+d+e} tu)\\[3mm] 
                        &=& x_{b+c+2d+e}( \gamma_3 tu).
\end{array}\eqno{(031)}
$$

$$
\begin{array}{rcl}
[x_d(u),x_{a+b+c+d+e}(t)] &=& x_{a+b+c+2d+e}(N_{d,a+b+c+d+e} tu)\\[3mm]
                          &=& x_{a+b+c+2d+e}(  \alpha_{56}\gamma_3 tu).
\end{array}\eqno{(032)}
$$

$$
\begin{array}{rcl}
[x_d(u),x_{e+f}(t)] &=& x_{d+e+f}(N_{d,e+f} tu)\\[3mm]
                    &=& x_{d+e+f}( \gamma_4 tu).
\end{array}\eqno{(033)}
$$

$$
\begin{array}{rcl}
[x_d(u),x_{b+c+d+e+f}(t)] &=& x_{b+c+2d+e+f}(N_{d,b+c+d+e+f} tu)\\[3mm]
                          &=& x_{b+c+2d+e+f}( \gamma_5 tu).
\end{array}\eqno{(034)}
$$

$$
\begin{array}{rcl}
[x_d(u),x_{a+b+c+d+e+f}(t)] &=& x_{a+b+c+2d+e+f}(N_{d,a+b+c+d+e+f} tu)\\[3mm]
                            &=& x_{a+b+c+2d+e+f}( \alpha_{89}\gamma_5 tu).
\end{array}\eqno{(035)}
$$

$$
\begin{array}{rcl}
[x_d(u),x_{a+b+2c+2d+2e+f}(t)] &=& x_{a+b+2c+3d+2e+f}(N_{d,a+b+2c+2d+2e+f} tu)\\[3mm]
                               &=& x_{a+b+2c+3d+2e+f}(  \gamma_6 tu).
\end{array}\eqno{(036)}
$$

{\boldmath$c+d:$}

$$
\begin{array}{rcl}
[x_{c+d}(u),x_{b}(t)] &=& x_{b+c+d}(N_{c+d,b} tu)\\[3mm]
                      &=& x_{b+c+d}( \beta_{12}\gamma_1 tu).
\end{array}\eqno{(037)}
$$

$$
\begin{array}{rcl}
[x_{c+d}(u),x_{e}(t)] &=& x_{c+d+e}(N_{c+d,e} tu)\\[3mm]
                      &=& x_{c+d+e}( \beta_{13}\gamma_2 tu).
\end{array}\eqno{(038)}
$$

$$
\begin{array}{rcl}
[x_{c+d}(u),x_{b+d+e}(t)] &=& x_{b+c+2d+e}(N_{c+d,b+d+e} tu)\\[3mm]
                          &=& x_{b+c+2d+e}( -\beta_{14}\gamma_3 tu).
\end{array}\eqno{(039)}
$$

$$
\begin{array}{rcl}
[x_{c+d}(u),x_{a+b+c+d+e}(t)] &=& x_{a+b+2c+2d+e}(N_{c+d,a+b+c+d+e} tu)\\[3mm]
                              &=& x_{a+b+2c+2d+e}( \alpha_{56}\beta_{15}\gamma_3  tu).
\end{array}\eqno{(040)}
$$

$$
\begin{array}{rcl}
[x_{c+d}(u),x_{e+f}(t)] &=& x_{c+d+e+f}(N_{c+d,e+f} tu)\\[3mm]
                        &=& x_{c+d+e+f}( \beta_{16}\gamma_4 tu).
\end{array}\eqno{(041)}
$$

$$
\begin{array}{rcl}
[x_{c+d}(u),x_{b+d+e+f}(t)] &=& x_{b+c+2d+e+f}(N_{c+d,b+d+e+f} tu)\\[3mm]
                            &=& x_{b+c+2d+e+f}( -\beta_{17}\gamma_5  tu).
\end{array}\eqno{(042)}
$$

$$
\begin{array}{rcl}
[x_{c+d}(u),x_{a+b+c+d+e+f}(t)] &=& x_{a+b+2c+2d+e+f}(N_{c+d,a+b+c+d+e+f} tu)\\[3mm]
                                &=& x_{a+b+2c+2d+e+f}( \alpha_{89}\beta_{18}\gamma_5  tu).
 \end{array}\eqno{(043)}
$$

$$
\begin{array}{rcl}
[x_{c+d}(u),x_{a+b+c+2d+2e+f}(t)] &=& x_{a+b+2c+3d+2e+f}(N_{c+d,a+b+c+2d+2e+f} tu)\\[3mm]
                                  &=& x_{a+b+2c+3d+2e+f}( -\beta_{19}\gamma_6 tu).
\end{array}\eqno{(044)}
$$

{\boldmath$a+c+d:$}

$$
\begin{array}{rcl}
[x_{a+c+d}(u),x_{b}(t)] &=& x_{a+b+c+d}(N_{a+c+d,b} tu)\\[3mm]
                        &=& x_{a+b+c+d}( \alpha_{23}\beta_{12}\gamma_1 tu).
\end{array}\eqno{(045)}
$$

$$
\begin{array}{rcl}
[x_{a+c+d}(u),x_{e}(t)] &=& x_{a+c+d+e}(N_{a+c+d,e} tu)\\[3mm]
                        &=& x_{a+c+d+e}( \alpha_{24}\beta_{13}\gamma_2 tu).
\end{array}\eqno{(046)}
$$

$$
\begin{array}{rcl}
[x_{a+c+d}(u),x_{b+d+e}(t)] &=& x_{a+b+c+2d+e}(N_{a+c+d,b+d+e} tu)\\[3mm]
                            &=& x_{a+b+c+2d+e}( -\alpha_{26}\beta_{14}\gamma_3 tu).
\end{array}\eqno{(047)}
$$

$$
\begin{array}{rcl}
[x_{a+c+d}(u),x_{b+c+d+e}(t)] &=& x_{a+b+2c+2d+e}(N_{a+c+d,b+c+d+e} tu)\\[3mm]
                              &=& x_{a+b+2c+2d+e}( -\alpha_{26}\beta_{15}\gamma_3 tu).
\end{array}\eqno{(048)}
$$

$$
\begin{array}{rcl}
[x_{a+c+d}(u),x_{e+f}(t)] &=& x_{a+c+d+e+f}(N_{a+c+d,e+f} tu)\\[3mm]
                          &=& x_{a+c+d+e+f}( \alpha_{27}\beta_{16}\gamma_4 tu).
\end{array}\eqno{(049)}
$$

$$
\begin{array}{rcl}
[x_{a+c+d}(u),x_{b+d+e+f}(t)] &=& x_{a+b+c+2d+e+f}(N_{a+c+d,b+d+e+f} tu)\\[3mm] 
                              &=& x_{a+b+c+2d+e+f}( -\alpha_{29}\beta_{17}\gamma_5 tu).
 \end{array}\eqno{(050)}
$$

$$
\begin{array}{rcl}
[x_{a+c+d}(u),x_{b+c+d+e+f}(t)] &=& x_{a+b+2c+2d+e+f}(N_{a+c+d,b+c+d+e+f} tu)\\[3mm]
                                &=& x_{a+b+2c+2d+e+f}( -\alpha_{29}\beta_{18}\gamma_5 tu).
\end{array}\eqno{(051)}
$$

$$
\begin{array}{rcl}
[x_{a+c+d}(u),x_{b+c+2d+2e+f}(t)] &=& x_{a+b+2c+3d+2e+f}(N_{a+c+d,b+c+2d+2e+f} tu)\\[3mm]
                                  &=& x_{a+b+2c+3d+2e+f}( \alpha_{20}\beta_{19}\gamma_6 tu).
\end{array}\eqno{(052)}
$$

{\boldmath$b:$}

$$
\begin{array}{rcl}
[x_b(u),x_{d+e}(t)] &=& x_{b+d+e}(N_{b,d+e} tu)\\[3mm]
                    &=& x_{b+d+e}( \delta_1 tu).
\end{array}\eqno{(053)}
$$

$$
\begin{array}{rcl}
[x_b(u),x_{c+d+e}(t)] &=& x_{b+c+d+e}(N_{b,c+d+e} tu)\\[3mm]
                      &=& x_{b+c+d+e}( \beta_{34}\delta_1 tu).
\end{array}\eqno{(054)}
$$

$$
\begin{array}{rcl}
[x_b(u),x_{a+c+d+e}(t)] &=& x_{a+b+c+d+e}(N_{b,a+c+d+e} tu)\\[3mm]
                        &=& x_{a+b+c+d+e}( \alpha_{45}\beta_{34}\delta_1 tu).
\end{array}\eqno{(055)}
$$

$$
\begin{array}{rcl}
[x_b(u),x_{d+e+f}(t)] &=& x_{b+d+e+f}(N_{b,d+e+f} tu)\\[3mm]
                      &=& x_{b+d+e+f}( \delta_2 tu).
\end{array}\eqno{(056)}
$$

$$
\begin{array}{rcl}
[x_b(u),x_{c+d+e+f}(t)] &=& x_{b+c+d+e+f}(N_{b,c+d+e+f} tu)\\[3mm]
                        &=& x_{b+c+d+e+f}( \beta_{67}\delta_2 tu).
\end{array}\eqno{(057)}
$$

$$
\begin{array}{rcl}
[x_b(u),x_{a+c+d+e+f}(t)] &=& x_{a+b+c+d+e+f}(N_{b,a+c+d+e+f} tu)\\[3mm]
                          &=& x_{a+b+c+d+e+f}( \alpha_{78}\beta_{67}\delta_2 tu).
\end{array}\eqno{(058)}
$$

$$
\begin{array}{rcl}
[x_b(u),x_{a+b+2c+3d+2e+f}(t)] &=& x_{a+2b+2c+3d+2e+f}(N_{b,a+b+2c+3d+2e+f} tu)\\[3mm]
                               &=& x_{a+2b+2c+3d+2e+f}( \delta_3 tu).
\end{array}\eqno{(059)}
$$

{\boldmath$b+d:$}

$$
\begin{array}{rcl}
[x_{b+d}(u),x_{e}(t)] &=& x_{b+d+e}(N_{b+d,e} tu)\\[3mm]
                      &=& x_{b+d+e}( -\gamma_{12}\delta_1 tu).
\end{array}\eqno{(060)}
$$

$$
\begin{array}{rcl}
[x_{b+d}(u),x_{c+d+e}(t)] &=& x_{b+c+2d+e}(N_{b+d,c+d+e} tu)\\[3mm]
                          &=& x_{b+c+2d+e}( \beta_{34}\gamma_{13}\delta_1 tu).
\end{array}\eqno{(061)}
$$

$$
\begin{array}{rcl}
[x_{b+d}(u),x_{a+c+d+e}(t)] &=& x_{a+b+c+2d+e}(N_{b+d,a+c+d+e} tu)\\[3mm]
                            &=& x_{a+b+c+2d+e}( \alpha_{46}\beta_{34}\gamma_{13}\delta_1 tu).
\end{array}\eqno{(062)}
$$

$$
\begin{array}{rcl}
[x_{b+d}(u),x_{e+f}(t)] &=& x_{b+d+e+f}(N_{b+d,e+f} tu)\\[3mm]
                        &=& x_{b+d+e+f}( -\gamma_{14}\delta_2 tu).
\end{array}\eqno{(063)}
$$

$$
\begin{array}{rcl}
[x_{b+d}(u),x_{c+d+e+f}(t)] &=& x_{b+c+2d+e+f}(N_{b+d,c+d+e+f} tu)\\[3mm]
                            &=& x_{b+c+2d+e+f}( \beta_{67}\gamma_{15}\delta_2 tu).
\end{array} \eqno{(064)}
$$

$$
\begin{array}{rcl}
[x_{b+d}(u),x_{a+c+d+e+f}(t)] &=& x_{a+b+c+2d+e+f}(N_{b+d,a+c+d+e+f} tu)\\[3mm]
                              &=& x_{a+b+c+2d+e+f}( \alpha_{79}\beta_{67}\gamma_{15}\delta_2 tu).
\end{array} \eqno{(065)}
$$

$$
\begin{array}{rcl}
[x_{b+d}(u),x_{a+b+2c+2d+2e+f}(t)] &=& x_{a+2b+2c+3d+2e+f}(N_{b+d,a+b+2c+2d+2e+f} tu)\\[3mm]
                                   &=& x_{a+2b+2c+3d+2e+f}( -\gamma_{16}\delta_3 tu).
\end{array} \eqno{(066)}
$$

{\boldmath$b+c+d:$}

$$
\begin{array}{rcl}
[x_{b+c+d}(u),x_{e}(t)] &=& x_{b+c+d+e}(N_{b+c+d,e} tu)\\[3mm]
                        &=& x_{b+c+d+e}( -\beta_{24}\gamma_{12}\delta_1 tu).
\end{array}\eqno{(067)}
$$

$$
\begin{array}{rcl}
[x_{b+c+d}(u),x_{d+e}(t)] &=& x_{b+c+2d+e}(N_{b+c+d,d+e} tu)\\[3mm]
                          &=& x_{b+c+2d+e}( -\beta_{24}\gamma_{13}\delta_1 tu).
\end{array}\eqno{(068)}
$$

$$
\begin{array}{rcl}
[x_{b+c+d}(u),x_{a+c+d+e}(t)] &=& x_{a+b+2c+2d+e}(N_{b+c+d,a+c+d+e} tu)\\[3mm]
                              &=& x_{a+b+2c+2d+e}( \alpha_{46}\beta_{2345}\gamma_{13}\delta_1 tu).
\end{array}\eqno{(069)}
$$

$$
\begin{array}{rcl}
[x_{b+c+d}(u),x_{e+f}(t)] &=& x_{b+c+d+e+f}(N_{b+c+d,e+f} tu)\\[3mm]
                          &=& x_{b+c+d+e+f}( -\beta_{27}\gamma_{14}\delta_2 tu).
\end{array}\eqno{(070)}
$$

$$
\begin{array}{rcl}
[x_{b+c+d}(u),x_{d+e+f}(t)] &=& x_{b+c+2d+e+f}(N_{b+c+d,d+e+f} tu)\\[3mm]
                            &=& x_{b+c+2d+e+f}( -\beta_{27}\gamma_{15}\delta_2 tu).
\end{array}\eqno{(071)}
$$

$$
\begin{array}{rcl}
[x_{b+c+d}(u),x_{a+c+d+e+f}(t)] &=& x_{a+b+2c+2d+e+f}(N_{b+c+d,a+c+d+e+f} tu)\\[3mm]
                                &=& x_{a+b+2c+2d+e+f}( \alpha_{79}\beta_{2678}\gamma_{15}\delta_2 tu).
\end{array}\eqno{(072)}
$$

$$
\begin{array}{rcl}
[x_{b+c+d}(u),x_{a+b+c+2d+2e+f}(t)] &=& x_{a+2b+2c+3d+2e+f}(N_{b+c+d,a+b+c+2d+2e+f} tu)\\[3mm]
                                    &=& x_{a+2b+2c+3d+2e+f}( \beta_{29}\gamma_{16}\delta_3 tu).
\end{array}\eqno{(073)}
$$

{\boldmath$a+b+c+d:$}

$$
\begin{array}{rcl}
[x_{a+b+c+d}(u),x_{e}(t)] &=& x_{a+b+c+d+e}(N_{a+b+c+d,e} tu)\\[3mm]
                          &=&x_{a+b+c+d+e}( -\alpha_{35}\beta_{24}\gamma_{12}\delta_1 tu).
\end{array}\eqno{(074)}
$$

$$
\begin{array}{rcl}
[x_{a+b+c+d}(u),x_{d+e}(t)] &=& x_{a+b+c+2d+e}(N_{a+b+c+d,d+e} tu)\\[3mm]
                            &=& x_{a+b+c+2d+e}( -\alpha_{36}\beta_{24}\gamma_{13}\delta_1 tu).
\end{array}\eqno{(075)}
$$

$$
\begin{array}{rcl}
[x_{a+b+c+d}(u),x_{c+d+e}(t)] &=& x_{a+b+2c+2d+e}(N_{a+b+c+d,c+d+e} tu)\\[3mm]
                              &=& x_{a+b+2c+2d+e}( -\alpha_{36}\beta_{2345}\gamma_{13}\delta_1 tu).
\end{array}\eqno{(076)}
$$

$$
\begin{array}{rcl}
[x_{a+b+c+d}(u),x_{e+f}(t)] &=& x_{a+b+c+d+e+f}(N_{a+b+c+d,e+f} tu)\\[3mm]
                            &=& x_{a+b+c+d+e+f}( -\alpha_{38}\beta_{27}\gamma_{14}\delta_2 tu).
\end{array}\eqno{(077)}
$$

$$
\begin{array}{rcl}
[x_{a+b+c+d}(u),x_{d+e+f}(t)] &=& x_{a+b+c+2d+e+f}(N_{a+b+c+d,d+e+f} tu)\\[3mm]
                              &=& x_{a+b+c+2d+e+f}( -\alpha_{39}\beta_{27}\gamma_{15}\delta_2 tu).
\end{array}\eqno{(078)}
$$

$$
\begin{array}{rcl}
[x_{a+b+c+d}(u),x_{c+d+e+f}(t)] &=& x_{a+b+2c+2d+e+f}(N_{a+b+c+d,c+d+e+f} tu)\\[3mm]
                                &=& x_{a+b+2c+2d+e+f}( -\alpha_{39}\beta_{2678}\gamma_{15}\delta_2 tu).
\end{array}\eqno{(079)}
$$

$$
\begin{array}{rcl}
[x_{a+b+c+d}(u),x_{b+c+2d+2e+f}(t)] &=& x_{a+2b+2c+3d+2e+f}(N_{a+b+c+d,b+c+2d+2e+f} tu)\\[3mm]
                                    &=& x_{a+2b+2c+3d+2e+f}( -\alpha_{30}\beta_{29}\gamma_{16}\delta_3 tu).
\end{array}\eqno{(080)}
$$

{\boldmath$e:$}

$$
\begin{array}{rcl}
[x_e(u),x_{f}(t)] &=& x_{e+f}(N_{e,f} tu)\\[3mm]
                  &=& x_{e+f}( \epsilon_1 tu).
\end{array}\eqno{(081)}
$$

$$
\begin{array}{rcl}
[x_e(u),x_{b+c+2d+e+f}(t)] &=& x_{b+c+2d+2e+f}(N_{e,b+c+2d+e+f} tu)\\[3mm]
                           &=& x_{b+c+2d+2e+f}( \epsilon_2 tu).
\end{array}\eqno{(082)}
$$

$$
\begin{array}{rcl}
[x_e(u),x_{a+b+c+2d+e+f}(t)] &=& x_{a+b+c+2d+2e+f}(N_{e,a+b+c+2d+e+f} tu)\\[3mm]
                             &=& x_{a+b+c+2d+2e+f}( \alpha_{90}\epsilon_2 tu).
\end{array}\eqno{(083)}
$$

$$
\begin{array}{rcl}
[x_e(u),x_{a+b+2c+2d+e+f}(t)] &=& x_{a+b+2c+2d+2e+f}(N_{e,a+b+2c+2d+e+f} tu)\\[3mm]
                              &=& x_{a+b+2c+2d+2e+f}( \alpha_{90}\beta_{89}\epsilon_2tu).
\end{array}\eqno{(084)}
$$

{\boldmath$d+e:$}

$$
\begin{array}{rcl}
[x_{d+e}(u),x_{f}(t)] &=& x_{d+e+f}(N_{d+e,f} tu)\\[3mm]
                      &=& x_{d+e+f}( \gamma_{24}\epsilon_1 tu).
\end{array}\eqno{(085)}
$$

$$
\begin{array}{rcl}
[x_{d+e}(u),x_{b+c+d+e+f}(t)] &=& x_{b+c+2d+2e+f}(N_{d+e,b+c+d+e+f} tu)\\[3mm]
                              &=& x_{b+c+2d+2e+f}( -\gamma_{25}\epsilon_2 tu).
\end{array}\eqno{(086)}
$$

$$
\begin{array}{rcl}
[x_{d+e}(u),x_{a+b+c+d+e+f}(t)] &=& x_{a+b+c+2d+2e+f}(N_{d+e,a+b+c+d+e+f} tu)\\[3mm]
                                &=& x_{a+b+c+2d+2e+f}( -\alpha_{80}\gamma_{25}\epsilon_2 tu).
\end{array}\eqno{(087)}
$$

$$
\begin{array}{rcl}
[x_{d+e}(u),x_{a+b+2c+2d+e+f}(t)] &=& x_{a+b+2c+3d+2e+f}(N_{d+e,a+b+2c+2d+e+f} tu)\\[3mm]
                                  &=& x_{a+b+2c+3d+2e+f}( \alpha_{90}\beta_{89}\gamma_{26}\epsilon_2 tu).
\end{array}\eqno{(088)}
$$

{\boldmath$c+d+e:$}

$$
\begin{array}{rcl}
[x_{c+d+e}(u),x_{f}(t)] &=& x_{c+d+e+f}(N_{c+d+e,f} tu)\\[3mm]
                        &=& x_{c+d+e+f}( \beta_{36}\gamma_{24}\epsilon_1 tu).
\end{array}\eqno{(089)}
$$

$$
\begin{array}{rcl}
[x_{c+d+e}(u),x_{b+d+e+f}(t)] &=& x_{b+c+2d+2e+f}(N_{c+d+e,b+d+e+f} tu)\\[3mm]
                              &=& x_{b+c+2d+2e+f}( \beta_{37}\gamma_{25}\epsilon_2 tu).
\end{array}\eqno{(090)}
$$

$$
\begin{array}{rcl}
[x_{c+d+e}(u),x_{a+b+c+d+e+f}(t)] &=& x_{a+b+2c+2d+2e+f}(N_{c+d+e,a+b+c+d+e+f} tu)\\[3mm]
                                  &=& x_{a+b+2c+2d+2e+f}( -\alpha_{80}\beta_{39}\gamma_{25}\epsilon_2 tu).
\end{array}\eqno{(091)}
$$

$$
\begin{array}{rcl}
[x_{c+d+e}(u),x_{a+b+c+2d+e+f}(t)] &=& x_{a+b+2c+3d+2e+f}(N_{c+d+e,a+b+c+2d+e+f} tu)\\[3mm]
                                   &=& x_{a+b+2c+3d+2e+f}( -\alpha_{90}\beta_{39}\gamma_{26}\epsilon_2tu).
\end{array}\eqno{(092)}
$$

{\boldmath$a+c+d+e:$}

$$
\begin{array}{rcl}
[x_{a+c+d+e}(u),x_{f}(t)] &=& x_{a+c+d+e+f}(N_{a+c+d+e,f} tu)\\[3mm]
                          &=& x_{a+c+d+e+f}( \alpha_{47}\beta_{36}\gamma_{24}\epsilon_1 tu).
\end{array}\eqno{(093)}
$$

$$
\begin{array}{rcl}
[x_{a+c+d+e}(u),x_{b+d+e+f}(t)] &=& x_{a+b+c+2d+2e+f}(N_{a+c+d+e,b+d+e+f} tu)\\[3mm]
                                &=& x_{a+b+c+2d+2e+f}( \alpha_{40}\beta_{37}\gamma_{25}\epsilon_2 tu).
\end{array}\eqno{(094)}
$$

$$
\begin{array}{rcl}
[x_{a+c+d+e}(u),x_{b+c+d+e+f}(t)] &=& x_{a+b+2c+2d+2e+f}(N_{a+c+d+e,b+c+d+e+f} tu)\\[3mm]
                                  &=& x_{a+b+2c+2d+2e+f}( \alpha_{40}\beta_{39}\gamma_{25}\epsilon_2 tu).
\end{array}\eqno{(095)}
$$

$$
\begin{array}{rcl}
[x_{a+c+d+e}(u),x_{b+c+2d+e+f}(t)] &=& x_{a+b+2c+3d+2e+f}(N_{a+c+d+e,b+c+2d+e+f} tu)\\[3mm]
                                   &=& x_{a+b+2c+3d+2e+f}( \alpha_{40}\beta_{39}\gamma_{26}\epsilon_2 tu).
\end{array}\eqno{(096)}
$$

{\boldmath$b+d+e:$}

$$
\begin{array}{rcl}
[x_{b+d+e}(u),x_{f}(t)] &=& x_{b+d+e+f}(N_{b+d+e,f} tu)\\[3mm]
                        &=& x_{b+d+e+f}( \gamma_{24}\delta_{12}\epsilon_1 tu).
\end{array}\eqno{(097)}
$$

$$
\begin{array}{rcl}
[x_{b+d+e}(u),x_{c+d+e+f}(t)] &=& x_{b+c+2d+2e+f}(N_{b+d+e,c+d+e+f} tu)\\[3mm]
                              &=& x_{b+c+2d+2e+f}( \beta_{67}\gamma_{25}\delta_{12}\epsilon_2 tu).
\end{array}\eqno{(098)}
$$

$$
\begin{array}{rcl}
[x_{b+d+e}(u),x_{a+c+d+e+f}(t)] &=& x_{a+b+c+2d+2e+f}(N_{b+d+e,a+c+d+e+f} tu)\\[3mm]
                                &=& x_{a+b+c+2d+2e+f}( \alpha_{70}\beta_{67}\gamma_{25}\delta_{12}\epsilon_2 tu).
\end{array}\eqno{(099)}
$$

$$
\begin{array}{rcl}
[x_{b+d+e}(u),x_{a+b+2c+2d+e+f}(t)] &=& x_{a+2b+2c+3d+2e+f}(N_{b+d+e,a+b+2c+2d+e+f} tu)\\[3mm]
                                    &=& x_{a+2b+2c+3d+2e+f}( \alpha_{90}\beta_{89}\gamma_{26}\delta_{13}\epsilon_2 tu).
\end{array}\eqno{(100)}
$$

{\boldmath$b+c+d+e:$}

$$
\begin{array}{rcl}
[x_{b+c+d+e}(u),x_{f}(t)] &=& x_{b+c+d+e+f}(N_{b+c+d+e,f} tu)\\[3mm] 
                          &=& x_{b+c+d+e+f}( \beta_{47}\gamma_{24}\delta_{12}\epsilon_1 tu).
\end{array}\eqno{(101)}
$$

$$
\begin{array}{rcl}
[x_{b+c+d+e}(u),x_{d+e+f}(t)] &=& x_{b+c+2d+2e+f}(N_{b+c+d+e,d+e+f} tu)\\[3mm]
                              &=& x_{b+c+2d+2e+f}( -\beta_{47}\gamma_{25}\delta_{12}\epsilon_2 tu).
\end{array}\eqno{(102)}
$$

$$
\begin{array}{rcl}
[x_{b+c+d+e}(u),x_{a+c+d+e+f}(t)] &=& x_{a+b+2c+2d+2e+f}(N_{b+c+d+e,a+c+d+e+f} tu)\\[3mm]
                                  &=& x_{a+b+2c+2d+2e+f}( \alpha_{70}\beta_{4679}\gamma_{25}\delta_{12}\epsilon_2 tu).
\end{array}\eqno{(103)}
$$

$$
\begin{array}{rcl}
[x_{b+c+d+e}(u),x_{a+b+c+2d+e+f}(t)] &=& x_{a+2b+2c+3d+2e+f}(N_{b+c+d+e,a+b+c+2d+e+f} tu)\\[3mm]
                                     &=& x_{a+2b+2c+3d+2e+f}( -\alpha_{90}\beta_{49}\gamma_{26}\delta_{13}\epsilon_2 tu).
\end{array}\eqno{(104)}
$$

{\boldmath$a+b+c+d+e:$}

$$
\begin{array}{rcl}
[x_{a+b+c+d+e}(u),x_{f}(t)] &=& x_{a+b+c+d+e+f}(N_{a+b+c+d+e,f} tu)\\[3mm]
                            &=& x_{a+b+c+d+e+f}( \alpha_{58}\beta_{47}\gamma_{24}\delta_{12}\epsilon_1 tu).
\end{array}\eqno{(105)}
$$

$$
\begin{array}{rcl}
[x_{a+b+c+d+e}(u),x_{d+e+f}(t)] &=& x_{a+b+c+2d+2e+f}(N_{a+b+c+d+e,d+e+f} tu)\\[3mm]
                                &=& x_{a+b+c+2d+2e+f}( -\alpha_{50}\beta_{47}\gamma_{25}\delta_{12}\epsilon_2 tu).
\end{array}\eqno{(106)}
$$

$$
\begin{array}{rcl}
[x_{a+b+c+d+e}(u),x_{c+d+e+f}(t)] &=& x_{a+b+2c+2d+2e+f}(N_{a+b+c+d+e,c+d+e+f} tu)\\[3mm]
                                  &=& x_{a+b+2c+2d+2e+f}( -\alpha_{50}\beta_{4679}\gamma_{25}\delta_{12}\epsilon_2 tu).
\end{array}\eqno{(107)}
$$

$$
\begin{array}{rcl}
[x_{a+b+c+d+e}(u),x_{b+c+2d+e+f}(t)] &=& x_{a+2b+2c+3d+2e+f}(N_{a+b+c+d+e,b+c+2d+e+f} tu)\\[3mm]
                                     &=& x_{a+2b+2c+3d+2e+f}(\alpha_{50}\beta_{49}\gamma_{26}\delta_{13}\epsilon_2 tu).
\end{array}\eqno{(108)}
$$

{\boldmath$b+c+2d+e:$}

$$
\begin{array}{rcl}
[x_{b+c+2d+e}(u),x_{f}(t)] &=& x_{b+c+2d+e+f}(N_{b+c+2d+e,f} tu)\\[3mm]
                           &=& x_{b+c+2d+e+f}( \beta_{47}\gamma_{2345}\delta_{12}\epsilon_1 tu).
\end{array}\eqno{(109)}
$$

$$
\begin{array}{rcl}
[x_{b+c+2d+e}(u),x_{e+f}(t)] &=& x_{b+c+2d+2e+f}(N_{b+c+2d+e,e+f} tu)\\[3mm]
                             &=& x_{b+c+2d+2e+f}( \beta_{47}\gamma_{2345}\delta_{12}\epsilon_2 tu).
\end{array}\eqno{(110)}
$$

$$
\begin{array}{rcl}
[x_{b+c+2d+e}(u),x_{a+c+d+e+f}(t)] &=& x_{a+b+2c+3d+2e+f}(N_{b+c+2d+e,a+c+d+e+f} tu)\\[3mm]
                                   &=& x_{a+b+2c+3d+2e+f}( \alpha_{70}\beta_{4679}\gamma_{2356}\delta_{12}\epsilon_2 tu).
\end{array}\eqno{(111)}
$$

$$
\begin{array}{rcl}
[x_{b+c+2d+e}(u),x_{a+b+c+d+e+f}(t)] &=& x_{a+2b+2c+3d+2e+f}(N_{b+c+2d+e,a+b+c+d+e+f} tu)\\[3mm]
                                     &=& x_{a+2b+2c+3d+2e+f}( \alpha_{80}\beta_{49}\gamma_{2356}\delta_{13}\epsilon_2 tu).
\end{array}\eqno{(112)}
$$

{\boldmath$a+b+c+2d+e:$}

$$
\begin{array}{rcl}
[x_{a+b+c+2d+e}(u),x_{f}(t)] &=& x_{a+b+c+2d+e+f}(N_{a+b+c+2d+e,f} tu)\\[3mm]
                             &=& x_{a+b+c+2d+e+f}( \alpha_{69}\beta_{47}\gamma_{2345}\delta_{12}\epsilon_1 tu).
\end{array}\eqno{(113)}
$$

$$
\begin{array}{rcl}
[x_{a+b+c+2d+e}(u),x_{e+f}(t)] &=& x_{a+b+c+2d+2e+f}(N_{a+b+c+2d+e,e+f} tu)\\[3mm]
                               &=& x_{a+b+c+2d+2e+f}( \alpha_{60}\beta_{47}\gamma_{2345}\delta_{12}\epsilon_2 tu).
\end{array}\eqno{(114)}
$$

$$
\begin{array}{rcl}
[x_{a+b+c+2d+e}(u),x_{c+d+e+f}(t)] &=& x_{a+b+2c+3d+2e+f}(N_{a+b+c+2d+e,c+d+e+f} tu)\\[3mm]
                                   &=& x_{a+b+2c+3d+2e+f}( -\alpha_{60}\beta_{4679}\gamma_{2356}\delta_{12}\epsilon_2 tu).
\end{array}\eqno{(115)}
$$

$$
\begin{array}{rcl}
[x_{a+b+c+2d+e}(u),x_{b+c+d+e+f}(t)] &=& x_{a+2b+2c+3d+2e+f}(N_{a+b+c+2d+e,b+c+d+e+f} tu)\\[3mm]
                                     &=& x_{a+2b+2c+3d+2e+f}( - \alpha_{60}\beta_{49}\gamma_{2356}\delta_{13}\epsilon_2 tu).
\end{array}\eqno{(116)}
$$

{\boldmath$a+b+2c+2d+e:$}

$$
\begin{array}{rcl}
[x_{a+b+2c+2d+e}(u),x_{f}(t)] &=& x_{a+b+2c+2d+e+f}(N_{a+b+2c+2d+e,f} tu)\\[3mm]
                              &=& x_{a+b+2c+2d+e+f}( \alpha_{69}\beta_{4578}\gamma_{2345}\delta_{12}\epsilon_1 tu).
\end{array}\eqno{(117)}
$$

$$
\begin{array}{rcl}
[x_{a+b+2c+2d+e}(u),x_{e+f}(t)] &=& x_{a+b+2c+2d+2e+f}(N_{a+b+2c+2d+e,e+f} tu)\\[3mm]
                                &=& x_{a+b+2c+2d+2e+f}( \alpha_{60}\beta_{4579}\gamma_{2345}\delta_{12}\epsilon_2 tu).
\end{array}\eqno{(118)}
$$

$$
\begin{array}{rcl}
[x_{a+b+2c+2d+e}(u),x_{d+e+f}(t)] &=& x_{a+b+2c+3d+2e+f}(N_{a+b+2c+2d+e,d+e+f} tu)\\[3mm]
                                  &=& x_{a+b+2c+3d+2e+f}( \alpha_{60}\beta_{4579}\gamma_{2356}\delta_{12}\epsilon_2 tu).
\end{array}\eqno{(119)}
$$

$$
\begin{array}{rcl}
[x_{a+b+2c+2d+e}(u),x_{b+d+e+f}(t)] &=& x_{a+2b+2c+3d+2e+f}(N_{a+b+2c+2d+e,b+d+e+f} tu)\\[3mm]
                                    &=& x_{a+2b+2c+3d+2e+f}( \alpha_{60}\beta_{4579}\gamma_{2356}\delta_{13}\epsilon_2 tu).
\end{array}\eqno{(120)}
$$

\vskip3mm
\centerline{\textit{Negative roots}}

{\boldmath$-a:$}

$$
\begin{array}{rcl} 
[x_{-a}(u),x_{-c}(t)] &=& x_{-a-c}(-N_{a,c} tu)\\[3mm] 
                      &=& x_{-a-c}( -\alpha_1 tu).
\end{array} \eqno{(001)}
$$

$$
\begin{array}{rcl} 
[x_{-a}(u),x_{-c-d}(t)] &=& x_{-a-c-d}(-N_{a,c+d} tu)\\[3mm]
                        &=& x_{-a-c-d}( -\alpha_2 tu).
\end{array} \eqno{(002)}
$$

$$
\begin{array}{rcl} 
[x_{-a}(u),x_{-b-c-d}(t)] &=& x_{-a-b-c-d}(-N_{a,b+c+d} tu)\\[3mm]
                          &=& x_{-a-b-c-d}( -\alpha_3 tu).
\end{array} \eqno{(003)}
$$

$$
\begin{array}{rcl} 
[x_{-a}(u),x_{-c-d-e}(t)] &=& x_{-a-c-d-e}(-N_{a,c+d+e} tu)\\[3mm]
                          &=& x_{-a-c-d-e}( -\alpha_4 tu).
\end{array} \eqno{(004)}
$$

$$
\begin{array}{rcl} 
[x_{-a}(u),x_{-b-c-d-e}(t)] &=& x_{-a-b-c-d-e}(-N_{a,b+c+d+e} tu)\\[3mm]
                            &=& x_{-a-b-c-d-e}( -\alpha_5 tu).
\end{array} \eqno{(005)}
$$

$$
\begin{array}{rcl} 
[x_{-a}(u),x_{-b-c-2d-e}(t)] &=& x_{-a-b-c-2d-e}(-N_{a,b+c+2d+e} tu)\\[3mm]
                             &=& x_{-a-b-c-2d-e}( -\alpha_6 tu).
\end{array} \eqno{(006)}
$$

$$
\begin{array}{rcl} 
[x_{-a}(u),x_{-c-d-e-f}(t)] &=& x_{-a-c-d-e-f}(-N_{a,c+d+e+f} tu)\\[3mm]
                            &=& x_{-a-c-d-e-f}( -\alpha_7 tu).
\end{array} \eqno{(007)}
$$

$$
\begin{array}{rcl} 
[x_{-a}(u),x_{-b-c-d-e-f}(t)] &=& x_{-a-b-c-d-e-f}(-N_{a,b+c+d+e+f} tu)\\[3mm] 
                              &=& x_{-a-b-c-d-e-f}( -\alpha_8 tu).
\end{array} \eqno{(008)}
$$

$$
\begin{array}{rcl} 
[x_{-a}(u),x_{-b-c-2d-e-f}(t)] &=& x_{-a-b-c-2d-e-f}(-N_{a,b+c+2d+e+f} tu) \\[3mm]
                               &=& x_{-a-b-c-2d-e-f}( -\alpha_9 tu).
\end{array} \eqno{(009)}
$$

$$
\begin{array}{rcl} 
[x_{-a}(u),x_{-b-c-2d-2e-f}(t)] &=& x_{-a-b-c-2d-2e-f}(-N_{a,b+c+2d+2e+f} tu)\\[3mm]
                                &=& x_{-a-b-c-2d-2e-f}( -\alpha_0 tu).
\end{array} \eqno{(010)}
$$

{\boldmath$-c:$}

$$
\begin{array}{rcl} 
[x_{-c}(u),x_{-d}(t)] &=& x_{-c-d}(-N_{c,d} tu)\\[3mm] 
                      &=& x_{-c-d}( -\beta_1 tu).
\end{array} \eqno{(011)}
$$

$$
\begin{array}{rcl} 
[x_{-c}(u),x_{-b-d}(t)] &=& x_{-b-c-d}(-N_{c,b+d} tu)\\[3mm] 
                        &=& x_{-b-c-d}( -\beta_2 tu).
\end{array} \eqno{(012)}
$$

$$
\begin{array}{rcl} 
[x_{-c}(u),x_{-d-e}(t)] &=& x_{-c-d-e}(-N_{c,d+e} tu)\\[3mm] 
                        &=& x_{-c-d-e}( -\beta_3 tu).
\end{array} \eqno{(013)}
$$

$$
\begin{array}{rcl} 
[x_{-c}(u),x_{-b-d-e}(t)] &=& x_{-b-c-d-e}(-N_{c,b+d+e} tu)\\[3mm] 
                          &=& x_{-b-c-d-e}( -\beta_4 tu).
\end{array} \eqno{(014)}
$$

$$
\begin{array}{rcl} 
[x_{-c}(u),x_{-a-b-c-2d-e}(t)] &=& x_{-a-b-2c-2d-e}(-N_{c,a+b+c+2d+e} tu)\\[3mm] 
                               &=& x_{-a-b-2c-2d-e}( -\beta_5 tu).
\end{array} \eqno{(015)}
$$

$$
\begin{array}{rcl} 
[x_{-c}(u),x_{-d-e-f}(t)] &=& x_{-c-d-e-f}(-N_{c,d+e+f} tu)\\[3mm] 
                          &=& x_{-c-d-e-f}( -\beta_6 tu).
\end{array} \eqno{(016)}
$$

$$
\begin{array}{rcl} 
[x_{-c}(u),x_{-b-d-e-f}(t)] &=& x_{-b-c-d-e-f}(-N_{c,b+d+e+f} tu)\\[3mm] 
                            &=& x_{-b-c-d-e-f}( -\beta_7 tu).
\end{array} \eqno{(017)}
$$

$$
\begin{array}{rcl} 
[x_{-c}(u),x_{-a-b-c-2d-e-f}(t)] &=& x_{-a-b-2c-2d-e-f}(-N_{c,a+b+c+2d+e+f} tu)\\[3mm] 
                                 &=& x_{-a-b-2c-2d-e-f}( -\beta_8 tu).
\end{array} \eqno{(018)}
$$

$$
\begin{array}{rcl} 
[x_{-c}(u),x_{-a-b-c-2d-2e-f}(t)] &=& x_{-a-b-2c-2d-2e-f}(-N_{c,a+b+c+2d+2e+f} tu)\\[3mm] 
                                  &=& x_{-a-b-2c-2d-2e-f}( -\beta_9 tu).
\end{array} \eqno{(019)}
$$

{\boldmath$-a-c:$}

$$
\begin{array}{rcl} 
[x_{-a-c}(u),x_{-d}(t)] &=& x_{-a-c-d}(-N_{a+c,d} tu)\\[3mm] 
                        &=& x_{-a-c-d}( -\alpha_{12}\beta_1 tu).
\end{array} \eqno{(020)}
$$

$$
\begin{array}{rcl} 
[x_{-a-c}(u),x_{-b-d}(t)] &=& x_{-a-b-c-d}(-N_{a+c,b+d} tu)\\[3mm] 
                          &=& x_{-a-b-c-d}( -\alpha_{13}\beta_2 tu).
\end{array} \eqno{(021)}
$$

$$
\begin{array}{rcl} 
[x_{-a-c}(u),x_{-d-e}(t)] &=& x_{-a-c-d-e}(-N_{a+c,d+e} tu)\\[3mm] 
                          &=& x_{-a-c-d-e}( -\alpha_{14}\beta_3 tu).
\end{array} \eqno{(022)}
$$

$$
\begin{array}{rcl} 
[x_{-a-c}(u),x_{-b-d-e}(t)] &=& x_{-a-b-c-d-e}(-N_{a+c,b+d+e} tu)\\[3mm] 
                            &=& x_{-a-b-c-d-e}( -\alpha_{15}\beta_4 tu).
\end{array} \eqno{(023)}
$$

$$
\begin{array}{rcl} 
[x_{-a-c}(u),x_{-b-c-2d-e}(t)] &=& x_{-a-b-2c-2d-e}(-N_{a+c,b+c+2d+e} tu)\\[3mm] 
                               &=& x_{-a-b-2c-2d-e}( \alpha_{16}\beta_5 tu).
\end{array} \eqno{(024)}
$$

$$
\begin{array}{rcl} 
[x_{-a-c}(u),x_{-d-e-f}(t)] &=& x_{-a-c-d-e-f}(-N_{a+c,d+e+f} tu)\\[3mm] 
                            &=& x_{-a-c-d-e-f}( -\alpha_{17}\beta_6 tu).
\end{array} \eqno{(025)}
$$

$$
\begin{array}{rcl} 
[x_{-a-c}(u),x_{-b-d-e-f}(t)] &=& x_{-a-b-c-d-e-f}(-N_{a+c,b+d+e+f} tu)\\[3mm] 
                              &=& x_{-a-b-c-d-e-f}( -\alpha_{18}\beta_7 tu).
\end{array} \eqno{(026)}
$$

$$
\begin{array}{rcl} 
[x_{-a-c}(u),x_{-b-c-2d-e-f}(t)] &=& x_{-a-b-2c-2d-e-f}(-N_{a+c,b+c+2d+e+f} tu)\\[3mm] 
                                 &=& x_{-a-b-2c-2d-e-f}( \alpha_{19}\beta_8 tu).
 \end{array} \eqno{(027)}
$$

$$
\begin{array}{rcl}
 [x_{-a-c}(u),x_{-b-c-2d-2e-f}(t)] &=& x_{-a-b-2c-2d-2e-f}(-N_{a+c,b+c+2d+2e+f} tu)\\[3mm]  
                                   &=& x_{-a-b-2c-2d-2e-f}( \alpha_{10}\beta_9 tu).
\end{array} \eqno{(028)}
$$

{\boldmath$-d:$}

$$
\begin{array}{rcl} 
[x_{-d}(u),x_{-b}(t)] &=& x_{-b-d}(-N_{d,b} tu)\\[3mm] 
                      &=& x_{-b-d}( -\gamma_1 tu).
\end{array} \eqno{(029)}
$$

$$
\begin{array}{rcl} 
[x_{-d}(u),x_{-e}(t)] &=& x_{-d-e}(-N_{d,e} tu)\\[3mm] 
                      &=& x_{-d-e}( -\gamma_2 tu).
\end{array} \eqno{(030)}
$$

$$
\begin{array}{rcl} 
[x_{-d}(u),x_{-b-c-d-e}(t)] &=& x_{-b-c-2d-e}(-N_{d,b+c+d+e} tu)\\[3mm] 
                            &=& x_{-b-c-2d-e}( -\gamma_3 tu).
\end{array} \eqno{(031)}
$$

$$
\begin{array}{rcl} 
[x_{-d}(u),x_{-a-b-c-d-e}(t)] &=& x_{-a-b-c-2d-e}(-N_{d,a+b+c+d+e} tu)\\[3mm] 
                              &=& x_{-a-b-c-2d-e}(  -\alpha_{56}\gamma_3 tu).
\end{array} \eqno{(032)}
$$

$$
\begin{array}{rcl} 
[x_{-d}(u),x_{-e-f}(t)] &=& x_{-d-e-f}(-N_{d,e+f} tu)\\[3mm] 
                        &=& x_{-d-e-f}( -\gamma_4 tu).
\end{array} \eqno{(033)}
$$

$$
\begin{array}{rcl}
[x_{-d}(u),x_{-b-c-d-e-f}(t)] &=& x_{-b-c-2d-e-f}(-N_{d,b+c+d+e+f} tu)\\[3mm] 
                              &=& x_{-b-c-2d-e-f}( -\gamma_5 tu).
\end{array} \eqno{(034)}
$$

$$
\begin{array}{rcl} 
[x_{-d}(u),x_{-a-b-c-d-e-f}(t)] &=& x_{-a-b-c-2d-e-f}(-N_{d,a+b+c+d+e+f} tu)\\[3mm] 
                                &=& x_{-a-b-c-2d-e-f}( -\alpha_{89}\gamma_5 tu).\end{array} \eqno{(035)}
$$

$$
\begin{array}{rcl} 
[x_{-d}(u),x_{-a-b-2c-2d-2e-f}(t)] &=& x_{-a-b-2c-3d-2e-f}(-N_{d,a+b+2c+2d+2e+f} tu)\\[3mm] 
                                   &=& x_{-a-b-2c-3d-2e-f}(  -\gamma_6 tu).
 \end{array} \eqno{(036)}
$$

{\boldmath$-c-d:$}

$$
\begin{array}{rcl} 
[x_{-c-d}(u),x_{-b}(t)] &=& x_{-b-c-d}(-N_{c+d,b} tu)\\[3mm] 
                        &=& x_{-b-c-d}( -\beta_{12}\gamma_1 tu).
\end{array} \eqno{(037)}
$$

$$
\begin{array}{rcl} 
[x_{-c-d}(u),x_{-e}(t)] &=& x_{-c-d-e}(-N_{c+d,e} tu)\\[3mm] 
                        &=& x_{-c-d-e}( -\beta_{13}\gamma_2 tu).
\end{array} \eqno{(038)}
$$

$$
\begin{array}{rcl} 
[x_{-c-d}(u),x_{-b-d-e}(t)] &=& x_{-b-c-2d-e}(-N_{c+d,b+d+e} tu)\\[3mm] 
                            &=& x_{-b-c-2d-e}( \beta_{14}\gamma_3 tu).
\end{array} \eqno{(039)}
$$

$$
\begin{array}{rcl} 
[x_{-c-d}(u),x_{-a-b-c-d-e}(t)] &=& x_{-a-b-2c-2d-e}(-N_{c+d,a+b+c+d+e} tu)\\[3mm] 
                                &=& x_{-a-b-2c-2d-e}( -\alpha_{56}\beta_{15}\gamma_3  tu).
\end{array} \eqno{(040)}
$$

$$
\begin{array}{rcl} 
[x_{-c-d}(u),x_{-e-f}(t)] &=& x_{-c-d-e-f}(-N_{c+d,e+f} tu)\\[3mm] 
                          &=& x_{-c-d-e-f}( -\beta_{16}\gamma_4 tu).
\end{array} \eqno{(041)}
$$

$$
\begin{array}{rcl} 
[x_{-c-d}(u),x_{-b-d-e-f}(t)] &=& x_{-b-c-2d-e-f}(-N_{c+d,b+d+e+f} tu)\\[3mm] 
                              &=& x_{-b-c-2d-e-f}( \beta_{17}\gamma_5  tu).
\end{array} \eqno{(042)}
$$

$$
\begin{array}{rcl} 
[x_{-c-d}(u),x_{-a-b-c-d-e-f}(t)] &=& x_{-a-b-2c-2d-e-f}(-N_{c+d,a+b+c+d+e+f} tu)\\[3mm]
                                 &=& x_{-a-b-2c-2d-e-f}( -\alpha_{89}\beta_{18}\gamma_5  tu).
\end{array} \eqno{(043)}
$$

$$
\begin{array}{rcl} 
[x_{-c-d}(u),x_{-a-b-c-2d-2e-f}(t)] &=& x_{-a-b-2c-3d-2e-f}(-N_{c+d,a+b+c+2d+2e+f} tu)\\[3mm] 
                                    &=& x_{-a-b-2c-3d-2e-f}( \beta_{19}\gamma_6 tu).
\end{array}%
 \eqno{(120)}
$$

\vskip3mm
\centerline{\textit{Positive and negative roots}}

{\boldmath$a:$}

$$
\begin{array}{rcl}  
[x_{a}(u),x_{-a-c}(t)] &=& x_{-c}(N_{a,-a-c} tu) \\[3mm] 
                       &=& x_{-c}( -\alpha_1 tu).
\end{array}%
 \eqno{(240)}
$$

\newpage

\section{Special Case}

\subsection{The Table of the Structure Constant} It is true the following
\vskip4pt

\textbf{Theorem 2.} \textit{Let the values of all structure constants corresponding to extraspecial pairs be positive.
Then the values of arbitrary structure constants are equal to the values indicated in the following tables 5 and 6. }
\vskip5mm

\noindent
%
\vskip10pt
\centerline{Table 6, part 2 of 2.}
		
\subsection{List of Formulas}

\begin{center}
\textit{Positive roots}
\end{center}

{\boldmath$a:$}

$$
[x_a(u),x_{c}(t)]=x_{a+c}(tu).\eqno{(001)}
$$

$$
[x_a(u),x_{c+d}(t)]=x_{a+c+d}(tu).\eqno{(002)}
$$

$$
[x_a(u),x_{b+c+d}(t)]=x_{a+b+c+d}(tu).\eqno{(003)}
$$

$$
[x_a(u),x_{c+d+e}(t)]=x_{a+c+d+e}(tu).\eqno{(004)}
$$

$$
[x_a(u),x_{b+c+d+e}(t)]=x_{a+b+c+d+e}(tu).\eqno{(005)}
$$

$$
[x_a(u),x_{b+c+2d+e}(t)]=x_{a+b+c+2d+e}(tu).\eqno{(006)}
$$

$$
[x_a(u),x_{c+d+e+f}(t)]=x_{a+c+d+e+f}(tu).\eqno{(007)}
$$

$$
[x_a(u),x_{b+c+d+e+f}(t)]=x_{a+b+c+d+e+f}(tu).\eqno{(008)}
$$

$$
[x_a(u),x_{b+c+2d+e+f}(t)]=x_{a+b+c+2d+e+f}(tu).\eqno{(009)}
$$

$$
[x_a(u),x_{b+c+2d+2e+f}(t)]=x_{a+b+c+2d+2e+f}(tu).\eqno{(010)}
$$

{\boldmath$c:$}

$$
[x_c(u),x_{d}(t)]=x_{c+d}(tu).\eqno{(011)}
$$

$$
[x_c(u),x_{b+d}(t)]=x_{b+c+d}(tu).\eqno{(012)}
$$

$$
[x_c(u),x_{d+e}(t)]=x_{c+d+e}(tu).\eqno{(013)}
$$

$$
[x_c(u),x_{b+d+e}(t)]=x_{b+c+d+e}(tu).\eqno{(014)}
$$

$$
[x_c(u),x_{a+b+c+2d+e}(t)]=x_{a+b+2c+2d+e}(tu).\eqno{(015)}
$$

$$
[x_c(u),x_{d+e+f}(t)]=x_{c+d+e+f}(tu).\eqno{(016)}
$$

$$
[x_c(u),x_{b+d+e+f}(t)]=x_{b+c+d+e+f}(tu).\eqno{(017)}
$$

$$
[x_c(u),x_{a+b+c+2d+e+f}(t)]=x_{a+b+2c+2d+e+f}(tu).\eqno{(018)}
$$

$$
[x_c(u),x_{a+b+c+2d+2e+f}(t)]=x_{a+b+2c+2d+2e+f}(tu).\eqno{(019)}
$$

{\boldmath$a+c:$}

$$
[x_{a+c}(u),x_{d}(t)]=x_{a+c+d}(tu).\eqno{(020)}
$$

$$
[x_{a+c}(u),x_{b+d}(t)]=x_{a+b+c+d}(tu).\eqno{(021)}
$$

$$
[x_{a+c}(u),x_{d+e}(t)]=x_{a+c+d+e}(tu).\eqno{(022)}
$$

$$
[x_{a+c}(u),x_{b+d+e}(t)]=x_{a+b+c+d+e}(tu).\eqno{(023)}
$$

$$
[x_{a+c}(u),x_{b+c+2d+e}(t)]=x_{a+b+2c+2d+e}( -tu).\eqno{(024)}
$$

$$
[x_{a+c}(u),x_{d+e+f}(t)]=x_{a+c+d+e+f}(tu).\eqno{(025)}
$$

$$
[x_{a+c}(u),x_{b+d+e+f}(t)]=x_{a+b+c+d+e+f}(tu).\eqno{(026)}
$$

$$
[x_{a+c}(u),x_{b+c+2d+e+f}(t)]=x_{a+b+2c+2d+e+f}( -tu).\eqno{(027)}
$$

$$
[x_{a+c}(u),x_{b+c+2d+2e+f}(t)]=x_{a+b+2c+2d+2e+f}( -tu).\eqno{(028)}
$$

{\boldmath$d:$}

$$
[x_d(u),x_{b}(t)]=x_{b+d}(tu).\eqno{(029)}
$$

$$
[x_d(u),x_{e}(t)]=x_{d+e}(tu).\eqno{(030)}
$$

$$
[x_d(u),x_{b+c+d+e}(t)]=x_{b+c+2d+e}(tu).\eqno{(031)}
$$

$$
[x_d(u),x_{a+b+c+d+e}(t)]=x_{a+b+c+2d+e}(tu).\eqno{(032)}
$$

$$
[x_d(u),x_{e+f}(t)]=x_{d+e+f}(tu).\eqno{(033)}
$$

$$
[x_d(u),x_{b+c+d+e+f}(t)]=x_{b+c+2d+e+f}(tu).\eqno{(034)}
$$

$$
[x_d(u),x_{a+b+c+d+e+f}(t)]=x_{a+b+c+2d+e+f}(tu).\eqno{(035)}
$$

$$
[x_d(u),x_{a+b+2c+2d+2e+f}(t)]=x_{a+b+2c+3d+2e+f}(tu).\eqno{(036)}
$$

{\boldmath$c+d:$}

$$
[x_{c+d}(u),x_{b}(t)]=x_{b+c+d}(tu).\eqno{(037)}
$$

$$
[x_{c+d}(u),x_{e}(t)]=x_{c+d+e}(tu).\eqno{(038)}
$$

$$
[x_{c+d}(u),x_{b+d+e}(t)]=x_{b+c+2d+e}( -tu).\eqno{(039)}
$$

$$
[x_{c+d}(u),x_{a+b+c+d+e}(t)]=x_{a+b+2c+2d+e}(tu).\eqno{(040)}
$$

$$
[x_{c+d}(u),x_{e+f}(t)]=x_{c+d+e+f}(tu).\eqno{(041)}
$$

$$
[x_{c+d}(u),x_{b+d+e+f}(t)]=x_{b+c+2d+e+f}( -tu).\eqno{(042)}
$$

$$
[x_{c+d}(u),x_{a+b+c+d+e+f}(t)]=x_{a+b+2c+2d+e+f}(tu).\eqno{(043)}
$$

$$
[x_{c+d}(u),x_{a+b+c+2d+2e+f}(t)]=x_{a+b+2c+3d+2e+f}( -tu).\eqno{(044)}
$$

{\boldmath$a+c+d:$}

$$
[x_{a+c+d}(u),x_{b}(t)]=x_{a+b+c+d}(tu).\eqno{(045)}
$$

$$
[x_{a+c+d}(u),x_{e}(t)]=x_{a+c+d+e}(tu).\eqno{(046)}
$$

$$
[x_{a+c+d}(u),x_{b+d+e}(t)]=x_{a+b+c+2d+e}( -tu).\eqno{(047)}
$$

$$
[x_{a+c+d}(u),x_{b+c+d+e}(t)]=x_{a+b+2c+2d+e}( -tu).\eqno{(048)}
$$

$$
[x_{a+c+d}(u),x_{e+f}(t)]=x_{a+c+d+e+f}(tu).\eqno{(049)}
$$

$$
[x_{a+c+d}(u),x_{b+d+e+f}(t)]=x_{a+b+c+2d+e+f}( -tu).\eqno{(050)}
$$

$$
[x_{a+c+d}(u),x_{b+c+d+e+f}(t)]=x_{a+b+2c+2d+e+f}( -tu).\eqno{(051)}
$$

$$
[x_{a+c+d}(u),x_{b+c+2d+2e+f}(t)]=x_{a+b+2c+3d+2e+f}(tu).\eqno{(052)}
$$

{\boldmath$b:$}

$$
[x_b(u),x_{d+e}(t)]=x_{b+d+e}(tu).\eqno{(053)}
$$

$$
[x_b(u),x_{c+d+e}(t)]=x_{b+c+d+e}(tu).\eqno{(054)}
$$

$$
[x_b(u),x_{a+c+d+e}(t)]=x_{a+b+c+d+e}(tu).\eqno{(055)}
$$

$$
[x_b(u),x_{d+e+f}(t)]=x_{b+d+e+f}(tu).\eqno{(056)}
$$

$$
[x_b(u),x_{c+d+e+f}(t)]=x_{b+c+d+e+f}(tu).\eqno{(057)}
$$

$$
[x_b(u),x_{a+c+d+e+f}(t)]=x_{a+b+c+d+e+f}(tu).\eqno{(058)}
$$

$$
[x_b(u),x_{a+b+2c+3d+2e+f}(t)]=x_{a+2b+2c+3d+2e+f}(tu).\eqno{(059)}
$$

{\boldmath$b+d:$}

$$
[x_{b+d}(u),x_{e}(t)]=x_{b+d+e}( -tu).\eqno{(060)}
$$

$$
[x_{b+d}(u),x_{c+d+e}(t)]=x_{b+c+2d+e}(tu).\eqno{(061)}
$$

$$
[x_{b+d}(u),x_{a+c+d+e}(t)]=x_{a+b+c+2d+e}(tu).\eqno{(062)}
$$

$$
[x_{b+d}(u),x_{e+f}(t)]=x_{b+d+e+f}( -tu).\eqno{(063)}
$$

$$
[x_{b+d}(u),x_{c+d+e+f}(t)]=x_{b+c+2d+e+f}(tu).\eqno{(064)}
$$

$$
[x_{b+d}(u),x_{a+c+d+e+f}(t)]=x_{a+b+c+2d+e+f}(tu).\eqno{(065)}
$$

$$
[x_{b+d}(u),x_{a+b+2c+2d+2e+f}(t)]=x_{a+2b+2c+3d+2e+f}( -tu).\eqno{(066)}
$$

{\boldmath$b+c+d:$}

$$
[x_{b+c+d}(u),x_{e}(t)]=x_{b+c+d+e}( -tu).\eqno{(067)}
$$

$$
[x_{b+c+d}(u),x_{d+e}(t)]=x_{b+c+2d+e}( -tu).\eqno{(068)}
$$

$$
[x_{b+c+d}(u),x_{a+c+d+e}(t)]=x_{a+b+2c+2d+e}(tu).\eqno{(069)}
$$

$$
[x_{b+c+d}(u),x_{e+f}(t)]=x_{b+c+d+e+f}( -tu).\eqno{(070)}
$$

$$
[x_{b+c+d}(u),x_{d+e+f}(t)]=x_{b+c+2d+e+f}( -tu).\eqno{(071)}
$$

$$
[x_{b+c+d}(u),x_{a+c+d+e+f}(t)]=x_{a+b+2c+2d+e+f}(tu).\eqno{(072)}
$$

$$
[x_{b+c+d}(u),x_{a+b+c+2d+2e+f}(t)]=x_{a+2b+2c+3d+2e+f}(tu).\eqno{(073)}
$$

{\boldmath$a+b+c+d:$}

$$
[x_{a+b+c+d}(u),x_{e}(t)]=x_{a+b+c+d+e}( -tu).\eqno{(074)}
$$

$$
[x_{a+b+c+d}(u),x_{d+e}(t)]=x_{a+b+c+2d+e}( -tu).\eqno{(075)}
$$

$$
[x_{a+b+c+d}(u),x_{c+d+e}(t)]=x_{a+b+2c+2d+e}( -tu).\eqno{(076)}
$$

$$
[x_{a+b+c+d}(u),x_{e+f}(t)]=x_{a+b+c+d+e+f}( -tu).\eqno{(077)}
$$

$$
[x_{a+b+c+d}(u),x_{d+e+f}(t)]=x_{a+b+c+2d+e+f}( -tu).\eqno{(078)}
$$

$$
[x_{a+b+c+d}(u),x_{c+d+e+f}(t)]=x_{a+b+2c+2d+e+f}( -tu).\eqno{(079)}
$$

$$
[x_{a+b+c+d}(u),x_{b+c+2d+2e+f}(t)]=x_{a+2b+2c+3d+2e+f}( -tu).\eqno{(080)}
$$

{\boldmath$e:$}

$$
[x_e(u),x_{f}(t)]=x_{e+f}(tu).\eqno{(081)}
$$

$$
[x_e(u),x_{b+c+2d+e+f}(t)]=x_{b+c+2d+2e+f}(tu).\eqno{(082)}
$$

$$
[x_e(u),x_{a+b+c+2d+e+f}(t)]=x_{a+b+c+2d+2e+f}(tu).\eqno{(083)}
$$

$$
[x_e(u),x_{a+b+2c+2d+e+f}(t)]=x_{a+b+2c+2d+2e+f}(tu).\eqno{(084)}
$$

{\boldmath$d+e:$}

$$
[x_{d+e}(u),x_{f}(t)]=x_{d+e+f}(tu).\eqno{(085)}
$$

$$
[x_{d+e}(u),x_{b+c+d+e+f}(t)]=x_{b+c+2d+2e+f}( -tu).\eqno{(086)}
$$

$$
[x_{d+e}(u),x_{a+b+c+d+e+f}(t)]=x_{a+b+c+2d+2e+f}( -tu).\eqno{(087)}
$$

$$
[x_{d+e}(u),x_{a+b+2c+2d+e+f}(t)]=x_{a+b+2c+3d+2e+f}(tu).\eqno{(088)}
$$

{\boldmath$c+d+e:$}

$$
[x_{c+d+e}(u),x_{f}(t)]=x_{c+d+e+f}(tu).\eqno{(089)}
$$

$$
[x_{c+d+e}(u),x_{b+d+e+f}(t)]=x_{b+c+2d+2e+f}(tu).\eqno{(090)}
$$

$$
[x_{c+d+e}(u),x_{a+b+c+d+e+f}(t)]=x_{a+b+2c+2d+2e+f}( -tu).\eqno{(091)}
$$

$$
[x_{c+d+e}(u),x_{a+b+c+2d+e+f}(t)]=x_{a+b+2c+3d+2e+f}( -tu).\eqno{(092)}
$$

{\boldmath$a+c+d+e:$}

$$
[x_{a+c+d+e}(u),x_{f}(t)]=x_{a+c+d+e+f}(tu).\eqno{(093)}
$$

$$
[x_{a+c+d+e}(u),x_{b+d+e+f}(t)]=x_{a+b+c+2d+2e+f}(tu).\eqno{(094)}
$$

$$
[x_{a+c+d+e}(u),x_{b+c+d+e+f}(t)]=x_{a+b+2c+2d+2e+f}(tu).\eqno{(095)}
$$

$$
[x_{a+c+d+e}(u),x_{b+c+2d+e+f}(t)]=x_{a+b+2c+3d+2e+f}(tu).\eqno{(096)}
$$

{\boldmath$b+d+e:$}

$$
[x_{b+d+e}(u),x_{f}(t)]=x_{b+d+e+f}(tu).\eqno{(097)}
$$

$$
[x_{b+d+e}(u),x_{c+d+e+f}(t)]=x_{b+c+2d+2e+f}(tu).\eqno{(098)}
$$

$$
[x_{b+d+e}(u),x_{a+c+d+e+f}(t)]=x_{a+b+c+2d+2e+f}(tu).\eqno{(099)}
$$

$$
[x_{b+d+e}(u),x_{a+b+2c+2d+e+f}(t)]=x_{a+2b+2c+3d+2e+f}(tu).\eqno{(100)}
$$

{\boldmath$b+c+d+e:$}

$$
[x_{b+c+d+e}(u),x_{f}(t)]=x_{b+c+d+e+f}(tu).\eqno{(101)}
$$

$$
[x_{b+c+d+e}(u),x_{d+e+f}(t)]=x_{b+c+2d+2e+f}( -tu).\eqno{(102)}
$$

$$
[x_{b+c+d+e}(u),x_{a+c+d+e+f}(t)]=x_{a+b+2c+2d+2e+f}(tu).\eqno{(103)}
$$

$$
[x_{b+c+d+e}(u),x_{a+b+c+2d+e+f}(t)]=x_{a+2b+2c+3d+2e+f}( -tu).\eqno{(104)}
$$

{\boldmath$a+b+c+d+e:$}

$$
[x_{a+b+c+d+e}(u),x_{f}(t)]=x_{a+b+c+d+e+f}(tu).\eqno{(105)}
$$

$$
[x_{a+b+c+d+e}(u),x_{d+e+f}(t)]=x_{a+b+c+2d+2e+f}( -tu).\eqno{(106)}
$$

$$
[x_{a+b+c+d+e}(u),x_{c+d+e+f}(t)]=x_{a+b+2c+2d+2e+f}( -tu).\eqno{(107)}
$$

$$
[x_{a+b+c+d+e}(u),x_{b+c+2d+e+f}(t)]=x_{a+2b+2c+3d+2e+f}(tu).\eqno{(108)}
$$

{\boldmath$b+c+2d+e:$}

$$
[x_{b+c+2d+e}(u),x_{f}(t)]=x_{b+c+2d+e+f}(tu).\eqno{(109)}
$$

$$
[x_{b+c+2d+e}(u),x_{e+f}(t)]=x_{b+c+2d+2e+f}(tu).\eqno{(110)}
$$

$$
[x_{b+c+2d+e}(u),x_{a+c+d+e+f}(t)]=x_{a+b+2c+3d+2e+f}(tu).\eqno{(111)}
$$

$$
[x_{b+c+2d+e}(u),x_{a+b+c+d+e+f}(t)]=x_{a+2b+2c+3d+2e+f}(tu).\eqno{(112)}
$$

{\boldmath$a+b+c+2d+e:$}

$$
[x_{a+b+c+2d+e}(u),x_{f}(t)]=x_{a+b+c+2d+e+f}(tu).\eqno{(113)}
$$

$$
[x_{a+b+c+2d+e}(u),x_{e+f}(t)]=x_{a+b+c+2d+2e+f}(tu).\eqno{(114)}
$$

$$
[x_{a+b+c+2d+e}(u),x_{c+d+e+f}(t)]=x_{a+b+2c+3d+2e+f}( -tu).\eqno{(115)}
$$

$$
[x_{a+b+c+2d+e}(u),x_{b+c+d+e+f}(t)]=x_{a+2b+2c+3d+2e+f}( -tu).\eqno{(116)}
$$
{\boldmath$a+b+2c+2d+e:$}

$$
[x_{a+b+2c+2d+e}(u),x_{f}(t)]=x_{a+b+2c+2d+e+f}(tu).\eqno{(117)}
$$

$$
[x_{a+b+2c+2d+e}(u),x_{e+f}(t)]=x_{a+b+2c+2d+2e+f}(tu).\eqno{(118)}
$$

$$
[x_{a+b+2c+2d+e}(u),x_{d+e+f}(t)]=x_{a+b+2c+3d+2e+f}(tu).\eqno{(119)}
$$

$$
[x_{a+b+2c+2d+e}(u),x_{b+d+e+f}(t)]=x_{a+2b+2c+3d+2e+f}(tu).\eqno{(120)}
$$

\vskip5mm
\begin{center}
\textit{Negative roots}
\end{center}

{\boldmath$-a:$}

$$
[x_{-a}(u),x_{-c}(t)]=x_{-a-c}( -tu).\eqno{(001)}
$$

$$
[x_{-a}(u),x_{-c-d}(t)]=x_{-a-c-d}( -tu).\eqno{(002)}
$$

$$
[x_{-a}(u),x_{-b-c-d}(t)]=x_{-a-b-c-d}( -tu).\eqno{(003)}
$$

$$
[x_{-a}(u),x_{-c-d-e}(t)]=x_{-a-c-d-e}( -tu).\eqno{(004)}
$$

$$
[x_{-a}(u),x_{-b-c-d-e}(t)]=x_{-a-b-c-d-e}( -tu).\eqno{(005)}
$$

$$
[x_{-a}(u),x_{-b-c-2d-e}(t)]=x_{-a-b-c-2d-e}( -tu).\eqno{(006)}
$$

$$
[x_{-a}(u),x_{-c-d-e-f}(t)]=x_{-a-c-d-e-f}( -tu).\eqno{(007)}
$$

$$
[x_{-a}(u),x_{-b-c-d-e-f}(t)]=x_{-a-b-c-d-e-f}( -tu).\eqno{(008)}
$$

$$
[x_{-a}(u),x_{-b-c-2d-e-f}(t)]=x_{-a-b-c-2d-e-f}( -tu).\eqno{(009)}
$$

$$
[x_{-a}(u),x_{-b-c-2d-2e-f}(t)]=x_{-a-b-c-2d-2e-f}( -tu).\eqno{(010)}
$$

{\boldmath$-c:$}

$$
[x_{-c}(u),x_{-d}(t)]=x_{-c-d}( -tu).\eqno{(011)}
$$

$$
[x_{-c}(u),x_{-b-d}(t)]=x_{-b-c-d}( -tu).\eqno{(012)}
$$

$$
[x_{-c}(u),x_{-d-e}(t)]=x_{-c-d-e}( -tu).\eqno{(013)}
$$

$$
[x_{-c}(u),x_{-b-d-e}(t)]=x_{-b-c-d-e}( -tu).\eqno{(014)}
$$

$$
[x_{-c}(u),x_{-a-b-c-2d-e}(t)]=x_{-a-b-2c-2d-e}( -tu).\eqno{(015)}
$$

$$
[x_{-c}(u),x_{-d-e-f}(t)]=x_{-c-d-e-f}( -tu).\eqno{(016)}
$$

$$
[x_{-c}(u),x_{-b-d-e-f}(t)]=x_{-b-c-d-e-f}( -tu).\eqno{(017)}
$$

$$
[x_{-c}(u),x_{-a-b-c-2d-e-f}(t)]=x_{-a-b-2c-2d-e-f}( -tu).\eqno{(018)}
$$

$$
[x_{-c}(u),x_{-a-b-c-2d-2e-f}(t)]=x_{-a-b-2c-2d-2e-f}( -tu).\eqno{(019)}
$$

{\boldmath$-a-c:$}

$$
[x_{-a-c}(u),x_{-d}(t)]=x_{-a-c-d}( -tu).\eqno{(020)}
$$

$$
[x_{-a-c}(u),x_{-b-d}(t)]=x_{-a-b-c-d}( -tu).\eqno{(021)}
$$

$$
[x_{-a-c}(u),x_{-d-e}(t)]=x_{-a-c-d-e}( -tu).\eqno{(022)}
$$

$$
[x_{-a-c}(u),x_{-b-d-e}(t)]=x_{-a-b-c-d-e}( -tu).\eqno{(023)}
$$

$$
[x_{-a-c}(u),x_{-b-c-2d-e}(t)]=x_{-a-b-2c-2d-e}( tu).\eqno{(024)}
$$

$$
[x_{-a-c}(u),x_{-d-e-f}(t)]=x_{-a-c-d-e-f}( -tu).\eqno{(025)}
$$

$$
[x_{-a-c}(u),x_{-b-d-e-f}(t)]=x_{-a-b-c-d-e-f}( -tu).\eqno{(026)}
$$

$$
[x_{-a-c}(u),x_{-b-c-2d-e-f}(t)]=x_{-a-b-2c-2d-e-f}( tu).\eqno{(027)}
$$

$$
[x_{-a-c}(u),x_{-b-c-2d-2e-f}(t)]=x_{-a-b-2c-2d-2e-f}( tu).\eqno{(028)}
$$

{\boldmath$-d:$}

$$
[x_{-d}(u),x_{-b}(t)]=x_{-b-d}( -tu).\eqno{(029)}
$$

$$
[x_{-d}(u),x_{-e}(t)]=x_{-d-e}( -tu).\eqno{(030)}
$$

$$
[x_{-d}(u),x_{-b-c-d-e}(t)]=x_{-b-c-2d-e}( -tu).\eqno{(031)}
$$

$$
[x_{-d}(u),x_{-a-b-c-d-e}(t)]=x_{-a-b-c-2d-e}(  -tu).\eqno{(032)}
$$

$$
[x_{-d}(u),x_{-e-f}(t)]=x_{-d-e-f}( -tu).\eqno{(033)}
$$

$$
[x_{-d}(u),x_{-b-c-d-e-f}(t)]=x_{-b-c-2d-e-f}( -tu).\eqno{(034)}
$$

$$
[x_{-d}(u),x_{-a-b-c-d-e-f}(t)]=x_{-a-b-c-2d-e-f}( -tu).\eqno{(035)}
$$

$$
[x_{-d}(u),x_{-a-b-2c-2d-2e-f}(t)]=x_{-a-b-2c-3d-2e-f}( -tu).\eqno{(036)}
$$

{\boldmath$-c-d:$}

$$
[x_{-c-d}(u),x_{-b}(t)]=x_{-b-c-d}( -tu).\eqno{(037)}
$$

$$
[x_{-c-d}(u),x_{-e}(t)]=x_{-c-d-e}( -tu).\eqno{(038)}
$$

$$
[x_{-c-d}(u),x_{-b-d-e}(t)]=x_{-b-c-2d-e}(tu).\eqno{(039)}
$$

$$
[x_{-c-d}(u),x_{-a-b-c-d-e}(t)]=x_{-a-b-2c-2d-e}( -tu).\eqno{(040)}
$$

$$
[x_{-c-d}(u),x_{-e-f}(t)]=x_{-c-d-e-f}( -tu).\eqno{(041)}
$$

$$
[x_{-c-d}(u),x_{-b-d-e-f}(t)]=x_{-b-c-2d-e-f}(tu).\eqno{(042)}
$$

$$
[x_{-c-d}(u),x_{-a-b-c-d-e-f}(t)]=x_{-a-b-2c-2d-e-f}( -tu).\eqno{(043)}
$$

$$
[x_{-c-d}(u),x_{-a-b-c-2d-2e-f}(t)]=x_{-a-b-2c-3d-2e-f}(tu).\eqno{(044)}
$$

{\boldmath$-a-c-d:$}

$$
[x_{-a-c-d}(u),x_{-b}(t)]=x_{-a-b-c-d}( -tu).\eqno{(045)}
$$

$$
[x_{-a-c-d}(u),x_{-e}(t)]=x_{-a-c-d-e}( -tu).\eqno{(046)}
$$

$$
[x_{-a-c-d}(u),x_{-b-d-e}(t)]=x_{-a-b-c-2d-e}(tu).\eqno{(047)}
$$

$$
[x_{-a-c-d}(u),x_{-b-c-d-e}(t)]=x_{-a-b-2c-2d-e}(tu).\eqno{(048)}
$$

$$
[x_{-a-c-d}(u),x_{-e-f}(t)]=x_{-a-c-d-e-f}( -tu).\eqno{(049)}
$$

$$
[x_{-a-c-d}(u),x_{-b-d-e-f}(t)]=x_{-a-b-c-2d-e-f}(tu).\eqno{(050)}
$$

$$
[x_{-a-c-d}(u),x_{-b-c-d-e-f}(t)]=x_{-a-b-2c-2d-e-f}(tu).\eqno{(051)}
$$

$$
[x_{-a-c-d}(u),x_{-b-c-2d-2e-f}(t)]=x_{-a-b-2c-3d-2e-f}( -tu).\eqno{(052)}
$$

{\boldmath$-b:$}

$$
[x_{-b}(u),x_{-d-e}(t)]=x_{-b-d-e}( -tu).\eqno{(053)}
$$

$$
[x_{-b}(u),x_{-c-d-e}(t)]=x_{-b-c-d-e}( -tu).\eqno{(054)}
$$

$$
[x_{-b}(u),x_{-a-c-d-e}(t)]=x_{-a-b-c-d-e}( -tu).\eqno{(055)}
$$

$$
[x_{-b}(u),x_{-d-e-f}(t)]=x_{-b-d-e-f}( -tu).\eqno{(056)}
$$

$$
[x_{-b}(u),x_{-c-d-e-f}(t)]=x_{-b-c-d-e-f}( -tu).\eqno{(057)}
$$

$$
[x_{-b}(u),x_{-a-c-d-e-f}(t)]=x_{-a-b-c-d-e-f}( -tu).\eqno{(058)}
$$

$$
[x_{-b}(u),x_{-a-b-2c-3d-2e-f}(t)]=x_{-a-2b-2c-3d-2e-f}( -tu).\eqno{(059)}
$$

{\boldmath$-b-d:$}

$$
[x_{-b-d}(u),x_{-e}(t)]=x_{-b-d-e}(tu).\eqno{(060)}
$$

$$
[x_{-b-d}(u),x_{-c-d-e}(t)]=x_{-b-c-2d-e}( -tu).\eqno{(061)}
$$

$$
[x_{-b-d}(u),x_{-a-c-d-e}(t)]=x_{-a-b-c-2d-e}( -tu).\eqno{(062)}
$$

$$
[x_{-b-d}(u),x_{-e-f}(t)]=x_{-b-d-e-f}(tu).\eqno{(063)}
$$

$$
[x_{-b-d}(u),x_{-c-d-e-f}(t)]=x_{-b-c-2d-e-f}( -tu).\eqno{(064)}
$$

$$
[x_{-b-d}(u),x_{-a-c-d-e-f}(t)]=x_{-a-b-c-2d-e-f}( -tu).\eqno{(065)}
$$

$$
[x_{-b-d}(u),x_{-a-b-2c-2d-2e-f}(t)]=x_{-a-2b-2c-3d-2e-f}(tu).\eqno{(066)}
$$

{\boldmath$-b-c-d:$}

$$
[x_{-b-c-d}(u),x_{-e}(t)]=x_{-b-c-d-e}(tu).\eqno{(067)}
$$

$$
[x_{-b-c-d}(u),x_{-d-e}(t)]=x_{-b-c-2d-e}(tu).\eqno{(068)}
$$

$$
[x_{-b-c-d}(u),x_{-a-c-d-e}(t)]=x_{-a-b-2c-2d-e}( -tu).\eqno{(069)}
$$

$$
[x_{-b-c-d}(u),x_{-e-f}(t)]=x_{-b-c-d-e-f}(tu).\eqno{(070)}
$$

$$
[x_{-b-c-d}(u),x_{-d-e-f}(t)]=x_{-b-c-2d-e-f}(tu).\eqno{(071)}
$$

$$
[x_{-b-c-d}(u),x_{-a-c-d-e-f}(t)]=x_{-a-b-2c-2d-e-f}( -tu).\eqno{(072)}
$$

$$
[x_{-b-c-d}(u),x_{-a-b-c-2d-2e-f}(t)]=x_{-a-2b-2c-3d-2e-f}( -tu).\eqno{(073)}
$$

{\boldmath$-a-b-c-d:$}

$$
[x_{-a-b-c-d}(u),x_{-e}(t)]=x_{-a-b-c-d-e}(tu).\eqno{(074)}
$$

$$
[x_{-a-b-c-d}(u),x_{-d-e}(t)]=x_{-a-b-c-2d-e}(tu).\eqno{(075)}
$$

$$
[x_{-a-b-c-d}(u),x_{-c-d-e}(t)]=x_{-a-b-2c-2d-e}(tu).\eqno{(076)}
$$

$$
[x_{-a-b-c-d}(u),x_{-e-f}(t)]=x_{-a-b-c-d-e-f}(tu).\eqno{(077)}
$$

$$
[x_{-a-b-c-d}(u),x_{-d-e-f}(t)]=x_{-a-b-c-2d-e-f}(tu).\eqno{(078)}
$$

$$
[x_{-a-b-c-d}(u),x_{-c-d-e-f}(t)]=x_{-a-b-2c-2d-e-f}(tu).\eqno{(079)}
$$

$$
[x_{-a-b-c-d}(u),x_{-b-c-2d-2e-f}(t)]=x_{-a-2b-2c-3d-2e-f}(tu).\eqno{(080)}
$$

{\boldmath$-e:$}

$$
[x_{-e}(u),x_{-f}(t)]=x_{-e-f}( - tu).\eqno{(081)}
$$

$$
[x_{-e}(u),x_{-b-c-2d-e-f}(t)]=x_{-b-c-2d-2e-f}( -tu).\eqno{(082)}
$$

$$
[x_{-e}(u),x_{-a-b-c-2d-e-f}(t)]=x_{-a-b-c-2d-2e-f}( -tu).\eqno{(083)}
$$

$$
[x_{-e}(u),x_{-a-b-2c-2d-e-f}(t)]=x_{-a-b-2c-2d-2e-f}( -tu).\eqno{(084)}
$$

{\boldmath$-d-e:$}

$$
[x_{-d-e}(u),x_{-f}(t)]=x_{-d-e-f}( -tu).\eqno{(085)}
$$

$$
[x_{-d-e}(u),x_{-b-c-d-e-f}(t)]=x_{-b-c-2d-2e-f}(tu).\eqno{(086)}
$$

$$
[x_{-d-e}(u),x_{-a-b-c-d-e-f}(t)]=x_{-a-b-c-2d-2e-f}(tu).\eqno{(087)}
$$

$$
[x_{-d-e}(u),x_{-a-b-2c-2d-e-f}(t)]=x_{-a-b-2c-3d-2e-f}( -tu).\eqno{(088)}
$$

{\boldmath$-c-d-e:$}

$$
[x_{-c-d-e}(u),x_{-f}(t)]=x_{-c-d-e-f}( -tu).\eqno{(089)}
$$

$$
[x_{-c-d-e}(u),x_{-b-d-e-f}(t)]=x_{-b-c-2d-2e-f}( -tu).\eqno{(090)}
$$

$$
[x_{-c-d-e}(u),x_{-a-b-c-d-e-f}(t)]=x_{-a-b-2c-2d-2e-f}(tu).\eqno{(091)}
$$

$$
[x_{-c-d-e}(u),x_{-a-b-c-2d-e-f}(t)]=x_{-a-b-2c-3d-2e-f}(tu).\eqno{(092)}
$$

{\boldmath$-a-c-d-e:$}

$$
[x_{-a-c-d-e}(u),x_{-f}(t)]=x_{-a-c-d-e-f}( -tu).\eqno{(093)}
$$

$$
[x_{-a-c-d-e}(u),x_{-b-d-e-f}(t)]=x_{-a-b-c-2d-2e-f}( -tu).\eqno{(094)}
$$

$$
[x_{-a-c-d-e}(u),x_{-b-c-d-e-f}(t)]=x_{-a-b-2c-2d-2e-f}( -tu).\eqno{(095)}
$$

$$
[x_{-a-c-d-e}(u),x_{-b-c-2d-e-f}(t)]=x_{-a-b-2c-3d-2e-f}( -tu).\eqno{(096)}
$$

{\boldmath$-b-d-e:$}

$$
[x_{-b-d-e}(u),x_{-f}(t)]=x_{-b-d-e-f}( -tu).\eqno{(097)}
$$

$$
[x_{-b-d-e}(u),x_{-c-d-e-f}(t)]=x_{-b-c-2d-2e-f}( -tu).\eqno{(098)}
$$

$$
[x_{-b-d-e}(u),x_{-a-c-d-e-f}(t)]=x_{-a-b-c-2d-2e-f}( -tu).\eqno{(099)}
$$

$$
[x_{-b-d-e}(u),x_{-a-b-2c-2d-e-f}(t)]=x_{-a-2b-2c-3d-2e-f}( -tu).\eqno{(100)}
$$

{\boldmath$-b-c-d-e:$}

$$
[x_{-b-c-d-e}(u),x_{-f}(t)]=x_{-b-c-d-e-f}( -tu).\eqno{(101)}
$$

$$
[x_{-b-c-d-e}(u),x_{-d-e-f}(t)]=x_{-b-c-2d-2e-f}(tu).\eqno{(102)}
$$

$$
[x_{-b-c-d-e}(u),x_{-a-c-d-e-f}(t)]=x_{-a-b-2c-2d-2e-f}( -tu).\eqno{(103)}
$$

$$
[x_{-b-c-d-e}(u),x_{-a-b-c-2d-e-f}(t)]=x_{-a-2b-2c-3d-2e-f}(tu).\eqno{(104)}
$$

{\boldmath$-a-b-c-d-e:$}

$$
[x_{-a-b-c-d-e}(u),x_{-f}(t)]=x_{-a-b-c-d-e-f}( -tu).\eqno{(105)}
$$

$$
[x_{-a-b-c-d-e}(u),x_{-d-e-f}(t)]=x_{-a-b-c-2d-2e-f}(tu).\eqno{(106)}
$$

$$
[x_{-a-b-c-d-e}(u),x_{-c-d-e-f}(t)]=x_{-a-b-2c-2d-2e-f}(tu).\eqno{(107)}
$$

$$
[x_{-a-b-c-d-e}(u),x_{-b-c-2d-e-f}(t)]=x_{-a-2b-2c-3d-2e-f}( -tu).\eqno{(108)}
$$

{\boldmath$-b-c-2d-e:$}

$$
[x_{-b-c-2d-e}(u),x_{-f}(t)]=x_{-b-c-2d-e-f}( -tu).\eqno{(109)}
$$

$$
[x_{-b-c-2d-e}(u),x_{-e-f}(t)]=x_{-b-c-2d-2e-f}( -tu).\eqno{(110)}
$$

$$
[x_{-b-c-2d-e}(u),x_{-a-c-d-e-f}(t)]=x_{-a-b-2c-3d-2e-f}( -tu).\eqno{(111)}
$$

$$
[x_{-b-c-2d-e}(u),x_{-a-b-c-d-e-f}(t)]=x_{-a-2b-2c-3d-2e-f}( -tu).\eqno{(112)}
$$

{\boldmath$-a-b-c-2d-e:$}

$$
[x_{-a-b-c-2d-e}(u),x_{-f}(t)]=x_{-a-b-c-2d-e-f}( -tu).\eqno{(113)}
$$

$$
[x_{-a-b-c-2d-e}(u),x_{-e-f}(t)]=x_{-a-b-c-2d-2e-f}( -tu).\eqno{(114)}
$$

$$
[x_{-a-b-c-2d-e}(u),x_{-c-d-e-f}(t)]=x_{-a-b-2c-3d-2e-f}(tu).\eqno{(115)}
$$

$$
[x_{-a-b-c-2d-e}(u),x_{-b-c-d-e-f}(t)]=x_{-a-2b-2c-3d-2e-f}(tu).\eqno{(116)}
$$
{\boldmath$-a-b-2c-2d-e:$}

$$
[x_{-a-b-2c-2d-e}(u),x_{-f}(t)]=x_{-a-b-2c-2d-e-f}( -tu).\eqno{(117)}
$$

$$
[x_{-a-b-2c-2d-e}(u),x_{-e-f}(t)]=x_{-a-b-2c-2d-2e-f}( -tu).\eqno{(118)}
$$

$$
[x_{-a-b-2c-2d-e}(u),x_{-d-e-f}(t)]=x_{-a-b-2c-3d-2e-f}( -tu).\eqno{(119)}
$$

$$
[x_{-a-b-2c-2d-e}(u),x_{-b-d-e-f}(t)]=x_{-a-2b-2c-3d-2e-f}( -tu).\eqno{(120)}
$$

\vskip5mm
\begin{center}
\textit{Positive and negative roots}
\end{center}

{\boldmath$a:$}

$$
[x_{a}(u),x_{-a-c}(t)]=x_{-c}( -tu).\eqno{(001)}
$$

$$
[x_{a}(u),x_{-a-c-d}(t)]=x_{-c-d}( -tu).\eqno{(002)}
$$

$$
[x_{a}(u),x_{-a-b-c-d}(t)]=x_{-b-c-d}( -tu).\eqno{(003)}
$$

$$
[x_{a}(u),x_{-a-c-d-e}(t)]=x_{-c-d-e}( -tu).\eqno{(004)}
$$

$$
[x_{a}(u),x_{-a-b-c-d-e}(t)]=x_{-b-c-d-e}( -tu).\eqno{(005)}
$$

$$
[x_{a}(u),x_{-a-b-c-2d-e}(t)]=x_{-b-c-2d-e}( -tu).\eqno{(006)}
$$

$$
[x_{a}(u),x_{-a-c-d-e-f}(t)]=x_{-c-d-e-f}( -tu).\eqno{(007)}
$$

$$
[x_{a}(u),x_{-a-b-c-d-e-f}(t)]=x_{-b-c-d-e-f}( -tu).\eqno{(008)}
$$

$$
[x_{a}(u),x_{-a-b-c-2d-e-f}(t)]=x_{-b-c-2d-e-f}( -tu).\eqno{(009)}
$$

$$
[x_{a}(u),x_{-a-b-c-2d-2e-f}(t)]=x_{-b-c-2d-2e-f}( -tu).\eqno{(010)}
$$

{\boldmath$c:$}

$$
[x_{c}(u),x_{-a-c}(t)]=x_{-a}( tu).\eqno{(011)}
$$

$$
[x_{c}(u),x_{-c-d}(t)]=x_{-d}( -tu).\eqno{(012)}
$$

$$
[x_{c}(u),x_{-b-c-d}(t)]=x_{-b-d}( -tu).\eqno{(013)}
$$

$$
[x_{c}(u),x_{-c-d-e}(t)]=x_{-d-e}( -tu).\eqno{(014)}
$$

$$
[x_{c}(u),x_{-b-c-d-e}(t)]=x_{-b-d-e}( -tu).\eqno{(015)}
$$

$$
[x_{c}(u),x_{-a-b-2c-2d-e}(t)]=x_{-a-b-c-2d-e}( -tu).\eqno{(016)}
$$

$$
[x_{c}(u),x_{-c-d-e-f}(t)]=x_{-d-e-f}( -tu).\eqno{(017)}
$$

$$
[x_{c}(u),x_{-b-c-d-e-f}(t)]=x_{-b-d-e-f}( -tu).\eqno{(018)}
$$

$$
[x_{c}(u),x_{-a-b-2c-2d-e-f}(t)]=x_{-a-b-c-2d-e-f}( -tu).\eqno{(019)}
$$

$$
[x_{c}(u),x_{-a-b-2c-2d-2e-f}(t)]=x_{-a-b-c-2d-2e-f}( -tu).\eqno{(020)}
$$

{\boldmath$a+c:$}

$$
[x_{a+c}(u),x_{-a-c-d}(t)]=x_{-d}(-tu).\eqno{(021)}
$$

$$
[x_{a+c}(u),x_{-a-b-c-d}(t)]=x_{-b-d}(-tu).\eqno{(022)}
$$

$$
[x_{a+c}(u),x_{-a-c-d-e}(t)]=x_{-d-e}(-tu).\eqno{(023)}
$$

$$
[x_{a+c}(u),x_{-a-b-c-d-e}(t)]=x_{-b-d-e}(-tu).\eqno{(024)}
$$

$$
[x_{a+c}(u),x_{-a-b-2c-2d-e}(t)]=x_{-b-c-2d-e}(tu).\eqno{(025)}
$$

$$
[x_{a+c}(u),x_{-a-c-d-e-f}(t)]=x_{-d-e-f}(-tu).\eqno{(026)}
$$

$$
[x_{a+c}(u),x_{-a-b-c-d-e-f}(t)]=x_{-b-d-e-f}(-tu).\eqno{(027)}
$$

$$
[x_{a+c}(u),x_{-a-b-2c-2d-e-f}(t)]=x_{-b-c-2d-e-f}(tu).\eqno{(028)}
$$

$$
[x_{a+c}(u),x_{-a-b-2c-2d-2e-f}(t)]=x_{-b-c-2d-2e-f}(tu).\eqno{(029)}
$$

{\boldmath$d:$}

$$
[x_{d}(u),x_{-c-d}(t)]=x_{-c}(tu).\eqno{(030)}
$$

$$
[x_{d}(u),x_{-a-c-d}(t)]=x_{-a-c}(tu).\eqno{(031)}
$$

$$
[x_{d}(u),x_{-b-d}(t)]=x_{-b}( -tu).\eqno{(032)}
$$

$$
[x_{d}(u),x_{-d-e}(t)]=x_{-e}( -tu).\eqno{(033)}
$$

$$
[x_{d}(u),x_{-b-c-2d-e}(t)]=x_{-b-c-d-e}( -tu).\eqno{(034)}
$$

$$
[x_{d}(u),x_{-a-b-c-2d-e}(t)]=x_{-a-b-c-d-e}( -tu).\eqno{(035)}
$$

$$
[x_{d}(u),x_{-d-e-f}(t)]=x_{-e-f}( -tu).\eqno{(036)}
$$

$$
[x_{d}(u),x_{-b-c-2d-e-f}(t)]=x_{-b-c-d-e-f}( -tu).\eqno{(037)}
$$

$$
[x_{d}(u),x_{-a-b-c-2d-e-f}(t)]=x_{-a-b-c-d-e-f}( -tu).\eqno{(038)}
$$

$$
[x_{d}(u),x_{-a-b-2c-3d-2e-f}(t)]=x_{-a-b-2c-2d-2e-f}( -tu).\eqno{(039)}
$$

{\boldmath$c+d:$}

$$
[x_{c+d}(u),x_{-a-c-d}(t)]=x_{-a}(tu).\eqno{(040)}
$$

$$
[x_{c+d}(u),x_{-b-c-d}(t)]=x_{-b}(-tu).\eqno{(041)}
$$

$$
[x_{c+d}(u),x_{-c-d-e}(t)]=x_{-e}(-tu).\eqno{(042)}
$$

$$
[x_{c+d}(u),x_{-b-c-2d-e}(t)]=x_{-b-d-e}(tu).\eqno{(043)}
$$

$$
[x_{c+d}(u),x_{-a-b-2c-2d-e}(t)]=x_{-a-b-c-d-e}( -tu).\eqno{(044)}
$$

$$
[x_{c+d}(u),x_{-c-d-e-f}(t)]=x_{-e-f}(-tu).\eqno{(045)}
$$

$$
[x_{c+d}(u),x_{-b-c-2d-e-f}(t)] = x_{-b-d-e-f}( tu).\eqno{(046)}
$$

$$
[x_{c+d}(u),x_{-a-b-2c-2d-e-f}(t)]=x_{-a-b-c-d-e-f}(-tu).\eqno{(047)}
$$

$$
[x_{c+d}(u),x_{-a-b-2c-3d-2e-f}(t)]=x_{-a-b-c-2d-2e-f}(tu).\eqno{(048)}
$$

{\boldmath$a+c+d:$}

$$
[x_{a+c+d}(u),x_{-a-b-c-d}(t)]=x_{-b}(-tu).\eqno{(049)}
$$

$$
[x_{a+c+d}(u),x_{-a-c-d-e}(t)]=x_{-e}(-tu).\eqno{(040)}
$$

$$
[x_{a+c+d}(u),x_{-a-b-c-2d-e}(t)]=x_{-b-d-e}(tu).\eqno{(051)}
$$

$$
[x_{a+c+d}(u),x_{-a-b-2c-2d-e}(t)]=x_{-b-c-d-e}(tu).\eqno{(052)}
$$

$$
[x_{a+c+d}(u),x_{-a-c-d-e-f}(t)]=x_{-e-f}(-tu).\eqno{(053)}
$$

$$
[x_{a+c+d}(u),x_{-a-b-c-2d-e-f}(t)]=x_{-b-d-e-f}(tu).\eqno{(054)}
$$

$$
[x_{a+c+d}(u),x_{-a-b-2c-2d-e-f}(t)]=x_{-b-c-d-e-f}(tu).\eqno{(055)}
$$

$$
[x_{a+c+d}(u),x_{-a-b-2c-3d-2e-f}(t)]=x_{-b-c-2d-2e-f}(-tu).\eqno{(056)}
$$

{\boldmath$b:$}

$$
[x_{b}(u),x_{-b-d}(t)]=x_{-d}(tu).\eqno{(057)}
$$

$$
[x_{b}(u),x_{-b-c-d}(t)]=x_{-c-d}(tu).\eqno{(058)}
$$

$$
[x_{b}(u),x_{-a-b-c-d}(t)]=x_{-a-c-d}(tu).\eqno{(059)}
$$

$$
[x_{b}(u),x_{-b-d-e}(t)]=x_{-d-e}(-tu).\eqno{(060)}
$$

$$
[x_{b}(u),x_{-b-c-d-e}(t)]=x_{-c-d-e}(-tu).\eqno{(061)}
$$

$$
[x_{b}(u),x_{-a-b-c-d-e}(t)]=x_{-a-c-d-e}(-tu).\eqno{(062)}
$$

$$
[x_{b}(u),x_{-b-d-e-f}(t)]=x_{-d-e-f}(-tu).\eqno{(063)}
$$

$$
[x_{b}(u),x_{-b-c-d-e-f}(t)]=x_{-c-d-e-f}(-tu).\eqno{(064)}
$$

$$
[x_{b}(u),x_{-a-b-c-d-e-f}(t)]=x_{-a-c-d-e-f}(-tu).\eqno{(065)}
$$

$$
[x_{b}(u),x_{-a-2b-2c-3d-2e-f}(t)]=x_{-a-b-2c-3d-2e-f}(-tu).\eqno{(066)}
$$

{\boldmath$b+d:$}

$$
[x_{b+d}(u),x_{-b-c-d}(t)]=x_{-c}(tu).\eqno{(067)}
$$

$$
[x_{b+d}(u),x_{-a-b-c-d}(t)]=x_{-a-c}(tu).\eqno{(068)}
$$

$$
[x_{b+d}(u),x_{-b-d-e}(t)]=x_{-e}(tu).\eqno{(069)}
$$

$$
[x_{b+d}(u),x_{-b-c-2d-e}(t)]=x_{-c-d-e}(-tu).\eqno{(070)}
$$

$$
[x_{b+d}(u),x_{-a-b-c-2d-e}(t)]=x_{-a-c-d-e}(-tu).\eqno{(071)}
$$

$$
[x_{b+d}(u),x_{-b-d-e-f}(t)]=x_{-e-f}(tu).\eqno{(072)}
$$

$$
[x_{b+d}(u),x_{-b-c-2d-e-f}(t)]=x_{-c-d-e-f}(-tu).\eqno{(073)}
$$

$$
[x_{b+d}(u),x_{-a-b-c-2d-e-f}(t)]=x_{-a-c-d-e-f}(-tu).\eqno{(074)}
$$

$$
[x_{b+d}(u),x_{-a-2b-2c-3d-2e-f}(t)]= x_{-a-b-2c-2d-2e-f}(tu).\eqno{(075)}
$$

{\boldmath$b+c+d:$}

$$
[x_{b+c+d}(u),x_{-a-b-c-d}(t)]=x_{-a}(tu).\eqno{(076)}
$$

$$
[x_{b+c+d}(u),x_{-b-c-d-e}(t)]=x_{-e}(tu).\eqno{(077)}
$$

$$
[x_{b+c+d}(u),x_{-b-c-2d-e}(t)]=x_{-d-e}(tu).\eqno{(078)}
$$

$$
[x_{b+c+d}(u),x_{-a-b-2c-2d-e}(t)]=x_{-a-c-d-e}(-tu).\eqno{(079)}
$$

$$
[x_{b+c+d}(u),x_{-b-c-d-e-f}(t)]=x_{-e-f}(tu).\eqno{(080)}
$$

$$
[x_{b+c+d}(u),x_{-b-c-2d-e-f}(t)]=x_{-d-e-f}(tu).\eqno{(081)}
$$

$$
[x_{b+c+d}(u),x_{-a-b-2c-2d-e-f}(t)]=x_{-a-c-d-e-f}(-tu).\eqno{(082)}
$$

$$
[x_{b+c+d}(u),x_{-a-2b-2c-3d-2e-f}(t)]=x_{-a-b-c-2d-2e-f}(-tu).\eqno{(083)}
$$

{\boldmath$a+b+c+d:$}

$$
[x_{a+b+c+d}(u),x_{-a-b-c-d-e}(t)]=x_{-e}(tu).\eqno{(084)}
$$

$$
[x_{a+b+c+d}(u),x_{-a-b-c-2d-e}(t)]=x_{-d-e}(tu).\eqno{(085)}
$$

$$
[x_{a+b+c+d}(u),x_{-a-b-2c-2d-e}(t)]= x_{-c-d-e}(tu).\eqno{(086)}
$$

$$
[x_{a+b+c+d}(u),x_{-a-b-c-d-e-f}(t)]=x_{-e-f}(tu).\eqno{(087)}
$$

$$
[x_{a+b+c+d}(u),x_{-a-b-c-2d-e-f}(t)]=x_{-d-e-f}(tu).\eqno{(088)}
$$

$$
[x_{a+b+c+d}(u),x_{-a-b-2c-2d-e-f}(t)]= x_{-c-d-e-f}(tu).\eqno{(089)}
$$

$$
[x_{a+b+c+d}(u),x_{-a-2b-2c-3d-2e-f}(t)]= x_{-b-c-2d-2e-f}(tu).\eqno{(090)}
$$

{\boldmath$e:$}

$$
[x_{e}(u),x_{-d-e}(t)]=x_{-d}(tu).\eqno{(091)}
$$

$$
[x_{e}(u),x_{-c-d-e}(t)]=x_{-c-d}(tu).\eqno{(092)}
$$

$$
[x_{e}(u),x_{-a-c-d-e}(t)]=x_{-a-c-d}(tu).\eqno{(093)}
$$

$$
[x_{e}(u),x_{-b-d-e}(t)]=x_{-b-d}(-tu).\eqno{(094)}
$$

$$
[x_{e}(u),x_{-b-c-d-e}(t)]=x_{-b-c-d}(-tu).\eqno{(095)}
$$

$$
[x_{e}(u),x_{-a-b-c-d-e}(t)]=x_{-a-b-c-d}(-tu).\eqno{(096)}
$$

$$
[x_{e}(u),x_{-e-f}(t)]=x_{-f}(-tu).\eqno{(097)}
$$

$$
[x_{e}(u),x_{-b-c-2d-2e-f}(t)]=x_{-b-c-2d-e-f}(-tu).\eqno{(098)}
$$

$$
[x_{e}(u),x_{-a-b-c-2d-2e-f}(t)]= x_{-a-b-c-2d-e-f}(-tu).\eqno{(099)}
$$

$$
[x_{e}(u),x_{-a-b-2c-2d-2e-f}(t)]= x_{-a-b-2c-2d-e-f}(-tu).\eqno{(100)}
$$

{\boldmath$d+e:$}

$$
[x_{d+e}(u),x_{-c-d-e}(t)]=x_{-c}(tu).\eqno{(101)}
$$

$$
[x_{d+e}(u),x_{-a-c-d-e}(t)]=x_{-a-c}(tu).\eqno{(102)}
$$

$$
[x_{d+e}(u),x_{-b-d-e}(t)]=x_{-b}(tu).\eqno{(103)}
$$

$$
[x_{d+e}(u),x_{-b-c-2d-e}(t)]=x_{-b-c-d}(-tu).\eqno{(104)}
$$

$$
[x_{d+e}(u),x_{-a-b-c-2d-e}(t)]=x_{-a-b-c-d}(-tu).\eqno{(105)}
$$

$$
[x_{d+e}(u),x_{-d-e-f}(t)]=x_{-f}(-tu).\eqno{(106)}
$$

$$
[x_{d+e}(u),x_{-b-c-2d-2e-f}(t)]=x_{-b-c-d-e-f}(tu).\eqno{(107)}
$$

$$
[x_{d+e}(u),x_{-a-b-c-2d-2e-f}(t)]= x_{-a-b-c-d-e-f}(tu).\eqno{(108)}
$$

$$
[x_{d+e}(u),x_{-a-b-2c-3d-2e-f}(t)]=x_{-a-b-2c-2d-e-f}(-tu).\eqno{(109)}
$$

{\boldmath$c+d+e:$}

$$
[x_{c+d+e}(u),x_{-a-c-d-e}(t)]=x_{-a}(tu).\eqno{(110)}
$$

$$
[x_{c+d+e}(u),x_{-b-c-d-e}(t)]=x_{-b}(tu).\eqno{(111)}
$$

$$
[x_{c+d+e}(u),x_{-b-c-2d-e}(t)]=x_{-b-d}(tu).\eqno{(112)}
$$

$$
[x_{c+d+e}(u),x_{-a-b-2c-2d-e}(t)]= x_{-a-b-c-d}(-tu).\eqno{(113)}
$$

$$
[x_{c+d+e}(u),x_{-c-d-e-f}(t)]=x_{-f}(-tu).\eqno{(114)}
$$

$$
[x_{c+d+e}(u),x_{-b-c-2d-2e-f}(t)]=x_{-b-d-e-f}(-tu).\eqno{(115)}
$$

$$
[x_{c+d+e}(u),x_{-a-b-2c-2d-2e-f}(t)]= x_{-a-b-c-d-e-f}(tu).\eqno{(116)}
$$

$$
[x_{c+d+e}(u),x_{-a-b-2c-3d-2e-f}(t)]= x_{-a-b-c-2d-e-f}(tu).\eqno{(117)}
$$

{\boldmath$a+c+d+e:$}

$$
[x_{a+c+d+e}(u),x_{-a-b-c-d-e}(t)]=x_{-b}(tu).\eqno{(118)}
$$

$$
[x_{a+c+d+e}(u),x_{-a-b-c-2d-e}(t)]=x_{-b-d}(tu).\eqno{(119)}
$$

$$
[x_{a+c+d+e}(u),x_{-a-b-2c-2d-e}(t)]= x_{-b-c-d}(tu).\eqno{(120)}
$$

$$
[x_{a+c+d+e}(u),x_{-a-c-d-e-f}(t)]=x_{-f}(-tu).\eqno{(121)}
$$

$$
[x_{a+c+d+e}(u),x_{-a-b-c-2d-2e-f}(t)]= x_{-b-d-e-f}(-tu).\eqno{(122)}
$$

$$
[x_{a+c+d+e}(u),x_{-a-b-2c-2d-2e-f}(t)]= x_{-b-c-d-e-f}(-tu).\eqno{(123)}
$$

$$
[x_{a+c+d+e}(u),x_{-a-b-2c-3d-2e-f}(t)]= x_{-b-c-2d-e-f}(-tu).\eqno{(124)}
$$

{\boldmath$b+d+e:$}

$$
[x_{b+d+e}(u),x_{-b-c-d-e}(t)]=x_{-c}(tu).\eqno{(125)}
$$

$$
[x_{b+d+e}(u),x_{-a-b-c-d-e}(t)]=x_{-a-c}(tu).\eqno{(126)}
$$

$$
[x_{b+d+e}(u),x_{-b-c-2d-e}(t)]=x_{-c-d}(-tu).\eqno{(127)}
$$

$$
[x_{b+d+e}(u),x_{-a-b-c-2d-e}(t)]=x_{-a-c-d}(-tu).\eqno{(128)}
$$

$$
[x_{b+d+e}(u),x_{-b-d-e-f}(t)]=x_{-f}(-tu).\eqno{(129)}
$$

$$
[x_{b+d+e}(u),x_{-b-c-2d-2e-f}(t)]= x_{-c-d-e-f}(-tu).\eqno{(130)}
$$

$$
[x_{b+d+e}(u),x_{-a-b-c-2d-2e-f}(t)]= x_{-a-c-d-e-f}(-tu).\eqno{(131)}
$$

$$
[x_{b+d+e}(u),x_{-a-2b-2c-3d-2e-f}(t)]= x_{-a-b-2c-2d-e-f}(-tu).\eqno{(132)}
$$

{\boldmath$b+c+d+e:$}

$$
[x_{b+c+d+e}(u),x_{-a-b-c-d-e}(t)]=x_{-a}(tu).\eqno{(133)}
$$

$$
[x_{b+c+d+e}(u),x_{-b-c-2d-e}(t)]=x_{-d}(tu).\eqno{(134)}
$$

$$
[x_{b+c+d+e}(u),x_{-a-b-2c-2d-e}(t)]=x_{-a-c-d}(-tu).\eqno{(135)}
$$

$$
[x_{b+c+d+e}(u),x_{-b-c-d-e-f}(t)]=x_{-f}(-tu).\eqno{(136)}
$$

$$
[x_{b+c+d+e}(u),x_{-b-c-2d-2e-f}(t)]=x_{-d-e-f}(tu).\eqno{(137)}
$$

$$
[x_{b+c+d+e}(u),x_{-a-b-2c-2d-2e-f}(t)]= x_{-a-c-d-e-f}(-tu).\eqno{(138)}
$$

$$
[x_{b+c+d+e}(u),x_{-a-2b-2c-3d-2e-f}(t)]= x_{-a-b-c-2d-e-f}(tu).\eqno{(139)}
$$

{\boldmath$a+b+c+d+e:$}

$$
[x_{a+b+c+d+e}(u),x_{-a-b-c-2d-e}(t)]=x_{-d}(tu).\eqno{(140)}
$$

$$
[x_{a+b+c+d+e}(u),x_{-a-b-2c-2d-e}(t)]=x_{-c-d}(tu).\eqno{(141)}
$$

$$
[x_{a+b+c+d+e}(u),x_{-a-b-c-d-e-f}(t)]= x_{-f}(-tu).\eqno{(142)}
$$

$$
[x_{a+b+c+d+e}(u),x_{-a-b-c-2d-2e-f}(t)]= x_{-d-e-f}(tu).\eqno{(143)}
$$

$$
[x_{a+b+c+d+e}(u),x_{-a-b-2c-2d-2e-f}(t)]= x_{-c-d-e-f}(tu).\eqno{(144)}
$$

$$
[x_{a+b+c+d+e}(u),x_{-a-2b-2c-3d-2e-f}(t)]= x_{-b-c-2d-e-f}(-tu).\eqno{(145)}
$$

{\boldmath$b+c+2d+e:$}

$$
[x_{b+c+2d+e}(u),x_{-a-b-c-2d-e}(t)]=x_{-a}(tu).\eqno{(146)}
$$

$$
[x_{b+c+2d+e}(u),x_{-a-b-2c-2d-e}(t)]=x_{-a-c}(-tu).\eqno{(147)}
$$

$$
[x_{b+c+2d+e}(u),x_{-b-c-2d-e-f}(t)]=x_{-f}(-tu).\eqno{(148)}
$$

$$
[x_{b+c+2d+e}(u),x_{-b-c-2d-2e-f}(t)]= x_{-e-f}(-tu).\eqno{(149)}
$$

$$
[x_{b+c+2d+e}(u),x_{-a-b-2c-3d-2e-f}(t)]= x_{-a-c-d-e-f}(-tu).\eqno{(150)}
$$

$$
[x_{b+c+2d+e}(u),x_{-a-2b-2c-3d-2e-f}(t)]= x_{-a-b-c-d-e-f}(-tu).\eqno{(151)}
$$

{\boldmath$a+b+c+2d+e:$}

$$
[x_{a+b+c+2d+e}(u),x_{-a-b-2c-2d-e}(t)]=x_{-c}(tu).\eqno{(152)}
$$

$$
[x_{a+b+c+2d+e}(u),x_{-a-b-c-2d-e-f}(t)]= x_{-f}(-tu).\eqno{(153)}
$$

$$
[x_{a+b+c+2d+e}(u),x_{-a-b-c-2d-2e-f}(t)]= x_{-e-f}(-tu).\eqno{(154)}
$$

$$
[x_{a+b+c+2d+e}(u),x_{-a-b-2c-3d-2e-f}(t)]= x_{-c-d-e-f}(tu).\eqno{(155)}
$$

$$
[x_{a+b+c+2d+e}(u),x_{-a-2b-2c-3d-2e-f}(t)]= x_{-b-c-d-e-f}(tu).\eqno{(156)}
$$

{\boldmath$a+b+2c+2d+e:$}

$$
[x_{a+b+2c+2d+e}(u),x_{-a-b-2c-2d-e-f}(t)]= x_{-f}(-tu).\eqno{(157)}
$$

$$
[x_{a+b+2c+2d+e}(u),x_{-a-b-2c-2d-2e-f}(t)]= x_{-e-f}(-tu).\eqno{(158)}
$$

$$
[x_{a+b+2c+2d+e}(u),x_{-a-b-2c-3d-2e-f}(t)]= x_{-d-e-f}(-tu).\eqno{(159)}
$$

$$
[x_{a+b+2c+2d+e}(u),x_{-a-2b-2c-3d-2e-f}(t)]= x_{-b-d-e-f}(-tu).\eqno{(160)}
$$

{\boldmath$f:$}

$$
[x_{f}(u),x_{-e-f}(t)]=x_{-e}(tu).\eqno{(161)}
$$

$$
[x_{f}(u),x_{-d-e-f}(t)]=x_{-d-e}(tu).\eqno{(162)}
$$

$$
[x_{f}(u),x_{-c-d-e-f}(t)]=x_{-c-d-e}(tu).\eqno{(163)}
$$

$$
[x_{f}(u),x_{-b-d-e-f}(t)]=x_{-b-d-e}(tu).\eqno{(164)}
$$

$$
[x_{f}(u),x_{-a-c-d-e-f}(t)]=x_{-a-c-d-e}(tu).\eqno{(165)}
$$

$$
[x_{f}(u),x_{-b-c-d-e-f}(t)]=x_{-b-c-d-e}(tu).\eqno{(166)}
$$

$$
[x_{f}(u),x_{-a-b-c-d-e-f}(t)]=x_{-a-b-c-d-e}(tu).\eqno{(167)}
$$

$$
[x_{f}(u),x_{-b-c-2d-e-f}(t)]=x_{-b-c-2d-e}(tu).\eqno{(168)}
$$

$$
[x_{f}(u),x_{-a-b-c-2d-e-f}(t)]= x_{-a-b-c-2d-e}(tu).\eqno{(169)}
$$

$$
[x_{f}(u),x_{-a-b-2c-2d-e-f}(t)]= x_{-a-b-2c-2d-e}(tu).\eqno{(170)}
$$

{\boldmath$e+f:$}

$$
[x_{e+f}(u),x_{-d-e-f}(t)]=x_{-d}(tu).\eqno{(171)}
$$

$$
[x_{e+f}(u),x_{-c-d-e-f}(t)]=x_{-c-d}(tu).\eqno{(172)}
$$

$$
[x_{e+f}(u),x_{-a-c-d-e-f}(t)]=x_{-a-c-d}(tu).\eqno{(173)}
$$

$$
[x_{e+f}(u),x_{-b-d-e-f}(t)]=x_{-b-d}(-tu).\eqno{(174)}
$$

$$
[x_{e+f}(u),x_{-b-c-d-e-f}(t)]=x_{-b-c-d}(-tu).\eqno{(175)}
$$

$$
[x_{e+f}(u),x_{-a-b-c-d-e-f}(t)]=x_{-a-b-c-d}(-tu).\eqno{(176)}
$$

$$
[x_{e+f}(u),x_{-b-c-2d-2e-f}(t)]=x_{-b-c-2d-e}(tu).\eqno{(177)}
$$

$$
[x_{e+f}(u),x_{-a-b-c-2d-2e-f}(t)]= x_{-a-b-c-2d-e}(tu).\eqno{(178)}
$$

$$
[x_{e+f}(u),x_{-a-b-2c-2d-2e-f}(t)]= x_{-a-b-2c-2d-e}(tu).\eqno{(179)}
$$

{\boldmath$d+e+f:$}

$$
[x_{d+e+f}(u),x_{-c-d-e-f}(t)]=x_{-c}(tu).\eqno{(180)}
$$

$$
[x_{d+e+f}(u),x_{-a-c-d-e-f}(t)]=x_{-a-c}(tu).\eqno{(181)}
$$

$$
[x_{d+e+f}(u),x_{-b-d-e-f}(t)]=x_{-b}(tu).\eqno{(182)}
$$

$$
[x_{d+e+f}(u),x_{-b-c-2d-e-f}(t)]=x_{-b-c-d}(-tu).\eqno{(183)}
$$

$$
[x_{d+e+f}(u),x_{-a-b-c-2d-e-f}(t)]= x_{-a-b-c-d}(-tu).\eqno{(184)}
$$

$$
[x_{d+e+f}(u),x_{-b-c-2d-2e-f}(t)]= x_{-b-c-d-e}(-tu).\eqno{(185)}
$$

$$
[x_{d+e+f}(u),x_{-a-b-c-2d-2e-f}(t)]= x_{-a-b-c-d-e}(-tu).\eqno{(186)}
$$

$$
[x_{d+e+f}(u),x_{-a-b-2c-3d-2e-f}(t)]= x_{-a-b-2c-2d-e}(tu).\eqno{(187)}
$$

{\boldmath$c+d+e+f:$}

$$
[x_{c+d+e+f}(u),x_{-a-c-d-e-f}(t)]=x_{-a}(tu).\eqno{(188)}
$$

$$
[x_{c+d+e+f}(u),x_{-b-c-d-e-f}(t)]=x_{-b}(tu).\eqno{(189)}
$$

$$
[x_{c+d+e+f}(u),x_{-b-c-2d-e-f}(t)]=x_{-b-d}(tu).\eqno{(190)}
$$

$$
[x_{c+d+e+f}(u),x_{-b-c-2d-2e-f}(t)]=x_{-b-d-e}(tu).\eqno{(191)}
$$

$$
[x_{c+d+e+f}(u),x_{-a-b-2c-2d-e-f}(t)]= x_{-a-b-d-e}(-tu).\eqno{(192)}
$$

$$
[x_{c+d+e+f}(u),x_{-a-b-2c-2d-2e-f}(t)]= x_{-a-b-c-d-e}(-tu).\eqno{(193)}
$$

$$
[x_{c+d+e+f}(u),x_{-a-b-2c-3d-2e-f}(t)]= x_{-a-b-c-2d-e}(-tu).\eqno{(194)}
$$

{\boldmath$a+c+d+e+f:$}

$$
[x_{a+c+d+e+f}(u),x_{-a-b-c-d-e-f}(t)]=x_{-b}(tu).\eqno{(195)}
$$

$$
[x_{a+c+d+e+f}(u),x_{-a-b-c-2d-e-f}(t)]=x_{-b-d}(tu).\eqno{(196)}
$$

$$
[x_{a+c+d+e+f}(u),x_{-a-b-2c-2d-e-f}(t)]= x_{-b-c-d}(tu).\eqno{(197)}
$$

$$
[x_{a+c+d+e+f}(u),x_{-a-b-c-2d-2e-f}(t)]= x_{-b-d-e}(tu).\eqno{(198)}
$$

$$
[x_{a+c+d+e+f}(u),x_{-a-b-2c-2d-2e-f}(t)]= x_{-b-c-d-e}(tu).\eqno{(199)}
$$

$$
[x_{a+c+d+e+f}(u),x_{-a-b-2c-3d-2e-f}(t)] = x_{-b-c-2d-e}(tu).\eqno{(200)}
$$

{\boldmath$b+d+e+f:$}

$$
[x_{b+d+e+f}(u),x_{-b-c-d-e-f}(t)]=x_{-c}(tu).\eqno{(201)}
$$

$$
[x_{b+d+e+f}(u),x_{-a-b-c-d-e-f}(t)]=x_{-a-c}(tu).\eqno{(202)}
$$

$$
[x_{b+d+e+f}(u),x_{-b-c-2d-e-f}(t)]=x_{-c-d}(-tu).\eqno{(203)}
$$

$$
[x_{b+d+e+f}(u),x_{-a-b-c-2d-e-f}(t)]= x_{-a-c-d}(-tu).\eqno{(204)}
$$

$$
[x_{b+d+e+f}(u),x_{-b-c-2d-2e-f}(t)]=x_{-c-d-e}(tu).\eqno{(205)}
$$

$$
[x_{b+d+e+f}(u),x_{-a-b-c-2d-2e-f}(t)]= x_{-a-c-d-e}(tu).\eqno{(206)}
$$

$$
[x_{b+d+e+f}(u),x_{-a-2b-2c-3d-2e-f}(t)]= x_{-a-b-2c-2d-e}(tu).\eqno{(207)}
$$

{\boldmath$b+c+d+e+f:$}

$$
[x_{b+c+d+e+f}(u),x_{-a-b-c-d-e-f}(t)]=x_{-a}(tu).\eqno{(208)}
$$

$$
[x_{b+c+d+e+f}(u),x_{-b-c-2d-e-f}(t)]=x_{-d}(tu).\eqno{(209)}
$$

$$
[x_{b+c+d+e+f}(u),x_{-b-c-2d-2e-f}(t)]=x_{-d-e}(-tu).\eqno{(210)}
$$

$$
[x_{b+c+d+e+f}(u),x_{-a-b-2c-2d-e-f}(t)]= x_{-a-c-d}(-tu).\eqno{(211)}
$$

$$
[x_{b+c+d+e+f}(u),x_{-a-b-2c-2d-2e-f}(t)]= x_{-a-c-d-e}(tu).\eqno{(212)}
$$

$$
[x_{b+c+d+e+f}(u),x_{-a-2b-2c-3d-2e-f}(t)]= x_{-a-b-c-2d-e}(-tu).\eqno{(213)}
$$

{\boldmath$a+b+c+d+e+f:$}

$$
[x_{a+b+c+d+e+f}(u),x_{-a-b-c-2d-e-f}(t)]=x_{-d}(tu).\eqno{(214)}
$$

$$
[x_{a+b+c+d+e+f}(u),x_{-a-b-2c-2d-e-f}(t)]= x_{-c-d}(tu).\eqno{(215)}
$$

$$
[x_{a+b+c+d+e+f}(u),x_{-a-b-c-2d-2e-f}(t)]= x_{-d-e}(-tu).\eqno{(216)}
$$

$$
[x_{a+b+c+d+e+f}(u),x_{-a-b-2c-2d-2e-f}(t)]= x_{-c-d-e}(-tu).\eqno{(217)}
$$

$$
[x_{a+b+c+d+e+f}(u),x_{-a-2b-2c-3d-2e-f}(t)]= x_{-b-c-2d-e}(tu).\eqno{(218)}
$$

{\boldmath$b+c+2d+e+f:$}

$$
[x_{b+c+2d+e+f}(u),x_{-a-b-c-2d-e-f}(t)]=x_{-a}(tu).\eqno{(219)}
$$

$$
[x_{b+c+2d+e+f}(u),x_{-b-c-2d-2e-f}(t)]=x_{-e}(tu).\eqno{(220)}
$$

$$
[x_{b+c+2d+e+f}(u),x_{-a-b-2c-2d-e-f}(t)]= x_{-a-c}(-tu).\eqno{(221)}
$$

$$
[x_{b+c+2d+e+f}(u),x_{-a-b-2c-3d-2e-f}(t)]= x_{-a-c-d-e}(tu).\eqno{(222)}
$$

$$
[x_{b+c+2d+e+f}(u),x_{-a-2b-2c-3d-2e-f}(t)]= x_{-a-b-c-d-e}(tu).\eqno{(223)}
$$

{\boldmath$a+b+c+2d+e+f:$}

$$
[x_{a+b+c+2d+e+f}(u),x_{-a-b-2c-2d-e-f}(t)]=x_{-c}(tu).\eqno{(224)}
$$

$$
[x_{a+b+c+2d+e+f}(u),x_{-a-b-c-2d-2e-f}(t)]=x_{-e}(tu).\eqno{(225)}
$$

$$
[x_{a+b+c+2d+e+f}(u),x_{-a-b-2c-3d-2e-f}(t)]= x_{-c-d-e}(-tu).\eqno{(226)}
$$

$$
[x_{a+b+c+2d+e+f}(u),x_{-a-2b-2c-3d-2e-f}(t)]= x_{-b-c-d-e}(-tu).\eqno{(227)}
$$

{\boldmath$a+b+2c+2d+e+f:$}

$$
[x_{a+b+2c+2d+e+f}(u),x_{-a-b-2c-2d-2e-f}(t)]= x_{-e}(tu).\eqno{(228)}
$$

$$
[x_{a+b+2c+2d+e+f}(u),x_{-a-b-2c-3d-2e-f}(t)]= x_{-d-e}(tu).\eqno{(229)}
$$

$$
[x_{a+b+2c+2d+e+f}(u),x_{-a-2b-2c-3d-2e-f}(t)]= x_{-b-d-e}(tu).\eqno{(230)}
$$

{\boldmath$b+c+2d+2e+f:$}

$$
[x_{b+c+2d+2e+f}(u),x_{-a-b-c-2d-2e-f}(t)]=x_{-a}(tu).\eqno{(231)}
$$

$$
[x_{b+c+2d+2e+f}(u),x_{-a-b-2c-2d-2e-f}(t)]= x_{-a-c}(-tu).\eqno{(232)}
$$

$$
[x_{b+c+2d+2e+f}(u),x_{-a-b-2c-3d-2e-f}(t)]= x_{-a-c-d}(tu).\eqno{(233)}
$$

$$
[x_{b+c+2d+2e+f}(u),x_{-a-2b-2c-3d-2e-f}(t)]= x_{-a-b-c-d}(-tu).\eqno{(234)}
$$

{\boldmath$a+b+c+2d+2e+f:$}

$$
[x_{a+b+c+2d+2e+f}(u),x_{-a-b-2c-2d-2e-f}(t)]= x_{-c}(tu).\eqno{(235)}
$$

$$
[x_{a+b+c+2d+2e+f}(u),x_{-a-b-2c-3d-2e-f}(t)]= x_{-c-d}(-tu).\eqno{(236)}
$$

$$
[x_{a+b+c+2d+2e+f}(u),x_{-a-2b-2c-3d-2e-f}(t)]= x_{-b-c-d}(tu).\eqno{(237)}
$$

{\boldmath$a+b+2c+2d+2e+f:$}

$$
[x_{a+b+2c+2d+2e+f}(u),x_{-a-b-2c-3d-2e-f}(t)]= x_{-d}(tu).\eqno{(238)}
$$

$$
[x_{a+b+2c+2d+2e+f}(u),x_{-a-2b-2c-3d-2e-f}(t)]= x_{-b-d}(-tu).\eqno{(239)}
$$

{\boldmath$a+b+2c+3d+2e+f:$}

$$
[x_{a+b+2c+3d+2e+f}(u),x_{-a-2b-2c-3d-2e-f}(t)]= x_{-b}(tu).\eqno{(240)}
$$

\section{Graphs, Matrices, Tables}

Below we present the directed weighted graph $G(E_6^-,\Delta)$ for $\Delta=\{ a,b,c,d,e,f\},$ its adjacency matrix 
and matrix of weights, using them  we calculate the table of the number of paths and the table of numbers 
$K_{r,s}^{\Delta}.$ Information about those concepts 
can be found in [3].
\bigskip

\begin{center}
\thicklines

\end{center}

\newpage

Let $\Delta = \{a,b,c,d,e,f\}.$ 
\bigskip

\centerline{\textbf{Adjacency Matrix of the Graph $G(\Phi^-,\Delta)$}}

\vskip5pt

{\scriptsize
$$
\left(%
\right)
$$
}

\newpage

\centerline{\textbf{Table of the number of paths from root $s$ to root $r$}}
\vskip15pt

\begin{center}
%
\vskip4mm
\ \ \ \ \ \ \ \ \ \ Table 7, part 3 of 3.

\newpage
\centerline{\textbf{Weighted Graph Adjacency Matrix $G(\Phi^-,\Delta)$}}

{\scriptsize
$$
\left (%
\right )
$$
}

\newpage
\centerline{\textbf{The Table of the Number $K_{r,s}^\Delta$ $r,s\in\Phi^-$}}
\vskip20pt

\begin{center}
%
\vskip4mm

\ \ \ \ \ \ \ \ \ \ \ Table 8, part 3 of 3.

\newpage

\centerline{\textbf{The Table of the Number $K_{r,s}^\Delta$ $r,s\in\Phi^+$}}
\vskip20pt

\begin{center}
%
\vskip4mm

\ \ \ \ \ \ \ \ Table 9, part 3 of 3.

\newpage

Let $\Delta = \{b,c,d,e,f\}.$
\bigskip

\centerline{\textbf{Adjacency Matrix of the Graph $G(\Phi^-,\Delta)$}}
\vskip20pt

{\scriptsize
$$
\left(%
\right)
$$ 
}

\newpage

\centerline{\textbf{Table of the number of paths from root $s$ to root $r$}}
\vskip20pt

\begin{center}
%
\vskip4mm

\ \ \ \ \ \ \ \ \ \ \ \ \ \ \ \ \ \ \ \ \ \ \ \ \ \ \ \ \ \ \ \ \ \ \ \ \ \ \ Table 10, part 2 of 2.

\newpage

\centerline{\textbf{Weighted Graph Adjacency Matrix $G(\Phi^-,\Delta)$}}
\vskip20pt

{\scriptsize
$$
\left [%
\right ]
$$
}

\newpage

\centerline{\textbf{The Table of the Number $K_{r,s}^\Delta,$ $r,s\in\Phi^-$}}
\vskip2mm

%
\vskip4mm

\ \ \ \ \ \ \ \ \ \ \ \ \ \ \ \ \ \ \ \ \ \ \ \ \ \ \ \ \ \ \ \ \ \ \ \ \ \ \ \ \ \ \ Table 11, part 2 of 2.

\newpage

\centerline{\textbf{The Table of the Number $K_{r,s}^\Delta,$ $r,s\in\Phi^+$}}

%
\vskip4mm

\ \ \ \ \ \ \ \ \ \ \ \  \ \ \ \ \ \ \ \ \ \ \ \ \ \ \ \ \ \ \ \ \ \ \ \ \ \ \ \ \ \ Table 12, part 2 of 2.

\bigskip

Anna I. Polovinkina

Siberian Federal University

email: aipolovinkina@sfu-kras.ru

\medskip

Sergey G. Kolesnikov

Siberian Federal University

email: skolesnikov@sfu-kras.ru


\begin{thebibliography}{1}

\bibitem{Bur72}
{\sl Bourbaki N.} Lie groups and algebras. M.: Mir, 1972. 334 p.

\bibitem{Car72}
{\sl Carter R.} Simple groups of Lie type.-Ney York: Wiley and
Sons, 1972. 458 p.

\bibitem{EKL2023}
{\sl Kolesnikov S.G., Polovinkina A.I.} https://arxiv.org/abs/2312.03439

\end{thebibliography}
\end{document}